\documentclass[10pt]{article}

\usepackage{amssymb,amsmath,amsthm,amsbsy}
\usepackage[title]{appendix}
\usepackage[affil-it]{authblk}
\usepackage[font=footnotesize,labelfont=bf]{caption}
\usepackage[scale=.78]{geometry}
\usepackage{graphicx}
\usepackage{hyperref}
\usepackage{microtype}
\usepackage{physics}
\usepackage{tcolorbox}
\usepackage{url}
\usepackage[nameinlink]{cleveref}

\hypersetup{
    colorlinks=true,
    linkcolor=black,
    citecolor=black,
    urlcolor=black,
    pdftitle={Bayesian learning with Gaussian processes for low-dimensional representations of time-dependent nonlinear systems},
    pdfauthor={Shane A. McQuarrie, Anirban Chaudhuri, Karen E. Willcox, Mengwu Guo}
}

\newcommand{\trp}{^{\mathsf{T}}}
\crefformat{equation}{#2(#1)#3}
\numberwithin{equation}{section}

\usepackage[section]{algorithm}
\usepackage[noend]{algpseudocode}
\algnewcommand{\LineComment}[1]{\State \textcolor{gray}{\texttt{\# #1}}}
\algrenewcomment[1]{\hfill\textcolor{gray}{\texttt{\# #1}}}

\theoremstyle{definition}
\newtheorem{remark}{Remark}[section]

\title{\Large\textbf{Bayesian learning with Gaussian processes for low-dimensional representations of time-dependent nonlinear systems}}

\usepackage{fancyhdr}
\newcommand{\apdxref}[1]{\hyperref[#1]{Appendix~\ref{#1}}}

\title{\Large\textbf{Bayesian learning with Gaussian processes for low-dimensional representations of time-dependent nonlinear systems}}
\author[1]{Shane~A.~McQuarrie}
\affil[1]{\normalsize Scientific Machine Learning, Sandia National Laboratories}

\author[2]{Anirban Chaudhuri}
\affil[2]{\normalsize Oden Institute for Computational Engineering and Sciences, The University of Texas at Austin}

\author[2]{Karen~E.~Willcox}

\author[3]{Mengwu Guo\thanks{Corresponding author. E-mail: \href{mailto:mengwu.guo@math.lu.se}{\texttt{mengwu.guo@math.lu.se}}.}}
\affil[3]{\normalsize Centre for Mathematical Sciences, Lund University}

\date{}

\begin{document} 

\thispagestyle{empty}  

\maketitle


\begin{abstract}
\noindent
This work presents a data-driven method for learning low-dimensional time-dependent physics-based surrogate models whose predictions are endowed with uncertainty estimates. We use the operator inference approach to model reduction that poses the problem of learning low-dimensional model terms as a regression of state space data and corresponding time derivatives by minimizing the residual of reduced system equations. Standard operator inference models perform well with accurate training data that are dense in time, but producing stable and accurate models when the state data are noisy and/or sparse in time remains a challenge. Another challenge is the lack of uncertainty estimation for the predictions from the operator inference models. Our approach addresses these challenges by incorporating Gaussian process surrogates into the operator inference framework to~(1)~probabilistically describe uncertainties in the state predictions and~(2)~procure analytical time derivative estimates with quantified uncertainties. The formulation leads to a generalized least-squares regression and, ultimately, reduced-order models that are described probabilistically with a closed-form expression for the posterior distribution of the operators. The resulting probabilistic surrogate model propagates uncertainties from the observed state data to reduced-order predictions. We demonstrate the method is effective for constructing low-dimensional models of two nonlinear partial differential equations representing a compressible flow and a nonlinear diffusion-reaction process, as well as for estimating the parameters of a low-dimensional system of nonlinear ordinary differential equations representing compartmental models in epidemiology.
\end{abstract}

\vspace{2mm}

\section{Introduction}

Engineering design and analysis often involves many-query tasks, such as sensitivity analysis, uncertainty quantification, parameter estimation, and optimization, that require simulating a complex, expensive-to-simulate nonlinear system. To enable computationally efficient many-query simulations, we need data-driven surrogate models that can embed known physics to enable predictions outside of the training regime with respect to time and/or model inputs and build trust in such approximations by understanding the uncertainty in the predictions. This paper focuses on physics-based model reduction, the process of constructing low-dimensional systems as approximation models for high-dimensional processes, in the context of complex time-dependent dynamical systems. However, data-driven model reduction does not typically provide an assessment of uncertainty in the predictions. We develop an interpretable, physics-aware probabilistic surrogate that addresses the challenge of characterizing prediction uncertainty in data-driven reduced-order modeling for nonlinear dynamical systems with sparse and/or noisy data.

Model reduction strategies confine system dynamics to a low-dimensional linear subspace or nonlinear manifold and define evolution equations for the coordinates of the original state in this reduced space \cite{BGW2015pmorSurvery,quarteroni2014reduced}. Classical reduced-order models (ROMs) are \emph{intrusive} in the sense that a numerical solver for the system of interest, called the full-order model (FOM), is used explicitly in the construction of the ROM dynamics. In terms of implementation, intrusive ROMs are constructed by acting directly on code objects from the FOM. By contrast, \emph{non-intrusive} or \emph{data-driven} ROMs are constructed from system data alone, without needing access to any codes for solving the system of interest. This is a major advantage in scenarios without a FOM (e.g., experimental data), or in which the FOM can only be queried as a black box (e.g., proprietary, classified, or legacy codes). Furthermore, non-intrusive ROMs can learn from data that has been preprocessed with variable transformations or other advantageous alterations, whereas such preprocessing may severely complicate or even preclude intrusive ROM construction~\cite{QKPW2020liftAndLearn,SKHW2020romCombustion}. Although non-intrusive ROMs can be constructed in many scenarios where intrusive ROMs are unavailable, they usually lack the accuracy guarantees and/or uncertainty estimates that many intrusive ROMs enjoy, such as in \cite{hesthaven2016certified}. Data-driven ROMs, especially those of low intrinsic dimensions, do not usually come with rigorous accuracy guarantees and must therefore be equipped with probabilistic estimates to characterize uncertainty in the predictions, especially in the presence of sparse and/or noisy data. We focus on operator inference (OpInf) as the underlying ROM methodology in our work and build a probabilistic ROM.

Operator inference is a regression-based method for non-intrusively constructing interpretable ROMs from state data \cite{ghattas2021acta,kramer2024opinfsurvey,PW2016operatorInference}. A polynomial model form is chosen based on the governing dynamics (or a variable transformation of the dynamics \cite{QKPW2020liftAndLearn}), thereby embedding the physics structure of the problem into the ROM. The regression requires state data and their time derivatives. These time derivatives may be provided by a FOM, but more often they must be estimated from the state data, typically via finite differences. In the presence of sparse and/or noisy state data, the time derivative estimation and the OpInf regression problem in general suffer from accuracy issues. Uncertainty due to noise in the states, inaccuracies in the time derivative estimates, or model form error can be partially accounted for by introducing a regularization to the regression~\cite{mcquarrie2021combustion,mcquarrie2023popinf,qian2022reduced,sawant2023piregopinf} or by enforcing physics-motivated constraints~\cite{geng2024gradientopinf,gruber2023hamiltonianopinf,koike2024energyopinf,sharma2024lagrangian,sharma2022hamiltonianopinf}. Dealing with data noise remains a challenge, but some progress has been made in recent years. In \cite{uy2023activeopinf}, the re-projection scheme of \cite{Peherstorfer2020reprojection} is coupled with an active-learning paradigm for strategically querying a FOM with noisy states. The procedure partially filters the training data, produces a data-driven ROM that can be related to an intrusive ROM, and decreases the estimated bias in the predicted states by at least one order of magnitude.
In \cite{uy2023opinfrollout}, a nonconvex output-based optimization is used instead of a residual-based regression to learn a ROM with polynomial structure. The method is computationally expensive compared with traditional OpInf but was shown to perform well in several scenarios with sparse data and significant noise pollution. However, OpInf ROMs do not come with rigorous characterization of uncertainty in the predictions. Our previous work \cite{guo2022bayesopinf} interprets OpInf as a Bayesian linear regression and develops a framework for producing probabilistic OpInf ROMs. The method poses a Gaussian noise model for the derivative of the state but still relies on high-quality time derivative data. In this paper, we utilize Gaussian processes to pose a noise model for the state, not its time derivative, and simultaneously obtain time derivative data in an analytical fashion along with a Bayesian version of OpInf.

Gaussian process (GP) regression~\cite{rasmussen2005GPsforMLbook} is a powerful tool for function approximation from data. The nonparametric Bayesian nature of GPs facilitates uncertainty quantification in the predictions, while its kernel-method essence ensures robustness in dealing with noisy and/or scarce data. Notwithstanding, GPs suffer from poor performance in temporal predictions beyond training coverage due to a general lack of physical consistency, as GP~surrogate modeling is conducted in a black-box manner and no information about the governing equations is incorporated. Several recent strategies aim to overcome this challenge by synthesizing GP regression and physics-based models~\cite{chen2021pdegps,hansen2024learningconservation,kast2020non,raissi2017GPsforlinearPDEs,zhang2019model}. In data-driven model reduction specifically, GP~regression has been used to approximate the reduced trajectory over the time-parameter domain of interest~\cite{beckers2023pegpvae,cicci2023uncertainty,guo2019data,guo2021learning}. In this work, we combine GPs with OpInf to develop physics-aware probabilistic surrogates.

We develop a GP-based Bayesian OpInf method (GP-BayesOpInf) for predictive modeling of nonlinear dynamical systems that produces probabilistic predictions. Specifically, GPs are used for the likelihood definition in a Bayesian scheme of OpInf, combining the advantage of GP emulation in smoothing and interpolating imperfect data, and the predictive power of OpInf guaranteed by the preservation of formulation structure of the governing equations in projection-based model reduction. Based on the fact that Gaussianity is preserved through linear operations, differential-equation constraints of reduced states are embedded in the likelihood function. This does not require direct time-derivative approximations as in a Bayesian linear regression for equation fitting~\cite{guo2022bayesopinf} and thus avoids generating low-quality derivative data when training trajectories are noisy and/or scarce. Thereafter, the resulting posterior estimate of the reduced operators can be formulated as a generalized least squares, in which an improved time derivative approximation is enabled by kernel regression and endowed with quantifiable uncertainty. We derive a closed form for the posterior distribution of the operators in OpInf, thus providing an interpretable, low-dimensional probabilistic surrogate model. The resulting physics-aware probabilistic surrogate can then sample realizations of operators from this posterior and integrate the low-dimensional dynamical system for each realization to  provide predictions of the states of the dynamical system along with uncertainty in the predictions stemming from training data and modeling process.

The remainder of the paper is structured as follows.
\Cref{sec:background} provides a technical overview of the OpInf procedure for learning ROMs from state data.
\Cref{sec:opinf-bayes} develops a procedure that combines GP~regression, Bayesian inference, and OpInf to learn probabilistic ROMs. Implementation details and a complete algorithm are presented in \Cref{sec:procedure}. GP-BayesOpInf is demonstrated for two nonlinear partial differential equations~(PDEs) in \Cref{sec:results}. \Cref{sec:odeparamestimate} shows that an adjustment to GP-BayesOpInf results in a parameter estimation strategy for a large class of systems of ordinary differential equations (ODEs). Finally, \Cref{sec:conclusion} concludes the paper.
See \apdxref{appendix:nomenclature} for a list of the major mathematical notation used throughout.

\section{Reduced-order models via operator inference}
\label{sec:background}

Consider a high-dimensional system of ODEs,
\begin{align}
    \label{eq:fullorder}
    \frac{\textup{d}}{\textup{d}t}{\vb{q}}(t)
    = {\vb{f}}({\vb{q}}(t),{\vb{u}}(t))\,,
    \quad
    \vb{q}(t_0) = \vb{q}_0\,,
    \quad
    t\in [t_0, t_f]\,,
\end{align}
where $\vb{q}(t)\in\mathbb{R}^{N}$ is the state of the system at time~$t$,
$\vb{u}(t)\in\mathbb{R}^{p}$ represents inputs at time~$t$ (e.g., time-dependent source terms and/or boundary conditions),
$\vb{f}:\mathbb{R}^{N}\times\mathbb{R}^{p}\to\mathbb{R}^{N}$ specifies the evolutionary system dynamics,
$t_0\in\mathbb{R}$ is the initial time,
$\vb{q}_{0}\in\mathbb{R}^{N}$ is the initial condition,
and $t_f > t_0$ is the final time for which we want to know the state $\vb{q}(t)$.
High-dimensional ODE systems with large $N$ arise often from the spatial discretization of PDEs.
A ROM of \cref{eq:fullorder} is a system of ODEs with substantially fewer degrees of freedom whose solution---computed much more efficiently than directly solving \cref{eq:fullorder}---can be used to approximate the solution to \cref{eq:fullorder}.
The cost reduction resulting from using a ROM is crucial for applications requiring multiple solves of \cref{eq:fullorder}, e.g., with different input functions or initial conditions.

Projection-based ROMs introduce an affine approximation for the high-dimensional state,
\begin{align}
    \label{eq:state-approx}
    \vb{q}(t)
    \approx \vb{V}\hat{\vb{q}}(t) + \bar{\vb{q}},
\end{align}
where $\vb{V}\in\mathbb{R}^{N\times r}$ is the basis matrix, $\hat{\vb{q}}(t)\in\mathbb{R}^{r}$ is the reduced state at time $t$, $\bar{\vb{q}}\in\mathbb{R}^{N}$ is
a fixed reference vector, and $r\ll N$ is the dimension of the reduced state.
Both $\vb{V}$ and $\bar{\vb{q}}$, and hence the reduced state dimension $r$, are either inferred from training observations of the state $\vb{q}(t)$ or chosen based on problem-specific knowledge (e.g., $\bar{\vb{q}}$ may encode constant boundary conditions or represent a time-averaged solution).
We typically require $\vb{V}$ to have orthonormal columns, in which case dynamics for $\hat{\vb{q}}(t)$ are derived by substituting~\cref{eq:state-approx} into \cref{eq:fullorder} and multiplying both sides by $\vb{V}\trp$, resulting in the system
\begin{align}
    \label{eq:intrusive-reduced}
    \frac{\textrm{d}}{\textrm{d}t}\hat{\vb{q}}(t)
    = \hat{\vb{f}}(\hat{\vb{q}}(t),\vb{u}(t))
    := \vb{V}\trp\vb{f}(\vb{V}\hat{\vb{q}}(t) + \bar{\vb{q}}, \vb{u}(t))\,,
    \quad
    \hat{\vb{q}}(t_0) = \vb{V}\trp(\vb{q}_0 - \bar{\vb{q}})\,,
    \quad
    t\in [t_0, t_f]\,,
\end{align}
where $\hat{\vb{f}}:\mathbb{R}^{r}\times\mathbb{R}^{p}\to\mathbb{R}^{r}$.
For a wide class of problems, this approach preserves the dynamical structure of~\cref{eq:fullorder} in the ROM~\cref{eq:intrusive-reduced} \cite{BGW2015pmorSurvery}.
In particular, if $\vb{f}$ has a polynomial form, then so does $\hat{\vb{f}}$:
for example, defining $\vb{c}\in\mathbb{R}^{N}$,
$\vb{A}\in\mathbb{R}^{N\times N}$,
$\vb{H}\in\mathbb{R}^{N\times N^{2}}$, and
$\vb{B}\in\mathbb{R}^{N\times p}$,
writing
\begin{subequations}
\begin{align}
    \label{eq:polynomial-fullorder}
    \frac{\textup{d}}{\textup{d}t}\vb{q}(t) =
    \vb{f}(\vb{q}(t),\vb{u}(t))
    = \vb{c} + \vb{A}\vb{q}(t) + \vb{H}[\vb{q}(t)\otimes\vb{q}(t)] + \vb{B}\vb{u}(t),
\end{align}
there exist corresponding $\hat{\vb{c}}\in\mathbb{R}^{r}$,
$\hat{\vb{A}}\in\mathbb{R}^{r\times r}$,
$\hat{\vb{H}}\in\mathbb{R}^{r\times r^{2}}$, and
$\hat{\vb{B}}\in\mathbb{R}^{r\times p}$
such that
\begin{align}
    \label{eq:polynomial-reducedorder}
    \frac{\textup{d}}{\textup{d}t}\hat{\vb{q}}(t) =
    \hat{\vb{f}}(\hat{\vb{q}}(t),\vb{u}(t))
    = \hat{\vb{c}} + \hat{\vb{A}}\hat{\vb{q}}(t) + \hat{\vb{H}}[\hat{\vb{q}}(t)\otimes\hat{\vb{q}}(t)] + \hat{\vb{B}}\vb{u}(t).
\end{align}
\end{subequations}
In this work, we consider the scenario in which reduced-order operators such as $\hat{\vb{c}}$, $\hat{\vb{A}}$, $\hat{\vb{H}}$, and $\hat{\vb{B}}$ cannot be computed explicitly, as occurs when working with large, proprietary, or classified code bases, or in systems that must be transformed through a change of variables to first expose a polynomial structure similar to \cref{eq:polynomial-fullorder} \cite{QKPW2020liftAndLearn}.
We will see examples of the latter in \Cref{sec:results}.
In these cases, the polynomial structure of $\vb{f}$ is known, but $\vb{f}$ itself is computationally inaccessible and cannot be used directly to form~$\hat{\vb{f}}$.

Operator inference (OpInf) is a non-intrusive projection-based method for learning polynomial ROMs by minimizing the model residual with respect to a collection of training data \cite{kramer2024opinfsurvey,PW2016operatorInference}.
In general, OpInf constructs ROMs that can be written as
\begin{equation}
    \label{eq:reducedorder}
    \frac{\textup{d}}{\textup{d}t} \hat{\vb{q}}(t)
    = \hat{\vb{O}}\vb{d}(\hat{\vb{q}}(t),\vb{u}(t))
    \,,\quad
    \hat{\vb{q}}(t_0) = \hat{\vb{q}}_0\,, \quad t\in [t_0, t_f]\,,
\end{equation}
where $\hat{\vb{O}}\in\mathbb{R}^{r\times d(r, p)}$ is the operator matrix,
$\vb{d}:\mathbb{R}^{r}\times\mathbb{R}^{p}\to\mathbb{R}^{d(r,p)}$ is a known (typically nonlinear) function of the reduced state and input,
and $d(r, p) \in \mathbb{N}$ is the operator dimension.
The system~\cref{eq:reducedorder} is linear with respect to the operator matrix $\hat{\vb{O}}$ but is usually nonlinear in the reduced state $\hat{\vb{q}}(t)$ and/or the input $\vb{u}(t)$ due to $\vb{d}(\hat{\vb{q}}(t),\vb{u}(t))$.
For example, the quadratic system \cref{eq:polynomial-reducedorder} can be written as \cref{eq:reducedorder} with
$d(r,p) = 1 + r + r^2 + p$ and
\begin{equation}
    \hat{\vb{O}}
    = [~\hat{\vb{c}}~~\hat{\vb{A}}~~\hat{\vb{H}}~~\hat{\vb{B}}~] \in \mathbb{R}^{r \times d(r,p)},
    \qquad
    \vb{d}(\hat{\vb{q}}(t),\vb{u}(t))
    = \left[\begin{array}{c}
        1 \\
        \hat{\vb{q}}(t) \\
        \hat{\vb{q}}(t)\otimes \hat{\vb{q}}(t) \\
        \vb{u}(t)
    \end{array}\right]\in\mathbb{R}^{d(r,p)}.
\end{equation}
There are several ways to approach problems whose nonpolynomial nonlinear structure does not directly motivate the reduced structure~\cref{eq:reducedorder}. Many problems can be posed in polynomial form after a change of variables~\cite{QKPW2020liftAndLearn}, a strategy that has found application in several model reduction contexts (see, e.g.,~\cite{kramer2022nonlinearbalanced} for an example of variable transformation enabling balanced truncation ROMs). OpInf models have also been demonstrated to be effective in some scenarios where the true dynamics can only be partially represented through a polynomial form~\cite{farcas2023opinf4rde,mcquarrie2021combustion,qian2022reduced,SKHW2020romCombustion}. Additionally, a large class of nonpolynomial nonlinear terms may be incorporated via discrete empirical interpolation~\cite{benner2020opinfdeim}.

Given state-input pairs $\{(\vb{q}_j,\vb{u}_j)\}_{j=1}^{m-1}$ describing the solution of (2.1) at $m$ discrete time instances $t_0,\ldots,t_{m-1}$, a reference vector $\bar{\vb{q}}$, and a basis matrix $\vb{V}$, OpInf learns an operator matrix $\hat{\vb{O}}$ by compressing the observed states as $\hat{\vb{q}}_j = \vb{V}\trp(\vb{q}_j - \bar{\vb{q}}) \in\mathbb{R}^{r}$, then solving the regression
\begin{subequations}
\label{eq:opinf-all}
\begin{gather}
    \label{eq:opinf-regression-matrix}
    \min_{\hat{\vb{O}}}\sum_{j=0}^{m-1}\left\|
        \hat{\vb{O}}\vb{d}(\hat{\vb{q}}_{j},\vb{u}_j)
        - \dot{\hat{\vb{q}}}_j
    \right\|_{2}^{2}
    = \min_{\hat{\vb{O}}}\left\|
        \hat{\vb{O}}\vb{D}\trp
        - \dot{\hat{\vb{Q}}}
    \right\|_{F}^{2}
    = \min_{\hat{\vb{O}}}\left\|
        \vb{D}\hat{\vb{O}}\trp
        - \dot{\hat{\vb{Q}}}\trp
    \right\|_{F}^{2}\,,
    \\
    \label{eq:opinf-data-matrix}
    \dot{\hat{\vb{Q}}}
    = \left[\begin{array}{cccc}
        & & & \\
        \dot{\hat{\vb{q}}}_0 &
        \dot{\hat{\vb{q}}}_1 &
        \cdots &
        \dot{\hat{\vb{q}}}_{m-1}
        \\ & & &
    \end{array}\right]
    \in\mathbb{R}^{r\times m}\,,
    \qquad
    \vb{D}
    = \left[\begin{array}{c}
        \vb{d}(\hat{\vb{q}}_0, \vb{u}_0)\trp \\
        \vb{d}(\hat{\vb{q}}_1, \vb{u}_1)\trp \\
        \vdots \\
        \vb{d}(\hat{\vb{q}}_{m-1}, \vb{u}_{m-1})\trp \\
    \end{array}\right]
    \in\mathbb{R}^{m \times d(r, p)}\,.
\end{gather}
\end{subequations}
Here, $\dot{\hat{\vb{q}}}_j\in\mathbb{R}^{r}$ is an estimate of the time derivative $\frac{\textrm{d}}{\textrm{d}t}\hat{\vb{q}}(t) = \vb{V}\trp\frac{\textrm{d}}{\textrm{d}t}{\vb{q}}(t)$ at time $t_j$, usually computed from finite differences of the reduced snapshots $\hat{\vb{q}}_0,\ldots,\hat{\vb{q}}_{m-1}$, and
the matrix $\vb{D}$ is called the data matrix.
Because $\vb{V}$, $\bar{\vb{q}}$, and the observations $\{(\vb{q}_j,\vb{u}_j)\}$ are known, the only unknowns in \cref{eq:opinf-all} are the entries of the operator matrix $\hat{\vb{O}}$.
Unless additional constraints on the unknown $\hat{\vb{O}}$ are imposed, \cref{eq:opinf-all} is a linear least-squares problem since \cref{eq:reducedorder} is linear with respect to $\hat{\vb{O}}$.

As a regression for $\hat{\vb{O}}$, \cref{eq:opinf-regression-matrix} decouples into $r$ independent regressions,
\begin{align}
    \label{eq:single-regression}
    \min_{\hat{\vb{o}}_i}\left\|
        \vb{D}\hat{\vb{o}}_i
        - \vb{z}_i
        \right\|_{2}^{2}\,,
\end{align}
where $\hat{\vb{o}}_i \in \mathbb{R}^{d(r, p)}$ is the $i$-th row of $\hat{\vb{O}}$ and $\vb{z}_{i}\in\mathbb{R}^{k}$ is the $i$-th row of $\dot{\hat{\vb{Q}}}$.
It is often important to regularize the regression \cref{eq:single-regression} in order for the resulting $\hat{\vb{O}}$ to produce a stable ROM \cite{mcquarrie2021combustion,sawant2023piregopinf}.
The most common way to do this is with a Tikhonov regularizer, i.e.,
\begin{align}
    \label{eq:single-regression-regularized}
    \min_{\hat{\vb{o}}_i}\left\{
        \left\|
            \vb{D}\hat{\vb{o}}_i
            - \vb{z}_i
        \right\|_{2}^{2}
        + \left\|
            \vb*{\Gamma}_{i}\hat{\vb{o}}_{i}
        \right\|_{2}^{2}
    \right\}\,,
\end{align}
for some $\vb*{\Gamma}_{i}\in\mathbb{R}^{d(r,p)\times d(r,p)}$.
Alternatively, some scenarios call for constraints to be added to the regression, for example when preserving Hamiltonian or Lagrangian dynamics \cite{gruber2023hamiltonianopinf,sharma2024lagrangian,sharma2022hamiltonianopinf}.

\section{Bayesian operator inference with Gaussian processes}
\label{sec:opinf-bayes}

In deterministic OpInf \cref{eq:single-regression}, the training data are constructed from whatever (reduced) snapshot data are available, be they uniform or sparse in time, smooth or corrupted by observational noise.
Here, we perform GP~regression using the snapshot data and use the updated process to produce time derivative estimates and, ultimately, a probabilistic description for the reduced-order operator matrix $\hat{\vb{O}}$.
The overall strategy is abbreviated as GP-BayesOpInf.
See \apdxref{appendix:gps} for a brief overview of GPs.

\subsection{Bayesian inference}

In this section, we formulate a linear Bayesian inference for the operator matrix $\hat{\vb{O}}$ by separately considering each of its rows $\hat{\vb{o}}_{1},\ldots,\hat{\vb{o}}_{r}$, following the OpInf paradigm of learning each $\hat{\vb{o}}_i$ individually.

\subsubsection{Likelihood and prior}

We begin by modeling the reduced states $\hat{q}_1(t),\ldots,\hat{q}_r(t)$ with independent GPs:
\begin{equation}
    \label{eq:gp-single}
    \hat{q}_i
    \sim \mathcal{GP}(0,\kappa_i(\cdot,\cdot))\,,
    \qquad
    i=1,\ldots,r,
\end{equation}
in which the mean function is assumed to be zero and $\kappa_i:\mathbb{R}\times\mathbb{R}\to\mathbb{R}$
is a kernel function specifying the covariance over the GP.
If the kernel $\kappa_i$ is at least once differentiable in each argument, then the time-derivative of each $\hat{q}_i$ also follows a GP, and we have the joint processes~\cite{raissi2017GPsforlinearPDEs,rasmussen2005GPsforMLbook}
\begin{equation}
    \label{jointfunc-new}
    \left[\begin{array}{c}
        \hat{q}_i \\ \dot{\hat{q}}_i
    \end{array}\right]
    \sim \text{vec-}\mathcal{GP}\left(
    \left[\begin{array}{c}
        0 \\ 0
    \end{array}\right],
    \left[\begin{array}{cc}
        \kappa_i(\cdot,\cdot)
        &
        \partial_{2}\kappa_i(\cdot,\cdot)
        \\
        \partial_{1}\kappa_i(\cdot,\cdot)
        &
        \partial_{1}\partial_{2}\kappa_i(\cdot,\cdot)
        \end{array}\right]
    \right)\,,
    \quad
\end{equation}
in which $\partial_1$ and $\partial_2$ respectively denote the partial derivatives with respect to the first and the second arguments.

Let $\vb*{t} = (t_{0},\ldots,t_{m - 1})\trp \in \mathbb{R}^{m}$ denote the $m\in\mathbb{N}$ time instances for which we have (reduced) snapshot data $\hat{\vb{q}}_0,\ldots,\hat{\vb{q}}_{m-1}\in\mathbb{R}^{r}$.
Define $\vb*{Y}_i=\hat{q}_i(\vb*{t})$ to be the random $m$-vector whose entries are the Gaussian random variables produced by evaluating the GP $\hat{q}_i$ at each time instance in $\vb*{t}$, and let $\vb{y}_i\in\mathbb{R}^{m}$ denote the snapshot data for the $i$-th reduced state $\hat{q}_i$ at times $\vb*{t}$, that is,
\begin{align}
\left[\begin{array}{cccc}
        & & & \\
        \hat{\vb{q}}_{0} &
        \hat{\vb{q}}_{1}
        & \cdots &
        \hat{\vb{q}}_{m-1}
        \\ & & &
    \end{array}\right]
&= \left[\begin{array}{c}
        \vb{y}_{1}\trp \\
        \vdots \\
        \vb{y}_{r}\trp
    \end{array}\right]
    \in \mathbb{R}^{r\times m}.
\end{align}
Note that $\vb{y}_i$ is an observation of $\vb*{Y}_i$.
Likewise, let $\vb*{t}' = (t_{0}',\ldots,t_{m'-1}')\trp \in \mathbb{R}^{m'}$ be the $m'\in\mathbb{N}$ time instances (typically with $m' \ge m$ and $[t_{0}',t_{m'-1}']\subset[t_{0},t_{m-1}]$) where the GP is to produce estimates, and define $\vb*{Z}_{i} = \dot{\hat{q}}_i(\vb*{t}')$ to be the random $m'$-vector with entries given by the time derivative GP $\dot{\hat{q}}_{i}$ at times $\vb*{t}'$.

The entries of the observational data vector $\vb{y}_i$ are noisy, which we model by an additive white noise term with Gaussian distribution $\mathcal{N}(0,\chi_i)$.
With this setting, the joint distribution of $\vb*{Y}_i$ and $\vb*{Z}_i$ is deduced from~\cref{jointfunc-new} to be
\begin{align}
    \label{eq:joint-gp}
    \begin{bmatrix}
        \vb*{Y}_i \\ \vb*{Z}_i
    \end{bmatrix} & \sim \mathcal{N}\left(
    \begin{bmatrix}
        \vb{0} \\ \vb{0}
    \end{bmatrix},
    \begin{bmatrix}
        \vb{K}^{yy}_i & \vb{K}^{yz}_i \\
        \vb{K}^{zy}_i & \vb{K}^{zz}_i
    \end{bmatrix}
    \right)\,,
    &\text{with}&&
    \begin{aligned}
    \vb{K}_i^{yy}
    &= \kappa_i(\vb*{t},\vb*{t})+\chi_i\vb{I}_{m}
    \in \mathbb{R}^{m\times m}
    \,,\\
    \vb{K}_i^{zy}
    &= \partial_1\kappa_i(\vb*{t}',\vb*{t})
    = (\vb{K}_i^{yz})\trp
    \in \mathbb{R}^{m'\times m}
    \,,\\
    \vb{K}_i^{zz}
    &= \partial_1\partial_2\kappa_i(\vb*{t}',\vb*{t}')
    \in \mathbb{R}^{m'\times m'}\,,
    \end{aligned}
\end{align}
where $\vb{I}_{m}$ is the $m \times m$ identity matrix
and, for two vectors $\vb*{a} = (a_1,\ldots,a_{n_a})\trp$ and $\vb*{b} = (b_1,\ldots,b_{n_b})\trp$, the shorthand $\kappa(\vb*{a},\vb*{b})$ denotes the $n_{a}\times n_{b}$ matrix whose $(i,j)$-th entry is $\kappa(a_i,b_j)$.

\begin{remark}[Choosing $m'$]
The integer $m'$ must be selected \emph{a priori} by the practitioner. If $m'$ is too large, fitting and using each GP becomes unnecessarily expensive; if $m'$ is too small, the subsequent Bayesian inference may have insufficient data and produce unsatisfactory results.
For the numerical examples in \Cref{sec:results} and \Cref{sec:odeparamestimate}, we base $m'$ on the interval length $t_{m-1} - t_0$ and uniformly space the entries of $\vb*{t}'$ in the interval $[t_{0},t_{m-1}]$, i.e., $t_{0} = t'_{0} < t'_{1} < \ldots < t'_{m' - 1} = t_{m - 1}$ with $t'_{j+1} = t'_{j} + (t_{m - 1} - t_{0})/(m'-1)$.
\end{remark}

For a given reduced-order operator vector $\hat{\vb{o}}_{i} \in \mathbb{R}^{d(r, p)}$ and a data matrix $\tilde{\vb{D}}\in\mathbb{R}^{m' \times d(r, p)}$ corresponding to the estimation times $\vb*{t}'$, let $\vb*{R}_{i} := \vb*{Z}_{i} - \tilde{\vb{D}}\hat{\vb{o}}_i$ be the random $m'$-vector describing the residual for the $i$-th reduced-order system equation~\cref{eq:reducedorder} at times $\vb*{t}'$.
We may now define a likelihood function, not only on the direct state observations $\vb*{Y}_i = \vb{y}_i$, but also on the alignment of the time derivatives $\vb*{Z}_i$ with $\tilde{\vb{D}}\hat{\vb{o}}_i$, i.e., the satisfaction of residuals $\vb*{R}_{i} = \vb{0}$:
\begin{equation}
    \label{eq:likelihood-new}
    \pi_\texttt{like}(\hat{\vb{o}}_i) := p(\vb*{Y}_i = \vb{y}_i, \vb*{R}_i = \vb{0} \mid \hat{\vb{o}}_i,\tilde{\vb{D}},
    \vb*{\theta}_i
    )
     \propto \exp\left(-\frac{1}{2}
    \begin{bmatrix}
        \vb{y}_i \\ \tilde{\vb{D}}\hat{\vb{o}}_i
    \end{bmatrix}\trp
    \begin{bmatrix}
        \vb{K}^{yy}_i & \vb{K}^{yz}_i \\
        \vb{K}^{zy}_i & \vb{K}^{zz}_i
    \end{bmatrix}^{-1}
    \begin{bmatrix}
        \vb{y}_i \\ \tilde{\vb{D}}\hat{\vb{o}}_i
    \end{bmatrix}
    \right)\,,
\end{equation}
which is derived by substituting $\vb*{Z}_i = \tilde{\vb{D}}\hat{\vb{o}}_i$ (i.e., $\vb*{R}_i = \vb{0}$) into \cref{eq:joint-gp}.
Here, $\vb*{\theta}_{i}$ consists of the noise variance $\chi_i$ and the hyperparameters of the kernel $\kappa_i$.
The likelihood depends on $\vb*{\theta}_{i}$ through the definition of the matrices $\vb{K}_{i}^{yy}$, $\vb{K}_{i}^{yz}$, $\vb{K}_{i}^{zy}$, and $\vb{K}_{i}^{zz}$
in \cref{eq:joint-gp}.

\begin{remark}[Interpretation of the likelihood]
In the likelihood definition \eqref{eq:likelihood-new}, the direct observations $\vb*{Y}_i = \vb{y}_i$ introduce the factor of GP regression where noisy and scarce data can be smoothed and interpolated, while $\vb*{R}_i = \vb{0}$ incorporates the constraints of differential equations into the inference process, explicitly ensuring physical consistency in the resulting ROMs.
\end{remark}

\begin{remark}[Gaussianity preserved by affine reconstruction]
The independent reduced state GP representations $\hat{q}_1,\ldots,\hat{q}_r$ of~\cref{eq:gp-single} can be expressed together as a single vector-valued GP,
\begin{subequations}
\begin{align}
    \hat{\vb{q}}(t)
    \sim \text{vec-}\mathcal{GP}(\vb{0}, \text{diag}(\kappa_1(t,t'),\ldots,\kappa_r(t,t')))\,.
\end{align}
If these reduced states are used in the approximation~\cref{eq:state-approx}, the reconstruction of full-order states is given by
\begin{align}
    \vb{q}(t)
    \approx \vb{V}\hat{\vb{q}}(t) +\bar{\vb{q}}
    \sim \text{vec-}\mathcal{GP}(\bar{\vb{q}}, \vb{V}\text{diag}(\kappa_1(t,t'),\ldots,\kappa_r(t,t'))\vb{V}\trp)\,,
\end{align}
\end{subequations}
in which the covariance $\vb{V}\text{diag}(\kappa_1(t,t'),\ldots,\kappa_r(t,t'))\vb{V}\trp$ is a low-rank approximation of $\operatorname{Cov}[\vb{q}(t),\vb{q}(t')]$.
Hence, the GP-based modeling for the reduced states
provides a quantification of uncertainties in the states on the full-order level as well.
\end{remark}

We also define the prior on the reduced-order operators  $\hat{\vb{o}}_i$ to be independent normal distributions, i.e.,
\begin{equation}
    \label{eq:bayesian-prior}
    \pi_\texttt{prior}(\hat{\vb{o}}_i):=p(\hat{\vb{o}}_i \mid \vb*{\gamma}_i )
    \propto \exp\left( -\frac{1}{2} \|\vb*{\Gamma}_i\hat{\vb{o}}_i\|_2^2 \right)\,,
    \quad \text{with} \quad
    \vb*{\Gamma}_i
    := \text{diag}(\vb*{\gamma}_i)\,,
\end{equation}
where the entries of $\vb*{\gamma}_i\in \mathbb{R}_{+}^{d(r,p)}$ are nonnegative and are considered hyperparameters.

\subsubsection{Posterior}

Combining \cref{eq:likelihood-new}--\cref{eq:bayesian-prior}, Bayes' rule gives the posterior distribution
\begin{equation}
    \label{eq:posterior-inference}
    \pi_\texttt{post}(\hat{\vb{o}}_i):= p(\hat{\vb{o}}_i
        \mid \vb*{Y}_i = \vb{y}_i, \vb*{R}_i = \vb{0}, \tilde{\vb{D}},\vb*{\theta}_i,\vb*{\gamma}_i)
    \propto \pi_\texttt{like}(\hat{\vb{o}}_i)\,\pi_\texttt{prior}(\hat{\vb{o}}_i)\,.
\end{equation}
From Schur's formula for block matrix inversion \cite{abadir2005matrix}, one has
\begin{align}
    \label{eq:KsandWs}
    \begin{bmatrix}
        \vb{K}^{yy}_i & \vb{K}^{yz}_i \\
        \vb{K}^{zy}_i & \vb{K}^{zz}_i
        \end{bmatrix}^{-1} & =
    \begin{bmatrix}
        \vb{W}^{yy}_i & \vb{W}^{yz}_i \\
        \vb{W}^{zy}_i & \vb{W}^{zz}_i
        \end{bmatrix}
    \,,
    &&
    \begin{aligned}
    \vb{W}^{yy}_i
    &= \left(\vb{K}^{yy}_i - \vb{K}^{yz}_i(\vb{K}^{zz}_i)^{-1}\vb{K}^{zy}_i \right)^{-1}
    \,,\\
    \vb{W}^{zy}_i
    &= -\vb{K}^{zz}_i\vb{K}^{zy}_i (\vb{K}^{yy}_i)^{-1}= (\vb{W}^{yz}_i)\trp
    \,,\\
    \vb{W}^{zz}_i
    &= \left(\vb{K}^{zz}_i - \vb{K}^{zy}_i(\vb{K}^{yy}_i)^{-1}\vb{K}^{yz}_i \right)^{-1}\,.
    \end{aligned}
\end{align}
Therefore,
\begin{equation}
\begin{split}
\pi_\texttt{post}(\hat{\vb{o}}_i)
    & \propto \pi_\texttt{like}(\hat{\vb{o}}_i)\,\pi_\texttt{prior}(\hat{\vb{o}}_i)
    \\
    &\propto \exp\left(-\frac{1}{2}\left(
        \vb{y}_i\trp\vb{W}^{yy}_i\vb{y}_i
        + 2\hat{\vb{o}}_i\trp\tilde{\vb{D}}\trp\vb{W}_i^{zy}\vb{y}_i
        + \hat{\vb{o}}_i\trp\tilde{\vb{D}}\trp\vb{W}_i^{zz}\tilde{\vb{D}}\hat{\vb{o}}_i
        + \hat{\vb{o}}_i\trp\vb*{\Gamma}_i\trp \vb*{\Gamma}_i\hat{\vb{o}}_i
    \right)\right)
    \\
    &\propto \exp\left(-\frac{1}{2}\left(
        \hat{\vb{o}}_i
        - \vb*{\mu}_i \right)\trp\vb*{\Sigma}_i^{-1} \left(\hat{\vb{o}}_i - \vb*{\mu}_i\right)
    \right)\,,
    \\
    \text{i.e.,}\quad
    \pi_\texttt{post}(\hat{\vb{o}}_i) &
    = \mathcal{N}\left(\hat{\vb{o}}_i \mid \vb*{\mu}_i,\vb*{\Sigma}_i\right)\,.
\end{split}
\end{equation}
Here, the posterior mean of $\hat{\vb{o}}_i$ coincides with a \emph{maximum a posteriori} (MAP) estimator, written as
\begin{equation}
\label{eq:posterior-mean}
\begin{split}
    \vb*{\mu}_i
    =\left(\hat{\vb{o}}_i\right)_\text{MAP}
    &= -(\tilde{\vb{D}}\trp\vb{W}_i^{zz}\tilde{\vb{D}}+\vb*{\Gamma}_{i}\trp\vb*{\Gamma}_{i})^{-1}\tilde{\vb{D}}\trp\vb{W}_i^{zy}\vb{y}_i
    \\
    &= (\tilde{\vb{D}}\trp\vb{W}_i^{zz}\tilde{\vb{D}}+\vb*{\Gamma}_{i}\trp\vb*{\Gamma}_{i})^{-1}\tilde{\vb{D}}\trp\vb{W}^{zz}_i \vb{K}^{zy}_i (\vb{K}^{yy}_i)^{-1}\vb{y}_i
\,,
\end{split}
\end{equation}
and the posterior covariance matrix of $\hat{\vb{o}}_i$ is
\begin{equation}
    \label{eq:posterior-covariance}
    \vb*{\Sigma}_i = \left(\tilde{\vb{D}}\trp\vb{W}_i^{zz}\tilde{\vb{D}}+\vb*{\Gamma}_{i}\trp\vb*{\Gamma}_{i}\right)^{-1}\,.
\end{equation}
Note that
\begin{align}
    \label{eq:gpreconstruct}
    \tilde{\vb{z}}_i
    := \vb{K}^{zy}_i (\vb{K}^{yy}_i)^{-1}\vb{y}_i
    \in \mathbb{R}^{m'}
\end{align}
describes the GP estimate of the time derivative $\frac{\textrm{d}}{\textrm{d}t}\hat{q}_i$ evaluated over $\vb*{t}'$ (see \cref{eq:gp-timederivative}).
Thus, \cref{eq:posterior-mean} can be abbreviated as
$
    \vb*{\mu}_i
    = \vb*{\Sigma}_i \tilde{\vb{D}}\trp\vb{W}^{zz}_i \tilde{\vb{z}}_i\,.
$

\subsubsection{Least squares formulation}
\label{sec:leastsquares}

The posterior mean $\vb*{\mu}_i$ can be interpreted as a Tikhonov-regularized, generalized least-squares estimator of $\hat{\vb{o}}_i$ with the weight matrix
$
    \vb{W}_i^{zz},
$
whose inverse is the posterior covariance matrix of
$
    \tilde{\vb{z}}_i
$
through the GP regression of $(\vb*{t},\vb{y}_i)$.
Specifically, $\vb*{\mu}_{i}$ solves the problem
\begin{align}
    \label{eq:generalizedlstsq}
    \min_{\vb*{\eta}}\left\{
        \left\|\tilde{\vb{D}}\vb*{\eta} - \tilde{\vb{z}}_{i}\right\|_{\vb{W}_{i}^{zz}}^{2}
        + \left\|\vb*{\Gamma}_i\vb*{\eta}\right\|_{2}^{2}
    \right\},
\end{align}
where $\left\|\bullet\right\|_{\vb{W}_{i}^{zz}} = (\bullet\trp\vb{W}_{i}^{zz}\bullet)^{1/2}$.
Furthermore, an equivalent $2$-norm least-squares problem is given by
\begin{align}
    \label{eq:glstsq2norm}
    \min_{\vb*{\eta}}
    &\left\|\left[\begin{array}{c}
        \left(\vb{W}_{i}^{zz}\right)^{1/2}\tilde{\vb{D}}
        \\ \vb*{\Gamma}_i
    \end{array}\right]\vb*{\eta}
    - \left[\begin{array}{c}
        \left(\vb{W}_{i}^{zz}\right)^{1/2}\tilde{\vb{z}}_{i}
        \\ \vb{0}
    \end{array}\right]
    \right\|_{2}^{2}\,.
\end{align}
We use this least-squares formulation in the computation of $\vb*{\mu}_i$ as it is more stable than solving the generalized normal equations in \cref{eq:posterior-mean}. This requires computing $(\vb{W}_{i}^{zz})^{1/2}$, which can be obtained by taking an eigendecomposition
\begin{align}
    \label{eq:eigendecomposition}
    (\vb{W}_{i}^{zz})^{-1}
    = \vb{K}^{zz}_i - \vb{K}^{zy}_i(\vb{K}^{yy}_i)^{-1}\vb{K}^{yz}_i
    = \vb*{\Psi}_{i}\vb*{\Xi}_{i}\vb*{\Psi}_{i}\trp
    \,,
\end{align}
then setting $(\vb{W}_{i}^{zz})^{1/2} = \vb*{\Psi}_{i}\vb*{\Xi}_{i}^{-1/2}\vb*{\Psi}_{i}\trp$.

\begin{remark}[Regularization of the covariance computation]
\label{remark:gpreg}
In practice, computing \cref{eq:eigendecomposition} directly may be unstable due to the round-off errors introduced from solving for $\vb{C}_i := \vb{K}^{zy}_i(\vb{K}^{yy}_i)^{-1}\vb{K}^{yz}_i$. As a remedy, we symmetrize $\vb{C}_i$ and add a small amount of numerical regularization, replacing \cref{eq:eigendecomposition} with
\begin{align}
    (\vb{W}_{i}^{zz})^{-1}
    &= \vb{K}^{zz}_i - \frac{1}{2}\left(\vb{C}_i + \vb{C}_i\trp\right) + \tau\vb{I}_{m'}
    = \vb*{\Psi}_{i}\vb*{\Xi}_{i}\vb*{\Psi}_{i}\trp,
\end{align}
where $0 < \tau \ll 1$.
The constant $\tau$ should be just large enough that the eigendecomposition of the above matrix yields all positive eigenvalues.
The numerical examples in \Cref{sec:results} and \Cref{sec:odeparamestimate} use $\tau = 10^{-8}$.
\end{remark}

\begin{algorithm}[t]
\begin{algorithmic}[1]
\Procedure{OpPost}{
    \newline\phantom{---}
    Data matrix $\tilde{\vb{D}}\in\mathbb{R}^{m'\times d}$,
    \newline\phantom{---}
    GP time derivative estimates $\tilde{\vb{z}}_{1},\ldots,\tilde{\vb{z}}_{r}\in\mathbb{R}^{m'}$,
    \newline\phantom{---}
    Weight matrices $(\vb{W}_{1}^{zz})^{1/2},\ldots,(\vb{W}_{r}^{zz})^{1/2}\in\mathbb{R}^{m'\times m'}$,
\newline\phantom{---}
    Prior variances $\vb*{\gamma}_1,\ldots,\vb*{\gamma}_r \in \mathbb{R}^{d}$
    \newline
    }

    \For{$i = 1, \ldots, r$}
        \State $\vb*{\Gamma}_{i} \gets \operatorname{diag}(\vb*{\gamma}_{i})$
        \State $\vb*{\mu}_i \gets \operatorname{argmin}_{\vb*{\eta}}\big\{
            \big\|
                \tilde{\vb{D}}\vb*{\eta}
                - \tilde{\vb{z}}_{i}
            \big\|_{\vb{W}_{i}^{zz}}^{2}
            + \big\|\vb*{\Gamma}_i\vb*{\eta}\big\|_{2}^{2}
        \big\}$
            \Comment{Solve numerically via \cref{eq:glstsq2norm}}.
\State $\vb*{\Sigma}_{i} \gets (\tilde{\vb{D}}\trp\vb{W}_{i}^{zz}\tilde{\vb{D}} + \vb*{\Gamma}_{i}\trp\vb*{\Gamma}_{i})^{-1}$
\EndFor
    \State \textbf{return} $(\vb*{\mu}_1,\vb*{\Sigma}_1),\ldots,(\vb*{\mu}_r,\vb*{\Sigma}_r)$
\EndProcedure
\end{algorithmic}
\caption{Given data, compute the moments of the posterior operator distribution.}
\label{alg:OperatorPosterior}
\end{algorithm}

\subsection{Probabilistic reduced-order model predictions}

Now that the reduced operator vectors $\hat{\vb{o}}_i$ are estimated with probability distributions $\pi_\texttt{post}(\hat{\vb{o}}_i)$, $1 \leq i \leq r$, they may be used to produce predictions over time through the ROM~\cref{eq:reducedorder}.
In fact, the posterior distribution for the entire reduced operator matrix $\hat{\vb{O}}$ can be written as
\begin{subequations}
\begin{equation}
\begin{split}
    \pi_\texttt{post}(\hat{\vb{O}}) & = p\left(\hat{\vb{O}}
    ~\big |~ \left\{\vb*{Y}_i=\vb{y}_i,\vb*{R}_i=\vb{0}\right\}_{i=1}^{r}, \tilde{\vb{D}},\left\{\vb*{\theta}_i,\vb*{\gamma}_i\right\}_{i=1}^{r}\right) \\
    & \propto p\left(\left\{\vb*{Y}_i=\vb{y}_i,\vb*{R}_i=\vb{0}\right\}_{i=1}^{r} ~\big |~ \hat{\vb{O}}, \tilde{\vb{D}}, \left\{\vb*{\theta}_i \right\}_{i=1}^{r} \right) p\left(\hat{\vb{O}}~\big |~\left\{\vb*{\gamma}_i \right\}_{i=1}^{r} \right) \\
    & = \prod_{i=1}^{r} p\left(\vb*{Y}_i=\vb{y}_i,\vb*{R}_i=\vb{0} ~\big |~ \hat{\vb{o}}_i, \tilde{\vb{D}}, \vb*{\theta}_i \right) \prod_{i=1}^{r} p\left( \hat{\vb{o}}_i \mid \vb*{\gamma}_i \right) = \prod_{i=1}^{r}\pi_\texttt{like}(\hat{\vb{o}}_i) \pi_\texttt{prior}(\hat{\vb{o}}_i)\,.
\end{split}
\end{equation}
The likelihood separates into a product because (1) the OpInf regression~\cref{eq:opinf-regression-matrix}--\cref{eq:opinf-data-matrix} decouples the inference of the operator matrix $\hat{\vb{O}}$ into its rows $\hat{\vb{o}}_{1},\ldots,\hat{\vb{o}}_{r}$,
and (2) we model the reduced modes independently with individual one-dimensional GPs~\cref{eq:gp-single}.
The separation of the prior is a result of the independence in the prior definition~\eqref{eq:bayesian-prior}. Hence the posterior of $\hat{\vb{O}}$ can be treated row-by-row independently, i.e.,
\begin{equation}
    \pi_\texttt{post}(\hat{\vb{O}}) = \prod_{i=1}^{r}  \pi_\texttt{post}(\hat{\vb{o}}_i)\,,
\end{equation}
\end{subequations}
and $\hat{\vb{O}}$ is sampled from the product of the posterior distributions $\pi_\texttt{post}(\hat{\vb{o}}_1),\ldots,\pi_\texttt{post}(\hat{\vb{o}}_r)$.

For a given sample of $\hat{\vb{O}}$, we evaluate the reduced state $\hat{\vb{q}}(t)$ by solving the corresponding ROM~\cref{eq:reducedorder}. The reduced state $\hat{\vb{q}}(t)$ therefore has the following distribution:
\begin{equation}
\begin{split}
    p\left(\hat{\vb{q}}(t)
    ~\big |~ \left\{\vb*{Y}_i=\vb{y}_i,\vb*{R}_i=\vb{0}\right\}_{i=1}^{r}, \tilde{\vb{D}},\left\{\vb*{\theta}_i,\vb*{\gamma}_i\right\}_{i=1}^{r}\right) =~& \int p(\hat{\vb{q}}(t) \mid \hat{\vb{O}})
    \,\pi_\texttt{post}(\hat{\vb{O}} ) ~\dd\hat{\vb{O}}
    \\
    =~& \int p(\hat{\vb{q}}(t) \mid \hat{\vb{o}}_1,\ldots,\hat{\vb{o}}_r)
    \prod_{i=1}^{r} \pi_\texttt{post}(\hat{\vb{o}}_i )~\dd\hat{\vb{o}}_i\,,
\end{split}
\end{equation}
in which all the $\hat{\vb{o}}_i$'s are marginalized.
Here, the notation $p(\hat{\vb{q}}(t)|\,\hat{\vb{O}})$ indicates a deterministic model evaluation for a given instance of $\hat{\vb{O}}$; indeed, one can write the reduced state $\hat{\vb{q}}(t)$ as a function of $\hat{\vb{O}}$, i.e., $\hat{\vb{q}}(t;\hat{\vb{O}})$.
Note, however, that the distribution of $\hat{\vb{q}}$ is not necessarily Gaussian, even though the posterior of $\hat{\vb{O}}$ is, primarily because $\hat{\vb{q}}(t;\hat{\vb{O}})$ is nonlinear with respect to $\hat{\vb{O}}$.
Having constructed a posterior distribution for $\hat{\vb{O}}$, we use Monte Carlo sampling over the operator posterior to estimate the moments of $\hat{\vb{q}}(t)$.
Specifically, an individual sample $\hat{\vb{O}}_{k}$ defines an instance of the (deterministic) ROM~\cref{eq:reducedorder},
\begin{align*}
    \frac{\textup{d}}{\textup{d}t} \hat{\vb{q}}(t)
    = \hat{\vb{O}}_{k}\vb{d}(\hat{\vb{q}}(t),\vb{u}(t))
    \,,\quad
    \hat{\vb{q}}(t_0) = \hat{\vb{q}}_0\,, \quad t\in [t_0, t_f]\,,
\end{align*}
which is solved numerically for $\hat{\vb{q}}(t)$ over $[t_0,t_f]$ using standard time-stepping integration methods.
Each numerical solve is computationally inexpensive due to the low dimensionality of the system.
After repeating this procedure for $k=1,\ldots,K$, the collection of $K$ trajectories
is used to estimate the mean and higher-order moments of the reduced state at time $t$, thereby establishing an uncertainty indicator for the model predictions.

\subsection{Multiple trajectories}
\label{sec:multitrajectory}

In many situations we have data for $\ell>1$ trajectories of a single system due to variable initial conditions or input functions (but with the same set of system parameters, such as material properties).
In this scenario, observed states from each of the $\ell$ trajectories are used together to construct a global basis of rank $r$, then $r$ GPs---one for each reduced mode, as before---are trained for each trajectory.
The observation times $\vb*{t}$ may vary by trajectory, but we use the same estimation times $\vb*{t}'$ for each trajectory and mode.

The GP estimates for the $i$-th mode from each trajectory are used jointly to learn a posterior for the operator matrix row $\hat{\vb{o}}_{i}$.
Using a superscript $(l)$ to represent the $l$-th trajectory over $1\leq l \leq \ell$,
for each $i = 1, \ldots, r$ we are concerned with the posterior
\begin{equation}
    \begin{split}
    \pi_\texttt{post}(\hat{\vb{o}}_i) =  p\left(\hat{\vb{o}}_i ~\Big |~
        \left\{\vb*{Y}_i^{(l)} = \vb{y}_i^{(l)}, \vb*{R}_i^{(l)} = \vb{0}, ~\tilde{\vb{D}}^{(l)},\vb*{\theta}_i^{(l)}\right\}_{l=1}^{\ell},
\vb*{\gamma}_i\right)
    \propto
    \pi_\texttt{prior}(\hat{\vb{o}}_i)\prod_{l=1}^{\ell} \pi_\texttt{like}^{(l)}(\hat{\vb{o}}_i)\,,
    & \\
    \text{with}\quad \pi_\texttt{like}^{(l)}(\hat{\vb{o}}_i)= p\left(\vb*{Y}_i^{(l)} = \vb{y}_i^{(l)},\vb*{R}_i^{(l)} = \vb{0}
        ~\big |~ \hat{\vb{o}}_i\,,\tilde{\vb{D}}_i^{(l)}, \vb*{\theta}_i^{(l)}\right)\,, &
    \end{split}
\end{equation}
where $\tilde{\vb{D}}^{(l)}$ is formed from the GP state estimates for the $l$-th trajectory, i.e., $\tilde{\vb{y}}_1^{(l)},\ldots,\tilde{\vb{y}}_r^{(l)}$.
The posterior mean and covariance respectively become
\begin{subequations}
\begin{align}
    \label{mapestimate-multipletrajectories}
    \vb*{\mu}_i
    &= \left(\sum_{l = 1}^{\ell}(\tilde{\vb{D}}^{(l)})\trp \vb{W}_{i}^{zz\,(l)}\tilde{\vb{D}}^{(l)}+\vb*{\Gamma}_i\trp \vb*{\Gamma}_i
    \right)^{-1} \sum_{l = 1}^{\ell}(\tilde{\vb{D}}^{(l)})\trp \vb{W}_{i}^{zz\,(l)} \tilde{\vb{z}}_{i}^{(l)}\,,
    \\
    \vb{\Sigma}_{i}
     &= \left(\sum_{l = 1}^{\ell}(\tilde{\vb{D}}^{(l)})\trp \vb{W}_{i}^{zz\,(l)}\tilde{\vb{D}}^{(l)}+\vb*{\Gamma}_i\trp \vb*{\Gamma}_i
    \right)^{-1},
\end{align}
\end{subequations}
in which $\tilde{\vb{z}}_{i}^{(l)} := \vb{K}^{zy\,(l)}_{i} (\vb{K}^{yy\,(l)}_{i})^{-1}\vb{y}_{i}^{(l)}$
and
$
    \vb{W}_{i}^{zz\,(l)}
    := (\vb{K}^{zz\,(l)}_i - \vb{K}^{zy\,(l)}_i(\vb{K}^{yy\,(l)}_i)^{-1}\vb{K}^{yz\,(l)}_i)^{-1}
$
is the corresponding weight matrix.
In this case, the generalized least-squares problem \cref{eq:generalizedlstsq} for $\vb*{\mu}_{i}$ becomes
\begin{align}
    \min_{\vb*{\eta}}\left\{
    \sum_{l = 1}^{\ell}\left\|
        \tilde{\vb{D}}^{(l)}\vb*{\eta} - \tilde{\vb{z}}_{i}^{(l)}
    \right\|_{\vb{W}_{i}^{zz\,(l)}}^{2}
    + \left\|\vb{\Gamma}_{i}\vb*{\eta}\right\|_{2}^{2} \right\}.
\end{align}
We can write this in a form that is similar to \cref{eq:generalizedlstsq} by noting that
\begin{equation}
    \label{eq:multitrajectory-concatenation}
    \sum_{l = 1}^{\ell}\left\|\tilde{\vb{D}}^{(l)}\vb*{\eta} - \tilde{\vb{z}}_{i}^{(l)}\right\|_{\vb{W}_{i}^{zz\,(l)}}^{2} = \left\|
        \left[\begin{array}{c}
            \tilde{\vb{D}}^{(1)} \\ \vdots \\ \tilde{\vb{D}}^{(\ell)}
        \end{array}\right]
        \vb*{\eta}
        - \left[\begin{array}{c}
            \tilde{\vb{z}}_{i}^{(1)} \\ \vdots \\ \tilde{\vb{z}}_{i}^{(\ell)}
        \end{array}\right]
    \right\|_{\vb{W}_i}^{2}\,,
\end{equation}
where $\vb{W}_{i}$ is the block-diagonal matrix $\vb{W}_{i} = \operatorname{diag}(\vb{W}_{i}^{zz\,(1)},\cdots,\vb{W}_{i}^{zz\,(\ell)})$.
Hence, the multiple-trajectory case can be written as the original case with a concatenated data matrix, block diagonal weights, and concatenated time derivative estimates.

\begin{remark}[Parametric problems]
The strategy presented here assumes that the underlying reduced operators comprising $\hat{\vb{O}}$ are constant across all trajectories, which is the case when only initial conditions or input terms are altered.
Parametric problems in which the reduced operators vary by trajectory require additional specialized treatment (see, e.g., \cite{mcquarrie2023popinf}) and remain a subject for future work.
\end{remark}

\section{Algorithmic procedure}
\label{sec:procedure}

The inference to define the posterior parameters \cref{eq:posterior-mean}--\cref{eq:posterior-covariance} depends on the choice of kernel hyperparameters $\vb*{\theta}_1,\ldots,\vb*{\theta}_r$, the data matrix $\tilde{\vb{D}}$, and the prior variance parameters $\vb*{\gamma}_1,\ldots,\vb*{\gamma}_r$.
This section provides strategies for selecting each of these quantities and details the associated numerical algorithms.
For simplicity we return to the single-trajectory PDE setting, but the approach readily generalizes to cases with data from multiple trajectories
(see~\apdxref{appendix:multiple-trajectories}).

\subsection{Kernel hyperparameter selection}
\label{sec:gpfitting}

For each $i=1,\ldots,r$, we select $\vb*{\theta}_i$ through the one-dimensional GP regression of the (noisy, reduced) snapshot data $(\vb*{t}, \vb{y}_i)$.
Specifically, we approximately solve the simplified maximum marginal log-likelihood problem
\begin{align}
    \label{eq:maxmlhd}
    \begin{aligned}
    \vb*{\theta}_i^{*}
    &= \operatorname{argmax}_{\vb*{\theta}_i}
    \mathcal{L}(\vb*{\theta}_i; \vb*{t},\vb{y}_i)
    \\
    &= \operatorname{argmax}_{\vb*{\theta}_i}\left\{-\frac{1}{2}\vb{y}_i\trp\left(\kappa(\vb*{t},\vb*{t};\vb*{\theta}_i) + \chi_i\vb{I}_m\right)^{-1}\vb{y}_i
    -\frac{1}{2}\log\big|\kappa(\vb*{t},\vb*{t};\vb*{\theta}_i) + \chi_i\vb{I}_m\big|
    -\frac{m}{2}\log(2\pi)\right\}\,,
    \end{aligned}
\end{align}
where $\chi_i$ is included in $\vb*{\theta}_i$.
Exact solutions of this nonconvex optimization are difficult to obtain, but there are several open-source resources for obtaining satisfactory approximate solutions.
Our implementation utilizes the Python library \texttt{scikit-learn}~\cite{sklearn-api}, which applies the limited-memory bound-constrained BFGS algorithm~\cite{byrd1991lbfgsb} available from \texttt{scipy.optimize}~\cite{2020SciPy-NMeth}.

\subsection{Data selection}
\label{sec:data_selection}

\begin{algorithm}[t]
    \begin{algorithmic}[1]
    \Procedure{GPFit}{
    \newline\phantom{---}
        Observation times $\vb*{t}\in\mathbb{R}^{m}$,
        \newline\phantom{---}
        Observed data $\vb{y} \in \mathbb{R}^{m}$,
        \newline\phantom{---}
        Estimation times $\vb*{t}'\in\mathbb{R}^{m'}$,
        \newline\phantom{---}
        Numerical regularization $0 < \tau \ll 1$
        \newline
        }
        \State $\vb*{\theta}^{*} \gets \operatorname{argmax}_{\vb*{\theta}}\mathcal{L}(\vb*{\theta};\vb*{t},\vb{y})$
            \Comment{Maximize the marginal likelihood \cref{eq:maxmlhd}.}
        \State $\vb{K}^{yy} \gets \kappa(\vb*{t},\vb*{t};\vb*{\theta}^{*}) + \chi\vb{I}_{m}$
        \State $\vb{K}^{zy} \gets \partial_{1}\kappa(\vb*{t}',\vb*{t};\vb*{\theta}^{*})$
        \State $\vb{K}^{zz} \gets \partial_{1}\partial_{2}\kappa(\vb*{t}',\vb*{t}';\vb{\theta}^{*})$
        \State $\tilde{\vb{y}} \gets \kappa(\vb*{t}',\vb*{t};\vb*{\theta}^{*})(\vb{K}^{yy})^{-1}\vb{y}$
            \label{step:ytilde}
            \Comment{GP state estimate over $\vb*{t}'$.}
        \State $\tilde{\vb{z}} \gets \vb{K}^{zy} (\vb{K}^{yy})^{-1}\vb{y}$
            \Comment{GP time derivative estimate over $\vb*{t}'$.}
        \State $\vb{C} \gets \vb{K}^{zy}(\vb{K}^{yy})^{-1}(\vb{K}^{zy})\trp$
            \label{step:zyyyzy}
        \State $\vb*{\Psi}\vb*{\Xi}\vb*{\Psi}\trp \gets \vb{K}^{zz}
        - \frac{1}{2}\left(\vb{C} + \vb{C}\trp\right)
        + \tau\vb{I}_{m'}$
            \Comment{Eigendecomposition of $(\vb{W}^{zz})^{-1}$ (see \Cref{remark:gpreg}).}
            \label{step:eigendecomp}
        \State $(\vb{W}^{zz})^{1/2} \gets \vb*{\Psi}\vb*{\Xi}^{-1/2}\vb*{\Psi}\trp$
        \State \textbf{return} $\tilde{\vb{y}}$, $\tilde{\vb{z}}$, $(\vb{W}^{zz})^{1/2}$
            \Comment{State / time derivative estimates and weighting matrix.}
    \EndProcedure
    \end{algorithmic}
    \caption{Gaussian process regression with derivative estimation.}
    \label{alg:GP-Fit}
\end{algorithm}

The data matrix $\tilde{\vb{D}} \in \mathbb{R}^{m' \times d(r, p)}$ driving the Bayesian inference contains data at the estimation points $\vb*{t}'$, yet we only are provided with snapshot data corresponding to $\vb*{t}$.
To populate $\tilde{\vb{D}}$, we utilize the one-dimensional GP $\hat{q}_i$ of \cref{eq:gp-single} to estimate the state at times $\vb*{t}'$,
\begin{equation}
    \label{eq:gp-state-estimate}
    \tilde{\vb{y}}_{i}
    := \kappa_i(\vb*{t}',\vb*{t})\left(\kappa_i(\vb*{t},\vb*{t})+\chi_i\vb{I}_{m}\right)^{-1}\vb{y}_{i}
    = \kappa_i(\vb*{t}',\vb*{t})(\vb{K}_{i}^{yy})^{-1}\vb{y}_{i}\,,
\end{equation}
which provides a smooth interpolation of the trajectory.
This is especially beneficial for cases when the snapshot data $\vb{y}_1,\ldots,\vb{y}_r$ are noisy.
Collecting the filtered state estimates $\tilde{\vb{y}}_1,\ldots,\tilde{\vb{y}}_r$, we construct $\tilde{\vb{D}}$ akin to~\cref{eq:opinf-data-matrix}:
\begin{align}
    \label{eq:filtered-data-matrix}
    \tilde{\vb{D}}
    &= \left[\begin{array}{c}
        \vb{d}(\tilde{\vb{q}}_0, \vb{u}(t_0'))\trp \\
        \vb{d}(\tilde{\vb{q}}_1, \vb{u}(t_1'))\trp \\
        \vdots \\
        \vb{d}(\tilde{\vb{q}}_{m'-1}, \vb{u}(t_{m'  -1}'))\trp \\
    \end{array}\right]
    \in\mathbb{R}^{m' \times d(r, p)}\,,
    &
    \left[\begin{array}{cccc}
        & & & \\
        \tilde{\vb{q}}_{0} &
        \tilde{\vb{q}}_{1}
        & \cdots &
        \tilde{\vb{q}}_{m'-1}
        \\ & & &
    \end{array}\right]
&= \left[\begin{array}{c}
        \tilde{\vb{y}}_{1}\trp \\
        \vdots \\
        \tilde{\vb{y}}_{r}\trp
    \end{array}\right]
    \in \mathbb{R}^{r\times m'}.
\end{align}
From \cref{eq:gp-state-estimate}, we explicitly see the derivative relationship between the state estimates $\tilde{\vb{y}}_i$ and the time derivative estimates $\tilde{\vb{z}}_i$:
\begin{equation}
    \label{eq:gp-timederivative}
    \partial_1\tilde{\vb{y}}_i
    = \partial_1\kappa_i(\vb*{t}',\vb*{t})(\vb{K}_{i}^{yy})^{-1}\vb{y}_{i}
    = \vb{K}^{zy}_i (\vb{K}^{yy}_i)^{-1}\vb{y}_i
    = \tilde{\vb{z}}_i \,.
\end{equation}

\Cref{alg:GP-Fit} enumerates the steps for the GP~regression for a single reduced mode.
Apart from solving the optimization problem \cref{eq:maxmlhd}, the computational cost is $\mathcal{O}(m^3 + (m')^3)$, dominated by the expense of inverting $\vb{K}^{yy}$ and factorizing $(\vb{W}^{zz})^{-1}$.
These matrices are symmetric, hence we use Cholesky factorization and back substitution to apply $(\vb{K}^{yy})^{-1}$ in steps \ref{step:ytilde}--\ref{step:zyyyzy}
and use symmetric QR for the eigendecomposition of $(\vb{W}^{zz})^{-1}$ in step \ref{step:eigendecomp},
both via standard tools available in \texttt{scipy.linalg} \cite{2020SciPy-NMeth}.

\subsection{Prior variance selection}
\label{sec:regularization}

\begin{algorithm}[t]
\begin{algorithmic}[1]
\Procedure{OpInfError}{
    \newline\phantom{---}Prior variances $\vb*{\gamma}_1,\ldots,\vb*{\gamma}_r \in \mathbb{R}^{d}$;
    \newline\phantom{---}Regression arguments $\vb*{\alpha} = \{\tilde{\vb{D}},\tilde{\vb{z}}_{1},\ldots,\tilde{\vb{z}}_{r},(\vb{W}_{1}^{zz})^{1/2},\ldots,(\vb{W}_{r}^{zz})^{1/2}\}$,
\newline\phantom{---}GP state estimates $\tilde{\vb{Q}}\in\mathbb{R}^{r\times m'}$
    \newline\phantom{---}Input function $\vb{u}:\mathbb{R}\to\mathbb{R}^{p}$,
    \newline\phantom{---}GP estimation times $\vb*{t}' \in\mathbb{R}^{m'}$,
    \newline\phantom{---}Stability margin parameter $\varphi > 0$,
    \newline\phantom{---}Number of samples $n_s \in \mathbb{N}$,
    \newline\phantom{---}Final prediction time $t_f > t'_0$
    \newline}

    \State $(\vb*{\mu}_1,\vb*{\Sigma}_1),\ldots,(\vb*{\mu}_r,\vb*{\Sigma}_r) \gets
        \textproc{OpPost(}
            \vb*{\alpha},
\vb*{\gamma}_1,\ldots,\vb*{\gamma}_r
        \textproc{)}$
        \Comment{\Cref{alg:OperatorPosterior}.}
\State $\hat{\vb{q}}_0 \gets \tilde{\vb{Q}}_{:,0}$
        \Comment{Initial condition.}
    \State $B \gets \varphi \cdot \max_{ij}|\tilde{\vb{Q}}_{ij}|$
        \Comment{Bound for stability check.}
    \For{$k = 1, \ldots, n_s$}
        \State $\hat{\vb{O}}_{k} \gets~\text{sample from}~~[~\mathcal{N}(\vb*{\mu}_1,\vb*{\Sigma}_1)~~\cdots~~\mathcal{N}(\vb*{\mu}_r,\vb*{\Sigma}_r)~]\trp$
            \Comment{Sample operator matrix.}
        \State $\breve{\vb{Q}} \gets~\text{solve}~~\frac{\textrm{d}}{\textrm{d}t}\hat{\vb{q}}(t) = \hat{\vb{O}}_{k}\vb{d}(\hat{\vb{q}}(t),\vb{u}(t)),~\hat{\vb{q}}(t'_0) = \hat{\vb{q}}_0~~\text{over}~~t\in[t'_0,t_f]$
        \If{$\max_{ij}|\breve{\vb{Q}}_{ij}| > B$}
            \Comment{Detect and disqualify unstable ROMs.}
            \State \textbf{return} $\infty$
        \EndIf
        \State $\hat{\vb{Q}}_{k} \gets~\text{solve}~~\frac{\textrm{d}}{\textrm{d}t}\hat{\vb{q}}(t) = \hat{\vb{O}}_{k}\vb{d}(\hat{\vb{q}}(t),\vb{u}(t)),~\hat{\vb{q}}(t'_0) = \hat{\vb{q}}_0~~\text{over}~~t\in\vb*{t}'$
    \EndFor
    \State $\hat{\vb{Q}} = \frac{1}{n_s}\sum_{k=1}^{n_s}\hat{\vb{Q}}_{k}$
        \Comment{Sample mean of the ROM solutions.}
    \State \textbf{return} $\|\tilde{\vb{Q}} - \hat{\vb{Q}}\|$
        \Comment{Error compared to GP state estimates.}
\EndProcedure
\end{algorithmic}
\caption{Estimate the mean error of the Bayesian OpInf model for a given prior variance.}
\label{alg:OpInfError}
\end{algorithm}

The vector $\vb*{\gamma}_i\in\mathbb{R}^{d}$ parameterizing the prior distribution $\pi_\texttt{prior}(\hat{\vb{o}}_{i})$ plays the role of a Tikhonov regularizer in the generalized least-squares regression \cref{eq:generalizedlstsq} and strongly affects the stability of the resulting ROM.
To methodically select $\vb*{\gamma}_1,\ldots,\vb*{\gamma}_r$, we extend the approach introduced in \cite{mcquarrie2021combustion}, which balances model stability with accuracy in reproducing training data.
For a candidate set $\{\vb*{\gamma}_i\}_{i=1}^{r}$, we solve the corresponding Bayesian inference \cref{eq:posterior-mean}--\cref{eq:posterior-covariance} for $\{(\vb*{\mu}_i,\vb*{\Sigma}_i)\}_{i=1}^{r}$, sample $n_s\in\mathbb{N}$ times from the posterior operator distribution $\hat{\vb{O}}\sim[~\mathcal{N}(\vb*{\mu}_1,\vb*{\Sigma}_1)~~\cdots~~\mathcal{N}(\vb*{\mu}_r,\vb*{\Sigma}_r)~]\trp$, and solve the resulting reduced systems \cref{eq:reducedorder} over the full time domain of interest~$[t_0,t_f]$.
Recall that $t_f > t_{m-1}$, that is, the full simulation domain extends beyond where noisy data are observed.
If the magnitude of any of the reduced model solutions exceeds a predetermined bound computed from the training data, the candidate set is disqualified for failing to produce a sufficiently stable model; otherwise, the sample mean of the model solutions is compared to the filtered state data $\tilde{\vb{y}}_1,\ldots,\tilde{\vb{y}}_r$.
A grid search, followed by a derivative-free optimization, is used to select optimal $\{\vb*{\gamma}_i^{*}\}_{i=1}^{r}$ that minimizes the difference between the sample mean and the training data while not violating the stability criterion.
Details are provided in \Cref{alg:OpInfError}; in \Cref{sec:results}, the experiments use $n_s = 20$ posterior samples and set the stability margin constant to $\varphi = 5$.
A similar strategy was deployed in our previous work \cite{guo2022bayesopinf}; however, in that work, only the MAP points $\hat{\vb{O}} = [~\vb*{\mu}_1~~\cdots~~\vb*{\mu}_r~]\trp$ were used to evaluate candidate regularizers.
Our approach, only slightly more computationally expensive, tests stability for multiple posterior samples and hence produces more robust models than the previous approach.

Executing \Cref{alg:OpInfError} is computationally inexpensive because all operations are in the reduced state space.
However, minimizing the reconstruction error with respect to the prior variance vectors $\vb*{\gamma}_1,\ldots,\vb*{\gamma}_r\in\mathbb{R}^{d(r,p)}$ constitutes an optimization of $r\cdot d(r,p)$ positive variables, which may require a large number of evaluations of \Cref{alg:OpInfError}.
To reduce the cost of the optimization, it is beneficial to parameterize the prior variances and thus reduce the total degrees of freedom designating the prior.
First, setting $\vb*{\gamma}_1 = \cdots = \vb*{\gamma}_r = \vb*{\gamma} \in \mathbb{R}^{d(r, p)}$ reduces the number of unknowns to $d(r, p)$ positive constants.
This is reasonable as long as the reduced-order modes have approximately the same magnitude.
Second, $\vb*{\gamma}$ may be designed so that each operator in the model is penalized individually.
For example, for a quadratic ROM with operator matrix $\hat{\vb{O}} = [~\hat{\vb{H}}~~\hat{\vb{B}}~]$, we may define
\begin{align}
    \label{eq:regsplit}
    \vb*{\gamma} = [~
        \overbrace{\gamma_1~~\cdots~~\gamma_1}^{r^2~\text{times}}
        ~~
        \overbrace{\gamma_2~~\cdots~~\gamma_2}^{p~\text{times}}~~]\trp\,,
\end{align}
so that $\gamma_1 > 0$ corresponds to a penalization for the quadratic operator $\hat{\vb{H}}\in\mathbb{R}^{r\times r^2}$ and $\gamma_2 > 0$ penalizes the input operator $\hat{\vb{B}}\in\mathbb{R}^{r\times p}$.
This choice is similar to the regularization strategies of \cite{jain2021performance,mcquarrie2021combustion,qian2022reduced,sawant2023piregopinf} and results in a two-dimensional optimization, reasonably solvable via, e.g., the Nelder--Mead method \cite{nelder1965simplex}.
In our numerical experiments, we found it convenient to make the further simplification $\vb*{\gamma} = [~\gamma~~\cdots~~\gamma~]\trp$ for a single constant $\gamma > 0$, resulting in a one-dimensional optimization that can be carried out, e.g., via Brent's method \cite{Brent2002minimization}.
The numerical results were comparable when using the strategy of \cref{eq:regsplit} or allowing each entry of~$\vb*{\gamma}$ to vary.

\begin{algorithm}[t]
    \begin{algorithmic}[1]
    \Procedure{GPBayesOpInf}{
    \newline\phantom{---}
        Observation times $\vb*{t}\in\mathbb{R}^{m}$,
    \newline\phantom{---}
        Snapshots $\vb{Q}\in\mathbb{R}^{N \times m}$,
        \newline\phantom{---}
        Reduced dimension $r \ll N$,
    \newline\phantom{---}
        Estimation times $\vb*{t}'\in\mathbb{R}^{m'}$,
        \newline\phantom{---}
        Numerical regularization $0 < \tau \ll 1$,
        \newline\phantom{---}
        Input function $\vb{u}:\mathbb{R}\to\mathbb{R}^{p}$
        \newline
        }
    \LineComment{GP regression on reduced training data.}
        \State $\vb{V} \gets \texttt{basis}(\vb{Q}; r)$
            \label{step:basis}
            \Comment{Rank-$r$ basis (e.g., POD).}
        \State $[~\vb{y}_1~~\cdots~~\vb{y}_r~]\trp \gets \vb{V}\trp\vb{Q}$
            \label{step:compression}
            \Comment{Compress training data.}
        \For{$i = 1, \ldots, r$}
            \State $\tilde{\vb{y}}_i,\tilde{\vb{z}}_i,(\vb{W}^{zz}_i)^{1/2} \gets \textproc{GPFit(}\vb*{t},\vb{y}_i,\vb*{t}',\tau\textproc{)}$
                \Comment{\Cref{alg:GP-Fit}.}
        \EndFor

        \LineComment{OpInf regression to construct a probabilistic ROM.}
        \State $\tilde{\vb{Q}} \gets [~\tilde{\vb{y}}_1~~\cdots~~\tilde{\vb{y}}_r~]\trp$
            \label{step:Qtilde}
        \State $\tilde{\vb{D}} \gets \texttt{data\_matrix}(\tilde{\vb{Q}},\vb{u}(\vb*{t}'))$
            \Comment{Data matrix \cref{eq:filtered-data-matrix} from GP state estimates.}
        \State $\vb*{\alpha} \gets \{\tilde{\vb{D}},\tilde{\vb{z}}_{1},\ldots,\tilde{\vb{z}}_{r},(\vb{W}_{1}^{zz})^{1/2},\ldots,(\vb{W}_{r}^{zz})^{1/2}\}$
            \label{step:arguments}
            \Comment{Regression arguments.}
        \State $\vb*{\gamma}_{1}^{*},\ldots,\vb*{\gamma}_{r}^{*} \,\gets\,
    \operatorname{argmin}_{\vb*{\gamma}_{1},\ldots,\vb*{\gamma}_{r}}
            \textproc{OpInfError}(\vb*{\gamma}_{1},\ldots,\vb*{\gamma}_{r}; \vb*{\alpha},\tilde{\vb{Q}},\vb{u},\vb*{t}')$
            \Comment{\Cref{alg:OpInfError}.}
        \State $(\vb*{\mu}_1,\vb*{\Sigma}_1),\ldots,(\vb*{\mu}_r,\vb*{\Sigma}_r) \,\gets\, \textproc{OpPost(}\vb*{\alpha},\vb*{\gamma}_1^{*},\ldots,\vb*{\gamma}_r^{*}\textproc{)}$
            \Comment{\Cref{alg:OperatorPosterior}.}
        \State \textbf{return} $[~\mathcal{N}(\vb*{\mu}_1,\vb*{\Sigma}_1)~~\cdots~~\mathcal{N}(\vb*{\mu}_r,\vb*{\Sigma}_r)~]\trp$
    \EndProcedure
    \end{algorithmic}
    \caption{Bayesian Operator Inference with Gaussian processes (GP-BayesOpInf).}
    \label{alg:GP-BayesOpInf}
\end{algorithm}

\begin{figure}[t]
    \centering
    \includegraphics[width=\textwidth]{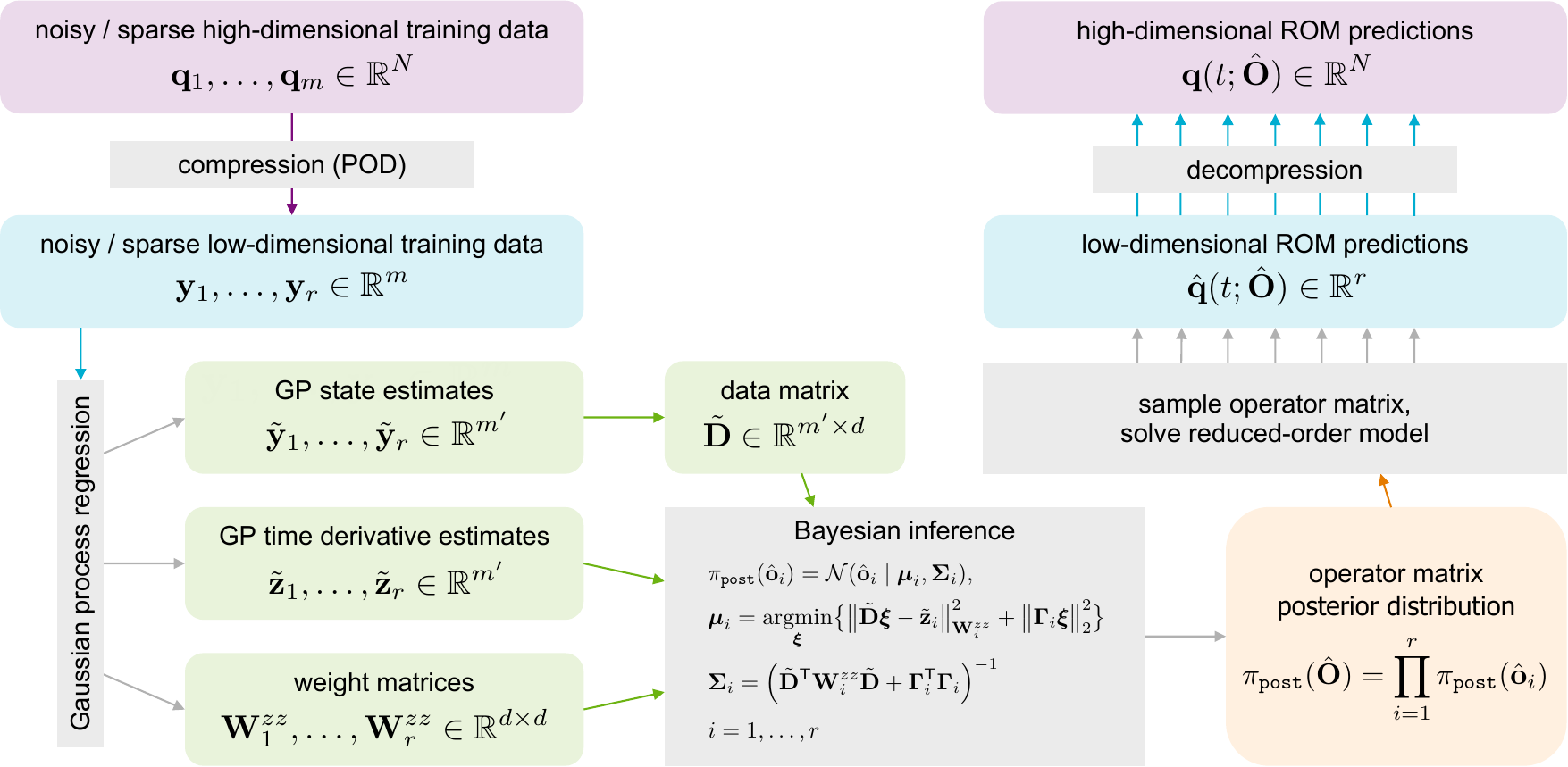}
    \vspace{.05cm}
    \caption{Flow diagram of GP-BayesOpInf via \Cref{alg:GP-BayesOpInf}.}
    \label{fig:GP-BayesOpInf-diagram}
\end{figure}

With all of the major ingredients now in hand, the entire computational procedure is detailed in \Cref{alg:GP-BayesOpInf} and illustrated in \Cref{fig:GP-BayesOpInf-diagram}.
To ensure that all reduced-order modes are of similar magnitudes, we center the (noisy) snapshots $\vb{q}_{0},\ldots,\vb{q}_{m-1}\in\mathbb{R}^{N}$ by setting
$
    \vb{Q}
    = [~
        (\vb{q}_0 - \bar{\vb{q}})
        ~~\cdots~~
        (\vb{q}_{m-1} - \bar{\vb{q}})
    ~]
    \in \mathbb{R}^{N\times m},
$
where $\bar{\vb{q}} = \frac{1}{m}\sum_{j=0}^{m-1}\vb{q}_{j} \in \mathbb{R}^{N}$ is the mean snapshot.
These centered snapshots are compressed to $r \ll N$ degrees of freedom in steps~\ref{step:basis}--\ref{step:compression} by introducing the $N\times r$ basis matrix $\vb{V}$.
We use proper orthogonal decomposition~(POD)~\cite{Berkooz1993,GPT1999vortexsheddingPOD,sirovich1987turbulence}, in which $\vb{V}$ is comprised of the first $r$ left singular vectors of $\vb{Q}$.
This choice is advantageous because POD has a hierarchical smoothing effect \cite{del2008proper,epps2010PODthresh,venturi2006PODperturb}, the dominant modes being smoothed more than the high-frequency modes (see \Cref{fig:euler-dims}).
After compression, $r$ GPs are trained via \Cref{alg:GP-Fit}, which are used in steps \ref{step:Qtilde}--\ref{step:arguments} to populate the generalized least-squares regression \cref{eq:generalizedlstsq} for each reduced mode.
Finally, the prior variances are chosen by minimizing the procedure of \Cref{alg:OpInfError}, and the posterior for the operator matrix is produced by applying \Cref{alg:OperatorPosterior} a final time.
With only minor adjustments, Algorithms \ref{alg:OpInfError}--\ref{alg:GP-BayesOpInf} can be adapted
to cases where multiple data trajectories are available, as discussed in \Cref{sec:multitrajectory} (see \apdxref{appendix:multiple-trajectories}).

\section{Numerical examples of GP-BayesOpInf}
\label{sec:results}

This section demonstrates GP-BayesOpInf via \Cref{alg:GP-BayesOpInf} for two nonlinear PDEs.
The first example uses a single data trajectory while the second incorporates multiple trajectories.
To guarantee sufficient smoothness for derivative evaluations, in each of the examples we use the squared exponential (radial basis function) kernel
\begin{equation}
    \kappa(t,t';\vb*{\theta})
    = \sigma^{2}\exp\left(-\frac{1}{2}\frac{(t-t')^2}{\ell^2}\right)\,,
\end{equation}
in which the signal variance $\sigma^2 > 0$ and the lengthscale $\ell > 0$ are the kernel hyperparameters.
Codes to reproduce the experiments are available at \url{https://github.com/sandialabs/GP-BayesOpInf}.

\subsection{Compressible Euler equations}

In this section, we show the effectiveness of the GP-BayesOpInf method as a probabilistic surrogate for compressible Euler equations with noisy and sparse data.

\subsubsection{Governing equations}

We consider the compressible Euler equations for an ideal gas over the one-dimensional spatial domain $[0, 2)$ with periodic boundary conditions, written in conservative form as
\begin{subequations}
\label{eq:euler-all}
\begin{align}
    \label{eq:euler-conservative}
    \frac{\partial}{\partial t}\left[\begin{array}{c}
        \rho \\ \rho v \\ \rho e
    \end{array}\right]
    &= -\frac{\partial}{\partial x}\left[\begin{array}{c}
        \rho v \\ \rho v^2 + p \\ (\rho e + p)v
    \end{array}\right],
    \\
    \label{eq:idealgaslaw}
    \rho e
    &= \frac{p}{\gamma - 1} + \frac{1}{2}\rho v^2,
    \\
    \label{eq:periodicBCs}
    \rho(0, t) = \rho(2, t),
    \quad
    v(0, t) &= v(2, t),
    \quad
    e(0, t) = e(2, t),
\end{align}
\end{subequations}
where
$v$ is the velocity, $\rho$ is the density, $p$ is the pressure, $\rho v$ is the specific momentum, $\rho e$ is the total energy, and $\gamma = 1.4$ is the heat capacity ratio.
The initial conditions are defined by setting $p(x, 0) = 10^5~\text{Pa}$ and constructing cubic splines for the density and velocity with interpolation values
\begin{align}
    \begin{aligned}
    \rho(0, 0) &= 22~\text{kg/m}^3,
    &
    \rho(\tfrac{2}{3}, 0) &= 20~\text{kg/m}^3,
    &
    \rho(\tfrac{4}{3}, 0) &= 24~\text{kg/m}^3,
    \\
    v(0, 0) &= 95~\text{m/s},
    &
    v(\tfrac{2}{3}, 0) &= 105~\text{m/s},
    &
    v(\tfrac{4}{3}, 0) &= 100~\text{m/s}.
    \end{aligned}
\end{align}
The initial energy $e(x,0)$ is defined from the initial conditions of $p$, $\rho$, and $v$ through the ideal gas law~\cref{eq:idealgaslaw}.

\begin{figure}[t]
    \centering
    \includegraphics[width=\textwidth]{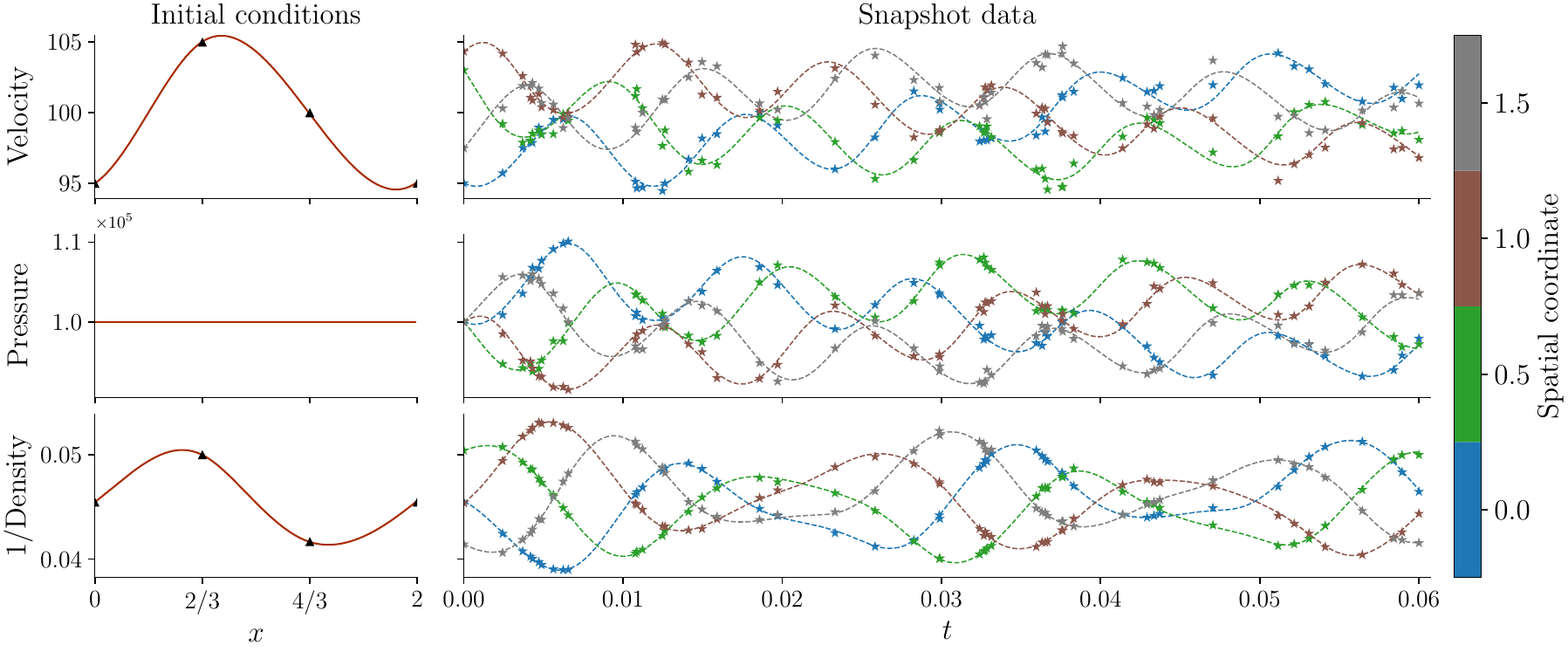}
    \vspace{-0.75cm}
    \caption{Initial conditions for the Euler system \cref{eq:euler-conservative}--\cref{eq:periodicBCs} (left) and noisy state observations (right) expressed in the specific volume variables $(v,p,\zeta)$. The initial conditions are constructed with periodic splines interpolating the marked points. On the right, the dashed lines show the numerical solution without noise as a function of time at four fixed points in space; the markers indicate $m = 50$ sparsely observed states with $\xi = 1\%$ relative noise. Though data are only shown here at select spatial locations, they are observed at every point in the spatial discretization.}
    \label{fig:euler-fomdata}
\end{figure}

To generate noisy and/or sparse data for this problem, we discretize the spatial derivatives in \cref{eq:euler-conservative} with first-order rightward finite differences using $n_x = 200$ degrees of freedom per variable, resulting in a system of $N = 3 n_x = 600$ ordinary differential equations.
For a given number of observations $m\in\mathbb{N}$ and a final observation time $t_{m-1} > 0$, we select observation times $\vb*{t} = [~t_0~~t_1~~\cdots~~t_{m-1}~]\trp\in\mathbb{R}^{m}$ by setting $t_0 = 0$ and sampling the $m - 2$ additional points $(t_1,\ldots,t_{m-2})$ from the uniform distribution over the interval~$(0, t_{m-1})$.
Note that this results in non-uniform spacing between the sampling times.
We then integrate the discretized system in time with the explicit adaptive-step Runge--Kutta method of order 5~\cite{dormand1980rk45} and employ quartic polynomial interpolation~\cite{shampine1986rkinterp} to obtain the solution at times $\vb*{t}$.
The data after the initial condition are polluted with additive Gaussian noise.
Specifically, for a given noise level $\xi \ge 0$, we observe

\begin{subequations}
\begin{align}
    \vb{q}^{(c)}_j
    := \left[\begin{array}{c}
        \vb{q}_{\rho}(t_j) + \vb*{\varepsilon}_{\rho} \\
        \vb{q}_{\rho v}(t_j) + \vb*{\varepsilon}_{\rho v} \\
        \vb{q}_{\rho e}(t_j) + \vb*{\varepsilon}_{\rho e} \end{array}\right]
    \in \mathbb{R}^{N},
    \quad
    \begin{aligned}
        \vb*{\varepsilon}_{\rho} \sim \mathcal{N}(\vb{0}, \nu_{\rho}^{2}\vb{I}_{n_x}),
        \\
        \vb*{\varepsilon}_{\rho v} \sim \mathcal{N}(\vb{0}, \nu_{\rho v}^{2}\vb{I}_{n_x}),
        \\
        \vb*{\varepsilon}_{\rho e} \sim \mathcal{N}(\vb{0}, \nu_{\rho e}^{2}\vb{I}_{n_x}),
    \end{aligned}
    \quad
    j = 1, \ldots, m-1,
\end{align}
where $\vb{q}_{\rho}(t_j),\vb{q}_{\rho v}(t_j),\vb{q}_{\rho e}(t_j)\in\mathbb{R}^{n_x}$ respectively denote the spatially discrete numerical solutions for $\rho$, $\rho v$, and $\rho e$ at time $t_j\in\vb*{t}$, and where $\nu_{\rho},\nu_{\rho v},\nu_{\rho e} > 0$ is the observed range of the indicated variable scaled by the noise level.
For example, collecting all numerical solutions for the density in the matrix $\vb{Q}_\rho = [~\vb{q}_\rho(t_0)~~\cdots\vb{q}_\rho(t_{m-1})~]$, we set $\nu_\rho = \xi(\max(\vb{Q}_{\rho}) - \min(\vb{Q}_{\rho}))$.
We also observe the given initial conditions without noise,
\begin{align}
    \vb{q}_{0}^{(c)}
    := \left[\begin{array}{c}
        \vb{q}_{\rho}(t_0) \\ \vb{q}_{\rho v}(t_0) \\ \vb{q}_{\rho e}(t_0)
    \end{array}\right] \in \mathbb{R}^{N}.
\end{align}
\end{subequations}
\Cref{fig:euler-fomdata} shows the initial conditions and, at select spatial locations, $m = 50$ snapshots with $t_{m-1} = 0.06~\text{s}$ and $\xi = 1\%$ additive noise, along with the nonnoised numerical solution to the system.

\subsubsection{ROM structure}
Before applying \Cref{alg:GP-BayesOpInf}, we apply a variable transformation to cast \cref{eq:euler-conservative} in a polynomial form that is suitable for OpInf~\cite{QKPW2020liftAndLearn}.
Using the ideal gas law \cref{eq:idealgaslaw}, changing from the conservative variables~$(\rho, \rho v, \rho e)$ to the specific volume variables $(v, p, 1/\rho)$ results in a quadratic system,
\begin{align}
    \label{eq:euler-lifted}
    \frac{\partial}{\partial t}\left[\begin{array}{c}
        v \\ p \\ \zeta
    \end{array}\right]
    = \left[\begin{array}{c}
        -v \frac{\partial v}{\partial x} - \zeta\frac{\partial p}{\partial x} \\
        -\gamma p \frac{\partial v}{\partial x} - v\frac{\partial p}{\partial x} \\
        -v \frac{\partial\zeta}{\partial x} + \zeta\frac{\partial v}{\partial x}
    \end{array}\right],
\end{align}
where $\zeta = 1/\rho$ is the specific volume.
Note that the variables $(v,p,\zeta)$ obey periodic boundary conditions as in~\cref{eq:periodicBCs}.
For each $j = 0,\ldots,m-1$, the noisy conservative-variable snapshot $\vb{q}^{(c)}_j$ is transformed to the specific volume variables and nondimensionalized
to produce the snapshot $\vb{q}_{j} \in \mathbb{R}^{N}$.
After spatial discretization, the system \cref{eq:euler-lifted} can be written in the quadratic form
\begin{align}
    \frac{\text{d}}{\text{d}t}\vb{q}(t)
    = \vb{H}[\vb{q}(t)\otimes\vb{q}(t)],
\end{align}
where $\vb{H}\in\mathbb{R}^{N\times N^{2}}$ and $\vb{q}(t) \in \mathbb{R}^{N}$ is the spatial discretization of the lifted variables $(v,p,\zeta)$ at time $t$.
Using a centered-POD approximation,
\begin{align}
    \label{eq:PODcentered}
    \vb{q}(t)
    \approx
    \vb{V}\hat{\vb{q}}(t) +
    \bar{\vb{q}}, \qquad \bar{\vb{q}} =
    \frac{1}{m}\sum_{j=0}^{m-1} \vb{q}_{j}
    \in\mathbb{R}^{N},
\end{align}
where $\vb{V}\in\mathbb{R}^{N\times r}$ has orthonormal columns and $\hat{\vb{q}}(t) \in \mathbb{R}^{r}$,
we obtain a ROM
\begin{align}
    \label{eq:euler-rom}
    \frac{\textrm{d}}{\textrm{d}t}\hat{\vb{q}}(t)
    = \hat{\vb{c}} + \hat{\vb{A}}\hat{\vb{q}}(t) + \hat{\vb{H}}[\hat{\vb{q}}(t)\otimes\hat{\vb{q}}(t)],
\end{align}
\newpage

\begin{figure}[!ht]
    \centering
    \includegraphics[width=\textwidth]{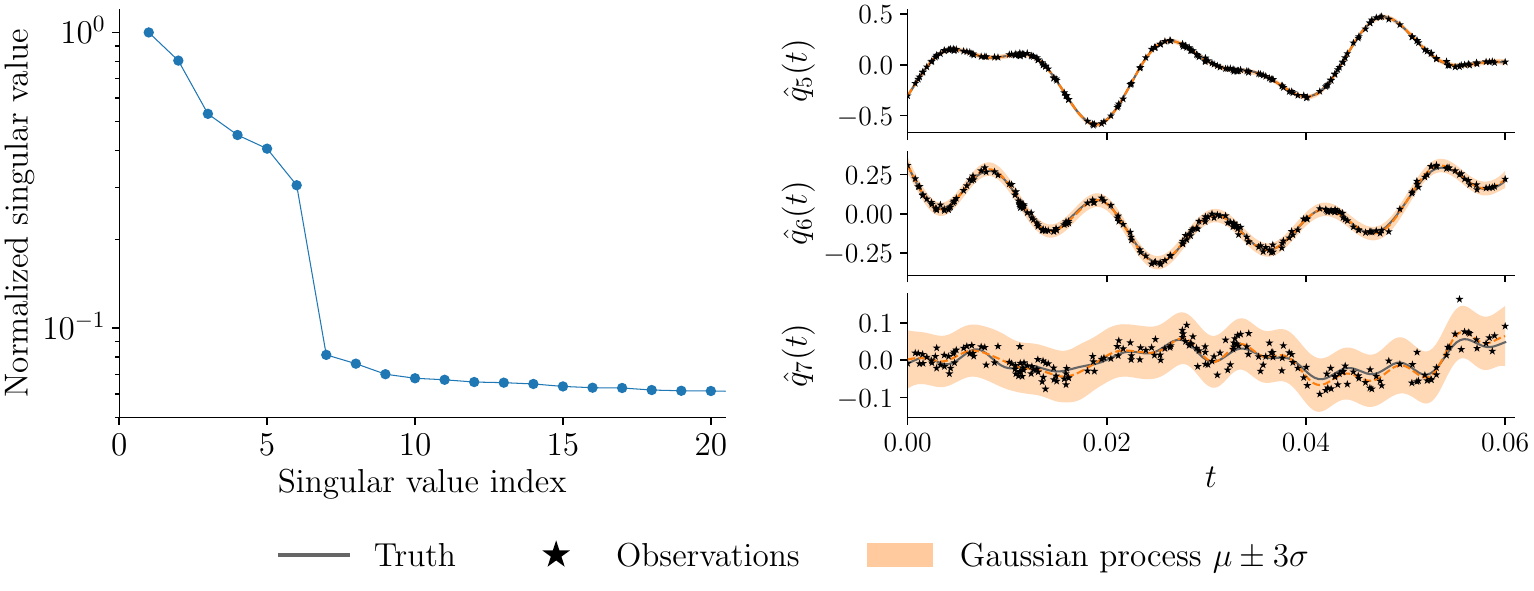}
    \vspace{-0.75cm}
    \caption{Singular value decay of $m=200$ snapshots with $\xi=3\%$ relative noise (left) and trained GPs for reduced modes $5$, $6$, and $7$, compressed via centered POD (right). The shaded regions are within three standard deviations of each GP mean. The POD compression filters much of the noise in the dominant modes, hence each GP exhibits a tight fit to the data until the spectral gap after $r=6$ modes.}
    \label{fig:euler-dims}
\end{figure}
\begin{figure}[!ht]
    \centering
    \includegraphics[width=\textwidth]{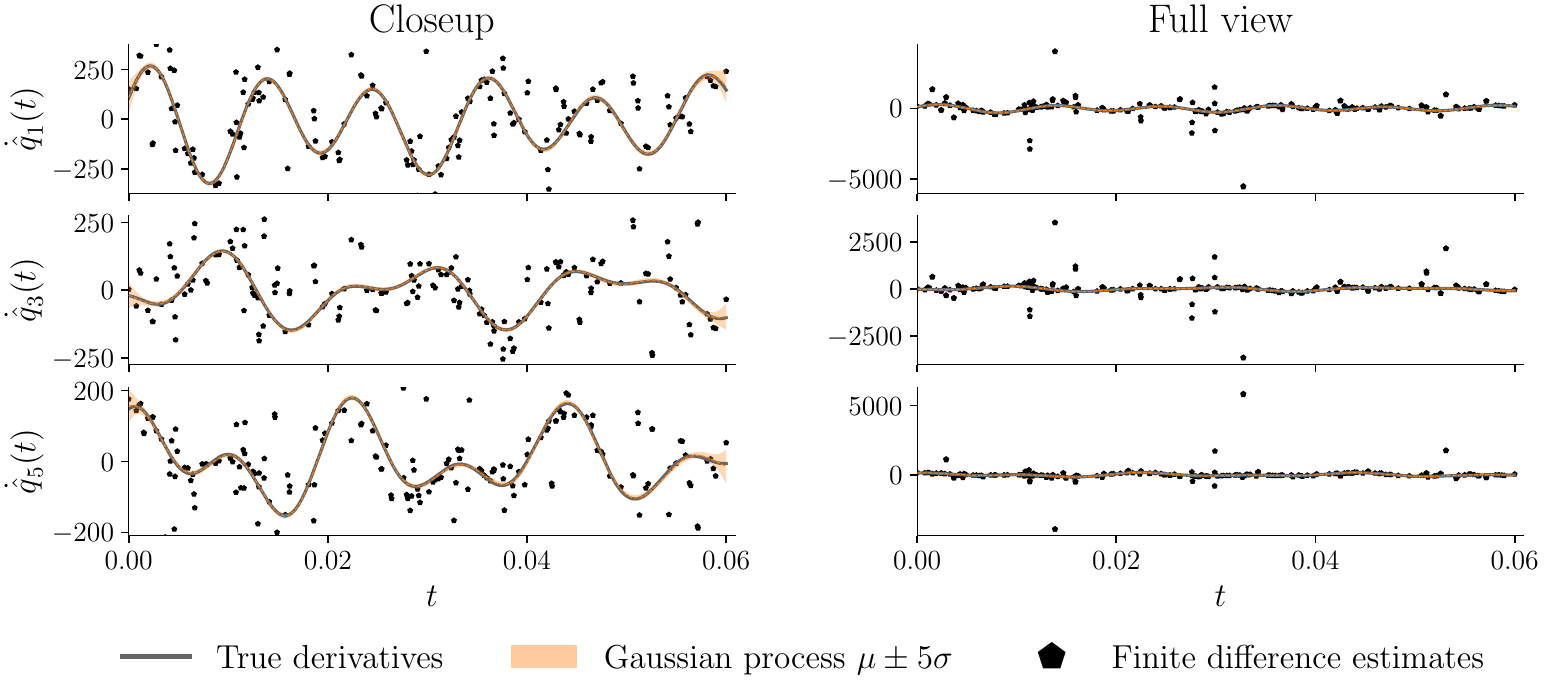}
    \vspace{-0.75cm}
    \caption{True reduced state time derivatives, as well as GP and finite difference estimates of the time derivatives, using $m = 200$ observations with $\xi = 3\%$ relative noise. The time derivative GPs have small variances and produce highly accurate estimates~(left), unlike the finite difference estimates, which in some instances differ from the truth by an order of magnitude~(right).}
    \label{fig:euler-ddts}
\end{figure}

\noindent
which is \cref{eq:reducedorder} with
$\hat{\vb{O}} = [~\vb{\hat{c}}~~\vb{\hat{A}}~~\vb{\hat{H}}~]\in\mathbb{R}^{r\times d(r)}$,
$\vb{d}(\hat{\vb{q}}) = [~1~~\hat{\vb{q}}\trp~~(\hat{\vb{q}}\otimes\hat{\vb{q}})\trp~]\trp \in \mathbb{R}^{d(r)}$,
and $d(r) = 1 + r + r^2$.
The nonquadratic terms $\hat{\vb{c}}\in\mathbb{R}^{r}$ and $\hat{\vb{A}}\in\mathbb{R}^{r\times r}$ appear due to the centering term $\bar{\vb{q}}$ \cite{swischuk2019thesis}.
Because the Kronecker product $\hat{\vb{q}}\otimes \hat{\vb{q}} = \text{vec}(\hat{\vb{q}}\hat{\vb{q}}\trp)\in \mathbb{R}^{r^2}$
has only $r(r+1)/2 < r^2$ degrees of freedom due to symmetry, in practice the quadratic operator $\hat{\vb{H}}\in\mathbb{R}^{r\times r^{2}}$ can be fully represented by an $r\times r(r+1)/2$ matrix that acts on a compressed version of the Kronecker product (see, e.g., Appendix~B of \cite{mcquarrie2023thesis}).

\begin{figure}[!ht]
    \centering
    \includegraphics[width=\textwidth]{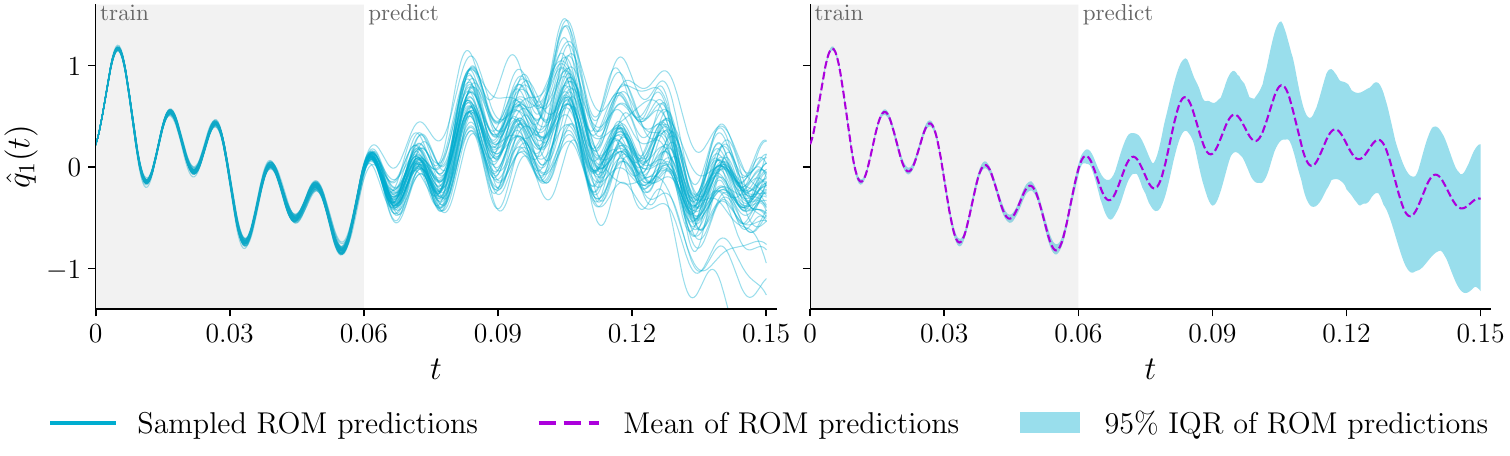}
    \vspace{-0.75cm}
    \caption{The dominant mode of the Bayesian ROM predictions for the Euler system \cref{eq:euler-conservative}--\cref{eq:periodicBCs} using $m=200$ observations over the time domain $[0, 0.06]$ with $\xi=3\%$ relative noise. Fifty samples are taken of the posterior operator distribution and the corresponding ROM predictions are displayed as a function of time (left). The uncertainty of the prediction distribution is characterized by computing the mean and the $95\%$ interquantile range (IQR) of the sample predictions at each point in time~(right).} \label{fig:euler-draws}
\end{figure}

\begin{figure}[!ht]
    \centering
    \includegraphics[width=\textwidth]{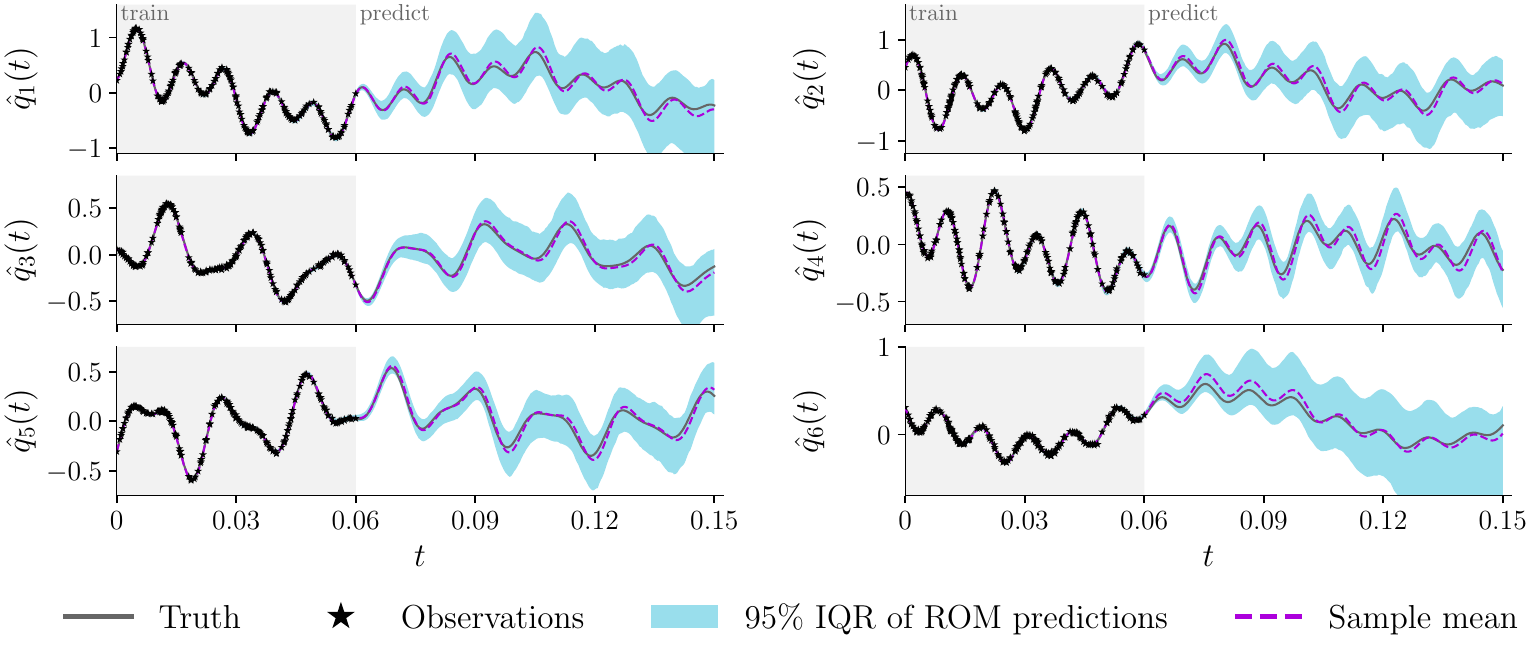}
    \vspace{-0.75cm}
    \caption{Bayesian ROM predictions for the Euler system \cref{eq:euler-conservative}--\cref{eq:periodicBCs} using $m=200$ observations over the time domain $[0, 0.06]$ with $\xi=3\%$ relative noise. The solid line shows the true numerical solution to the system without noise, expressed in the POD coordinates. The posterior operator distribution is sampled $500$ times.
    The pointwise mean of the sample predictions gives an accurate estimate of the true solution, with uncertainty described by the $95\%$ interquantile range.}
    \label{fig:euler-noisy-reduced}
\end{figure}

\subsubsection{Noisy training data}

We consider first a scenario in which the observed data are relatively abundant but corrupted with significant noise.
Let $m = 200$ with final observation time $t_{m-1} = 0.06~\text{s}$ and relative noise level $\xi = 3\%$.
This is a nontrivial amount of noise for this problem, especially given that the noise is applied to the conservative variables but the states are observed after a variable transformation.
Our goal is to predict, with uncertainty, solutions of \cref{eq:euler-conservative}--\cref{eq:periodicBCs} from $t_0 = 0$ to the final time $t_f = 0.15~\text{s}$.

\Cref{fig:euler-dims} shows the POD singular value decay of the snapshot data, which decreases sharply and levels out after $6$ modes.
Compressed snapshot data for mode $7$ and greater concentrate near zero and do not result in a tight GP fit, hence we set the ROM dimension to $r = 6$.
The GP estimation times $\vb*{t}'$ are defined to be $m' = 400$ uniformly spaced times in the interval $[t_0, t_{m-1}]$.
To further validate the quality of the training data produced by each GP regression, we solve~\cref{eq:euler-all} on a fine time grid and explicitly evaluate the right-hand side of~\cref{eq:euler-conservative}

\newpage

\begin{figure}[!ht]
    \centering
    \includegraphics[width=\textwidth]{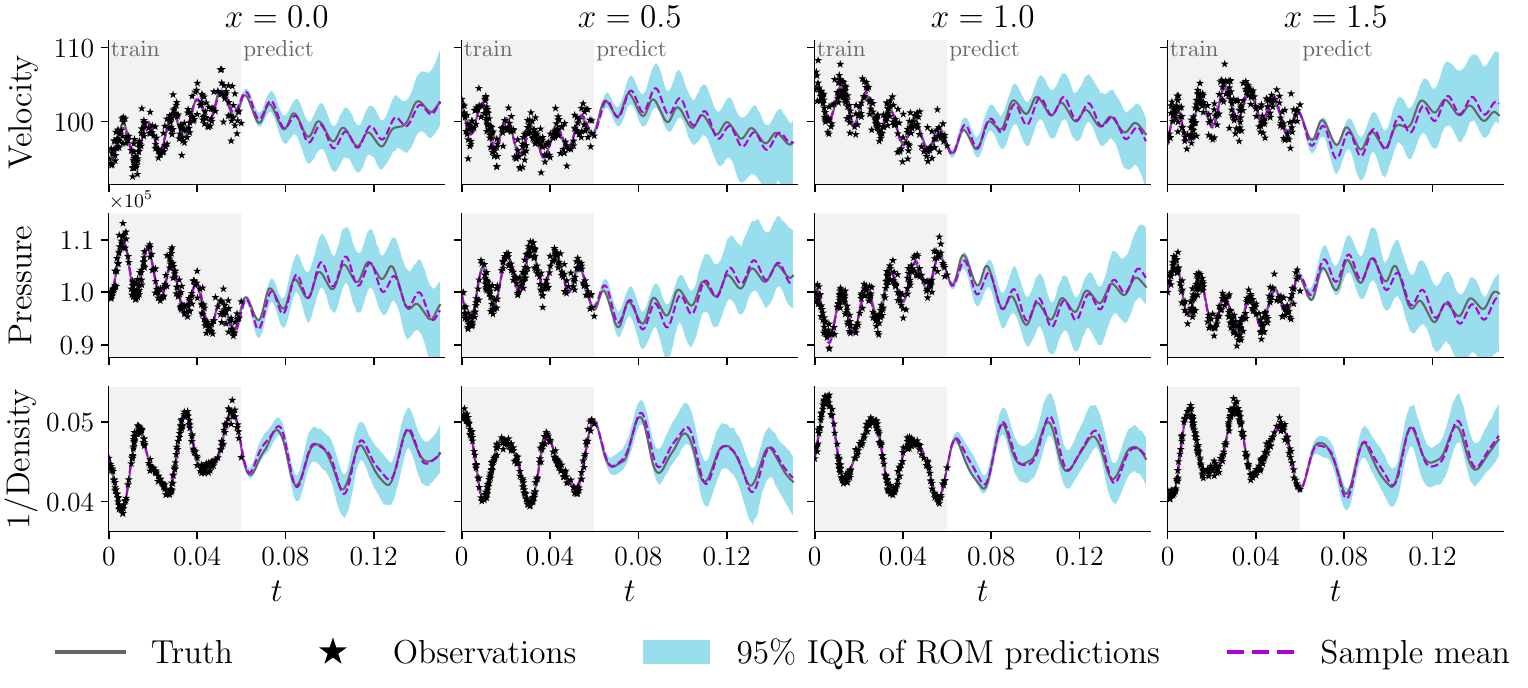}
    \vspace{-0.75cm}
    \caption{Bayesian ROM predictions of \Cref{fig:euler-noisy-reduced} expressed in the specific volume variables $(v,p,1/\rho)$. The solid line shows the true numerical solution to the system without noise. Results are shown at four spatial locations.\vspace{-0.5cm}}
    \label{fig:euler-noisy-full}
\end{figure}

\begin{figure}[!ht]
    \centering
    \includegraphics[width=\textwidth]{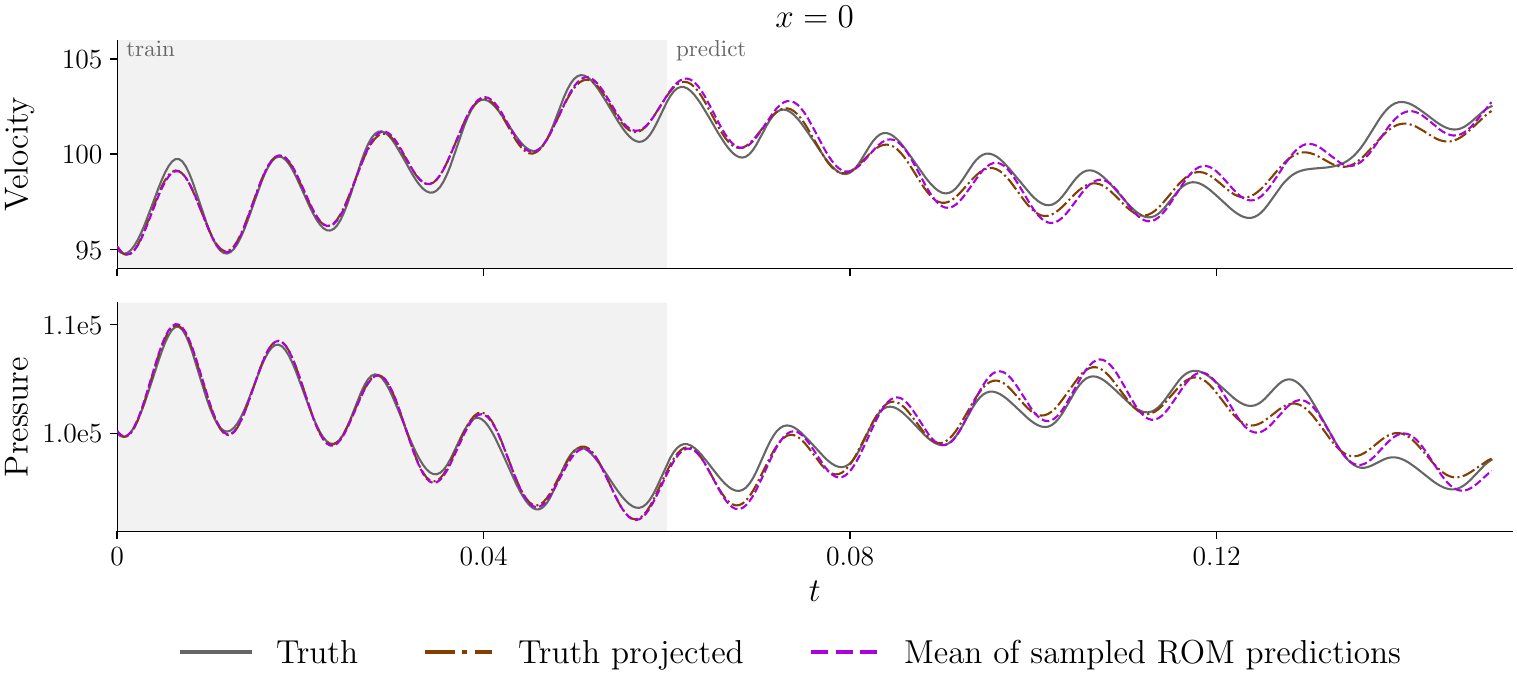}
    \vspace{-0.75cm}
    \caption{Closeup of the results in \Cref{fig:euler-noisy-full} at $x = 0$, as well as the projection of the true numerical solution. In the training regime, the mean of the sampled ROM predictions tightly match the projection of the truth, a consequence of \Cref{alg:OpInfError}. In the prediction regime, when the sample mean deviates from the truth it is generally close to the projection of the truth.}
    \label{fig:euler-noisy-closeup}
\end{figure}

\noindent
(without adding noise) to obtain time derivatives $(\partial_t\rho,\partial_t[\rho v],\partial_t[\rho e])$, which are then transformed to the specific-volume formulation $(\partial_t v, \partial_t p, \partial_t \xi)$ and compressed using the POD basis.
\Cref{fig:euler-ddts} plots these true time derivatives together with the GP estimates for the time derivatives which, over the time grid $\vb*{t}'$ and for the $i$-th reduced state, are given by the mean vector $\tilde{\vb{z}}_i = \vb{K}_{i}^{zy}(\vb{K}_{i}^{yy})^{-1}\vb{y}_i$ with covariance $(\vb{W}_{i}^{zz})^{-1} = \vb{K}_{i}^{zz} - \vb{K}_{i}^{zy}(\vb{K}_{i}^{yy})^{-1}\vb{K}_{i}^{yz}$. The GP estimates are highly faithful to the true time derivatives, and the associated variances so small that the region within $10$ standard deviations of the mean is only slightly visible. As a comparison, we also apply second-order finite differences to the noisy compressed snapshots to produce time derivative estimates, the standard practice for OpInf problems. Although the POD basis filters much of the noise from the high-dimensional observations (see \Cref{fig:euler-dims}), the remaining noise results in significant error in the finite difference estimates, rendering them unsuitable for populating the OpInf regression~\cref{eq:opinf-all}.

After applying \Cref{alg:GP-BayesOpInf} to obtain a posterior distribution for the reduced operators $\hat{\vb{O}}$, we sample from the posterior $500$ times and solve the corresponding ROMs \cref{eq:euler-rom} using the same Runge--Kutta scheme as before.
The uncertainty stemming from the data noise propagates to the ROM solutions and forms a prediction distribution.
To characterize the distribution, we compute the mean and the $2.5\%$ and $97.5\%$ quantiles of the samples at each time instance, see \Cref{fig:euler-draws}.
The individual draws are nearly identical in the training regime, but at the start of the prediction regime the width of the distribution increases noticeably and thereafter generally continues to increase in time at a relatively slow rate.

The results are compared to the true solution in \Cref{fig:euler-noisy-reduced} for the reduced modes and in \Cref{fig:euler-noisy-full} for the original state space in the (physically relevant) specific volume variables.
In all variables and throughout the spatial domain, each ROM solution closely matches the true solution during the training regime, which is expected since \Cref{alg:GP-BayesOpInf} minimizes the error of the sample mean and the GP estimates of the reduced training data.
Although the ROM solutions are not accurate individually in the prediction regime, in aggregate they form an accurate estimate of the projection of the true solution.
\Cref{fig:euler-noisy-closeup} shows a closeup of results from \Cref{fig:euler-noisy-full} at a single spatial location and includes the true solution $\vb{q}(t)$, the projection $\tilde{\vb{q}}(t) = \vb{V}\vb{V}\trp(\vb{q}(t) - \bar{\vb{q}}) + \bar{\vb{q}}$, and the pointwise mean of the sampled ROM predictions.
The error between the truth $\vb{q}(t)$ and the projection $\tilde{\vb{q}}(t)$ is due entirely to the low-dimensional state approximation \cref{eq:PODcentered}, not the GP modeling or the OpInf regression.
The sample mean is about as accurate as the projection of the truth, hence, in this experiment, the error in the sample mean is dominated by the fundamental limitations incurred by dimensionality reduction.

\begin{figure}[t]
    \centering
    \includegraphics[width=\textwidth]{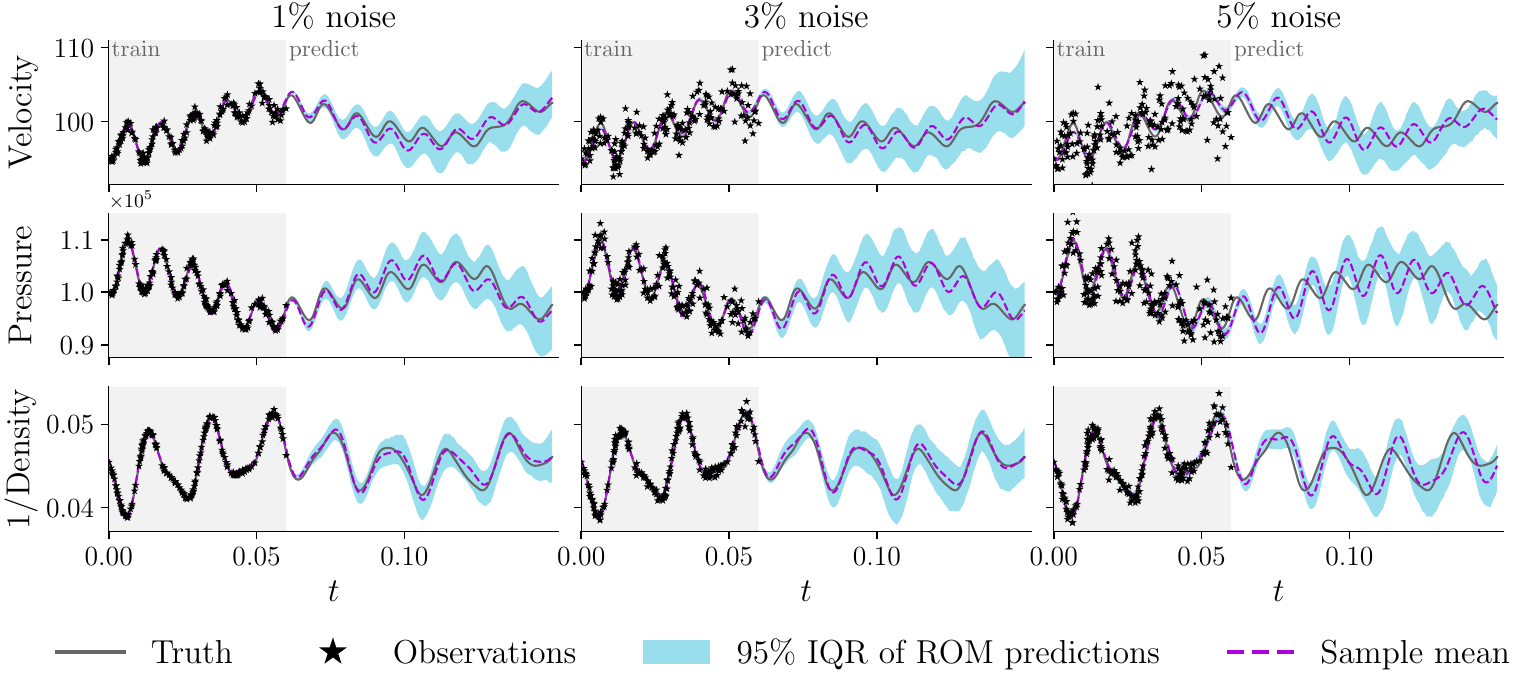}
    \vspace{-0.75cm}
    \caption{Pointwise Bayesian ROM predictions of the specific volume variables at $x = 0$ for the Euler system \cref{eq:euler-conservative}--\cref{eq:periodicBCs} using $m=200$ observations over the time domain $[0, 0.06]$ with varying noise levels.}
    \label{fig:euler-noisy-comparison}
\end{figure}

\Cref{fig:euler-noisy-comparison} repeats the numerical experiment with $\xi = 1\%$ and $\xi = 5\%$ relative noise and reports the results at a single point in space.
As the noise level $\xi$ increases, so does the width of the distribution of predicted ROM solutions.
The sample mean of the ROM predictions estimates the truth well for $\xi = 1\%$ and $\xi = 3\%$, but with $\xi = 5\%$ the ROM solutions are off phase from the truth in the prediction regime.
In the latter case, the large amount of noise causes the ROM solutions in the training regime (and, though not shown, the projection of the truth) to be less accurate, and at times the true solution is not well represented by the prediction distribution.
However, the solutions still roughly capture the trends of the system behavior, especially in the specific volume.

\begin{figure}[t]
    \centering
    \includegraphics[width=\textwidth]{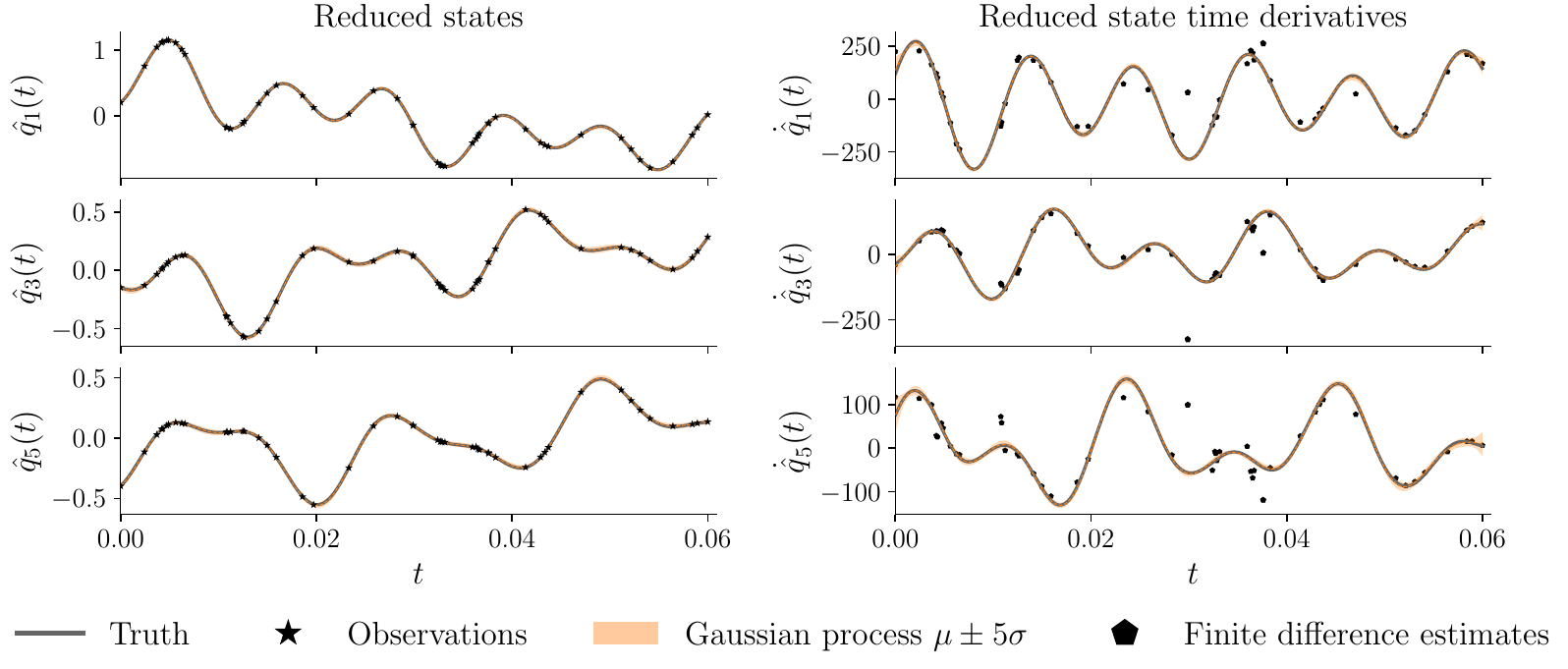}
    \vspace{-0.75cm}
    \caption{True reduced states, observations, and trained GPs for reduced modes 1, 3, and 5~(left), and corresponding reduced state time derivatives, GP estimates, and finite difference estimates~(right).}
    \label{fig:euler-sparse-fitandddts}
\end{figure}

\begin{figure}[!ht]
    \centering
    \includegraphics[width=\textwidth]{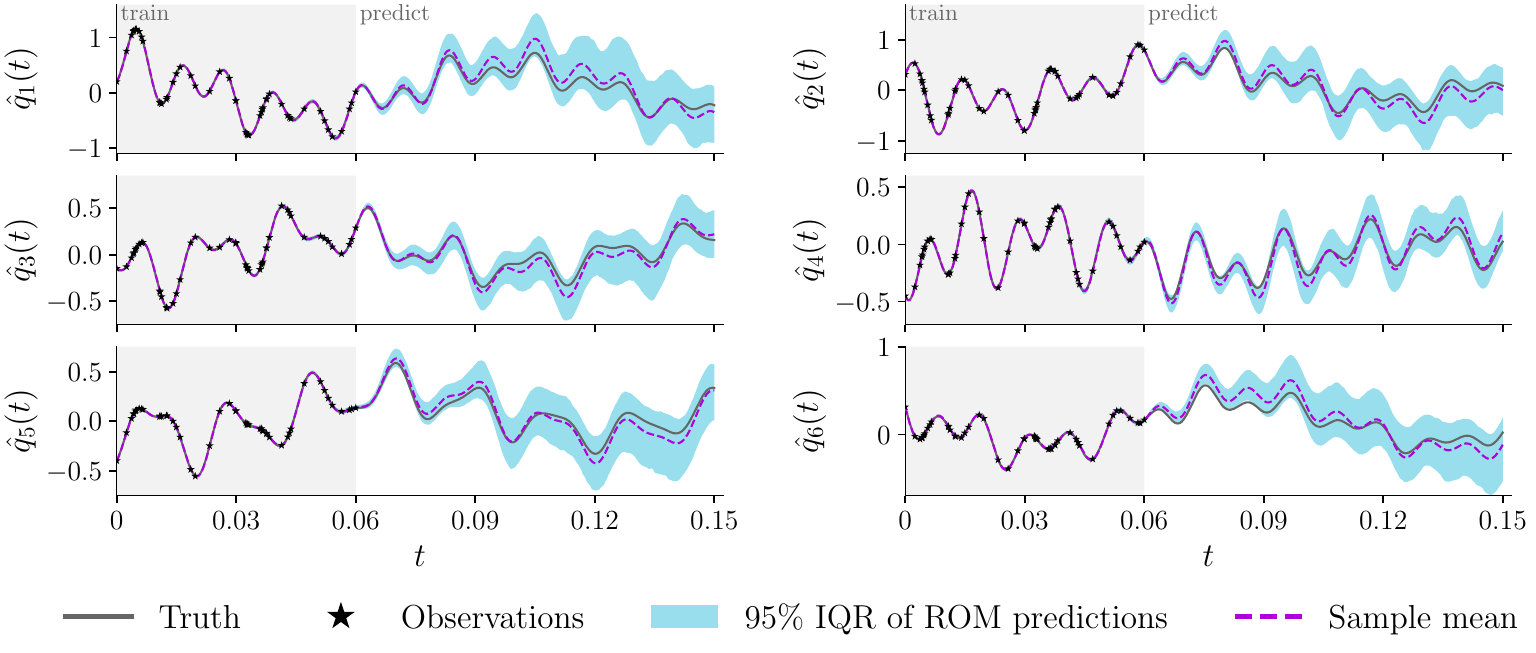}
    \vspace{-0.75cm}
    \caption{Bayesian ROM predictions for the Euler system \cref{eq:euler-conservative}--\cref{eq:periodicBCs} using $m=50$ observations over the time domain $[0, 0.06]$ with $\xi=1\%$ relative noise. The solid line shows the true numerical solution to the system without noise, expressed in the POD coordinates. The posterior is sampled $500$ times. The pointwise mean of the sample predictions gives an accurate estimate of the true solution, with uncertainty described by the $95\%$ interquantile range.}
    \label{fig:euler-sparse-reduced}
\end{figure}

\subsubsection{Sparse training data}

Next, consider a scenario in which the observed data are sparse in time.
Let $m = 50$, again with final observation time $t_{m-1} = 0.06~\text{s}$ but now with only $\xi = 1\%$ relative noise. This is the scenario depicted in \Cref{fig:euler-fomdata}.
\Cref{fig:euler-sparse-fitandddts} shows GP regressions of the compressed state data and predictions for the associated time derivatives. In this case, finite differences of the noisy snapshots result in estimates for the time derivatives that are much less unwieldy than in the previous experiment (see \Cref{fig:euler-ddts}), but they are still significantly less accurate than the GP time derivative estimates. Additionally, the GPs yield accurate training data for both the reduced state and the time derivative over the entire training time regime, not just at the sparse observation times.

\begin{figure}[t]
    \centering
    \includegraphics[width=\textwidth]{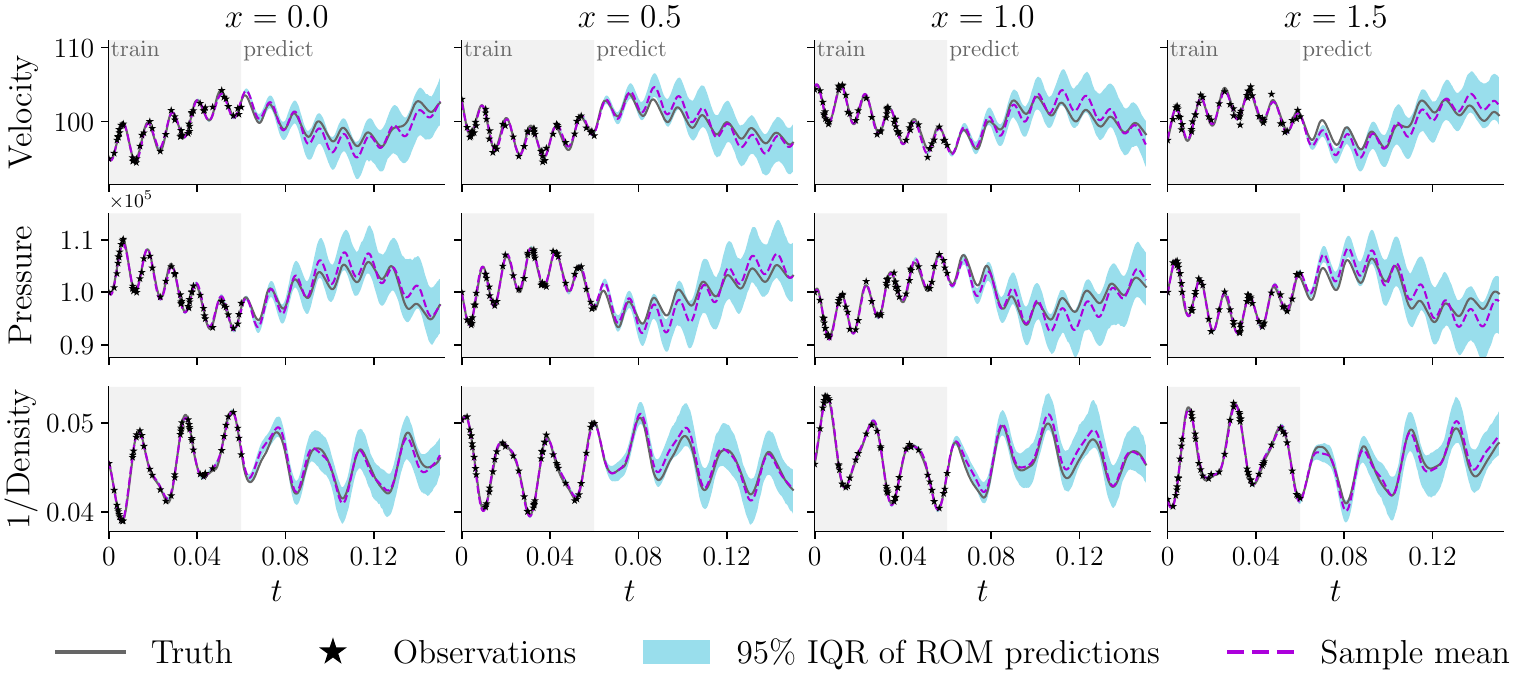}
    \vspace{-0.75cm}
    \caption{Bayesian ROM predictions of \Cref{fig:euler-sparse-reduced} expressed in the specific volume variables $(v,p,1/\rho)$ at four spatial locations.}
    \label{fig:euler-sparse-full}
\end{figure}

\begin{figure}[!ht]
    \centering
    \includegraphics[width=\textwidth]{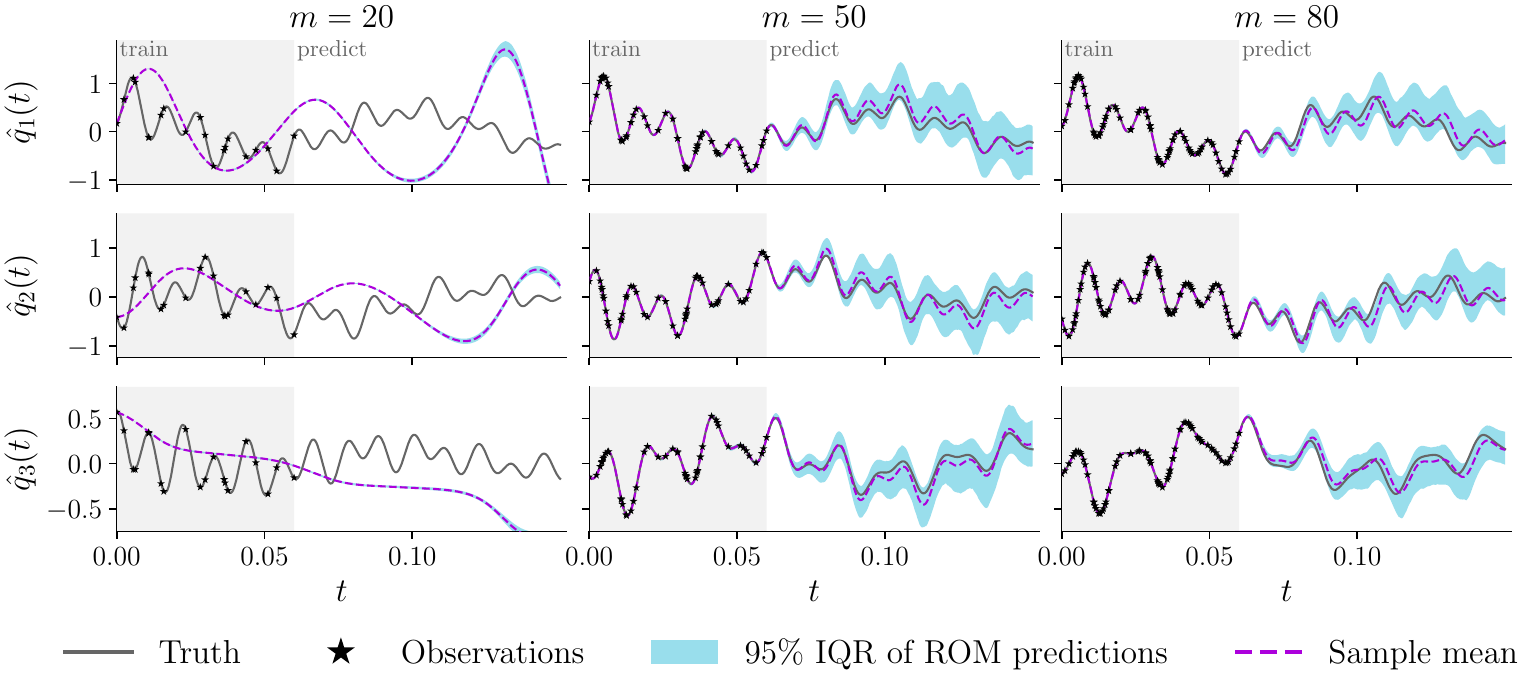}
    \vspace{-0.75cm}
    \caption{Bayesian ROM predictions for the Euler system \cref{eq:euler-conservative}--\cref{eq:periodicBCs} using a varied number of observations over the time domain $[0, 0.06]$ with $\xi=1\%$ relative noise. With only $m = 20$ observations \Cref{alg:GP-BayesOpInf} fails to produce an accurate model.}
    \label{fig:euler-sparse-comparison}
\end{figure}

Results with $500$ posterior samples are shown in \Cref{fig:euler-sparse-reduced} in the reduced space and \Cref{fig:euler-sparse-full} in the physical state space. The results are similar in quality to the previous experiments: the truth is reconstructed accurately in the training regime and the sample mean of the ROM predictions is a reasonable estimate of the truth in the prediction regime. To further probe the performance of our method on sparse data, we repeat the experiment with only $m = 20$ observations and again with $m = 80$ observations, all drawn uniformly over the same training time domain.
\Cref{fig:euler-sparse-comparison} displays the results at a single spatial location.
With only $m = 20$ data points, GP-BayesOpInf fails to produce accurate solutions and also does not produce a wide prediction distribution.
This is partially due to not having sufficient data to perform a good GP fit, but also because \Cref{alg:OpInfError} prioritizes stable but inaccurate ROM solutions over ROM solutions that are accurate for a limited time but which become unstable before the final simulation time.
On the other hand, using $m = 80$ training points produces a slightly more accurate sample mean and tighter prediction distribution than with $m=50$, i.e., more training data leads to improved accuracy.

\subsection{Nonlinear diffusion-reaction equation}

In this section, we use GP-BayesOpInf to construct probabilistic surrogates for a nonlinear diffusion-reaction equation using data from multiple trajectories.
Each trajectory uses the same initial conditions but a different instance of a parameterized nonlinear source term. The efficacy of the probabilistic surrogate is shown for prediction in both time and in parameters defining the source term with noisy, sparse training data.

\subsubsection{Governing equations}

Consider the following forced heat diffusion model with a cubic reaction term,
\begin{subequations}
\begin{align}
    \label{eq:heat-cubic}
    \begin{aligned}
    \frac{\partial q}{\partial t}
    &= \mu \frac{\partial^{2} q}{\partial x^{2}}
        - q^3
        + u,
    && x\in(0, 1), t\in(0, T],
    \\
    q(0, t) &= 0,
    \quad
    q(1, t) = 1,
    && t\in(0, T],
    \\
    q(x, 0) &= q_{0}(x),
    && x\in [0, 1],
    \end{aligned}
\end{align}
where $\mu = 0.005$ is the diffusion constant and where the source term $u$, parameterized by constants $a,b\in\mathbb{R}$, is given by
\begin{align}
    \label{eq:heat-forcing}
    u(x, t; a, b)
    = \frac{a\sin(2\pi t)}{1 + 100\left(x - \frac{1}{4}\right)^2}
    + \frac{b\sin(4\pi t)}{1 + 100\left(x - \frac{3}{4}\right)^2}\,.
\end{align}
The initial conditions are defined to be
\begin{align}
    \label{eq:heat-initial}
    q_0(x)
    &= x(1 - x)\left(
        6 (1-x)^2 e^{-x}
        - 10 e^{x} \sin(x/6)
    \right) + x,
\end{align}
\end{subequations}
a nonpolynomial function satisfying the boundary conditions.
In the following experiments, we use training data from $\ell = 5$ trajectories obtained by solving \cref{eq:heat-cubic}--\cref{eq:heat-initial} for the $(a,b)$ pairs $(-2,0)$, $(-1,-2)$, $(0,1)$, $(1,-1)$, and $(2,2)$, which form a coarse Latin square sampling of $[-2,2]\times[-2,2]$.

\subsubsection{ROM structure}
To choose a ROM structure for \cref{eq:heat-cubic}, we introduce the auxiliary variable $w = q^2$ to expose quadratic dynamics~\cite{qian2021thesis}.
As $\partial_t w = 2qq_t$, \cref{eq:heat-cubic} can be written equivalently as the PDE system
\begin{align}
    \label{eq:heat-lifted}
    \begin{aligned}
    \frac{\partial q}{\partial t}
    &= \mu \frac{\partial^{2} q}{\partial x^{2}}
        - qw
        + u,
&\qquad
    \frac{\partial w}{\partial t}
    &= 2\mu q\frac{\partial^{2} q}{\partial x^{2}}
        - w^2
        + qu.
\end{aligned}
\end{align}
Let $\vb{q}(t) = [~\vb{s}(t)\trp~~(\vb{s}(t)\ast\vb{s}(t))\trp~]\trp\in\mathbb{R}^{N}$ be the spatial discretization of $(q, w)$, where $N = 2n_x = 1000$ and $\ast$ is the Hadamard (element-wise) product.
The discretization of \cref{eq:heat-lifted} is a system of ODEs that can be written as
\begin{align}
    \frac{\text{d}}{\text{d}t}\vb{q}(t)
    = \vb{Aq}(t) + \vb{H}[\vb{q}(t)\otimes\vb{q}(t)] + \vb{Bu}(t) + \vb{N}[\vb{u}(t)\otimes\vb{q}(t)],
\end{align}
where $\vb{A}\in\mathbb{R}^{N\times N}$, $\vb{H}\in\mathbb{R}^{N\times N^2}$, $\vb{B}\in\mathbb{R}^{N\times 2}$, $\vb{N}\in\mathbb{R}^{N\times 2N}$, and $\vb{u}(t;a,b) = (a\sin(2\pi t), b\sin(4\pi t))\trp\in\mathbb{R}^{2}$ gathers the time-dependent components of the source term $u$.
Using a centered-POD approximation $\vb{q}(t) \approx \vb{V}\hat{\vb{q}}(t) + \bar{\vb{q}}$ as in \cref{eq:PODcentered}, we obtain the ROM structure
\begin{align}
    \label{eq:heat-rom}
    \frac{\textrm{d}}{\textrm{d}t}\hat{\vb{q}}(t)
    &= \hat{\vb{c}} + \hat{\vb{A}}\hat{\vb{q}}(t) + \hat{\vb{H}}[\hat{\vb{q}}(t)\otimes\hat{\vb{q}}(t)] + \hat{\vb{B}}\vb{u}(t) + \hat{\vb{N}}[\vb{u}(t)\otimes\hat{\vb{q}}(t)],
\end{align}
where $\hat{\vb{q}}(t)\in\mathbb{R}$, $\hat{\vb{c}}\in\mathbb{R}^{r}$, $\hat{\vb{A}}\in\mathbb{R}^{r\times r}$, $\hat{\vb{H}}\in\mathbb{R}^{r\times r^2}$, $\hat{\vb{B}}\in\mathbb{R}^{r\times 2}$, and $\hat{\vb{N}}\in\mathbb{R}^{r\times 2r}$.
This system is \cref{eq:reducedorder} with operator matrix $\hat{\vb{O}} = [~\hat{\vb{c}}~~\hat{\vb{A}}~~\hat{\vb{H}}~~\hat{\vb{B}}~~\hat{\vb{N}}~]\in\mathbb{R}^{r\times d(r,p)}$, data vector $\vb{d}(\hat{\vb{q}},\vb{u}) = [~1~~\hat{\vb{q}}\trp~~(\hat{\vb{q}}\otimes\hat{\vb{q}})\trp~~\vb{u}\trp~~(\vb{u}\otimes\hat{\vb{q}})\trp~]\trp\in\mathbb{R}^{d(r,p)}$, and dimension $d(r,p) = 1 + r + r^2 + p + rp$.
The constant term $\hat{\vb{c}}$ arises due to the centering with $\bar{\vb{q}}$.
Having selected this structure, the snapshots to be used in the procedure are $\vb{q}_{j} = [~\vb{s}_{j}\trp~~(\vb{s}_j\ast\vb{s}_j)\trp~]\trp\in\mathbb{R}^{N}$.
Because noise is not applied at the spatial boundary, the centered snapshots $\vb{q}_{j} - \bar{\vb{q}}$ obey homogeneous Dirichlet boundary conditions, hence all ROM predictions $\vb{V}\hat{\vb{q}}(t) + \bar{\vb{q}}$ automatically satisfy the boundary conditions of~\cref{eq:heat-cubic}.

\subsubsection{Data generation}

\begin{figure}[t]
    \centering
    \includegraphics[width=\textwidth]{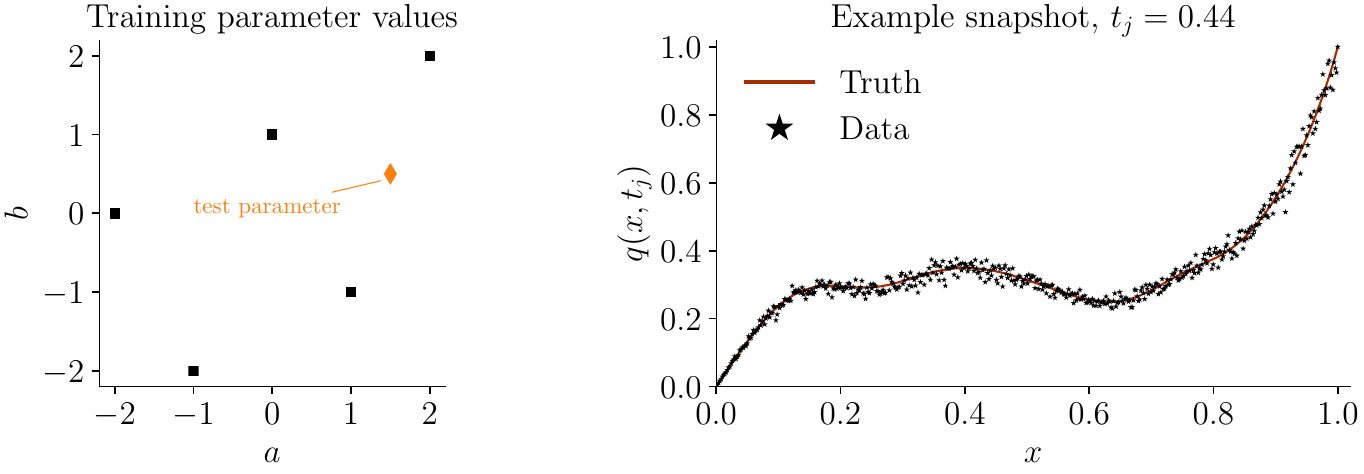}
    \vspace{-0.75cm}
    \caption{The five sets of input parameter values $(a,b)$ parameterizing the input function $u(x,t;a,b)$ of \cref{eq:heat-forcing} used in generating training data (left) and an example snapshot with $\xi = 5\%$ noise for $a-1$, $b=-2$ (right).}
    \label{fig:heat3-samples}
\end{figure}

Data for this experiment are generated by using second-order central finite differences to discretize the spatial derivative in~\cref{eq:heat-cubic} with $n_x = 500$ degrees of freedom, then time stepping the system with an implicit multi-step backward difference method of variable order (1 to 5) with quasi-constant step size~\cite{byrne1975odes,shampine1997ode}.
In this experiment, we observe only $m = 20$ snapshots for each choice of $(a, b)$, totaling $100$ snapshots, with different sampling times~$\vb*{t}$ for each trajectory.
In each case, we set $t_0 = 0$, $t_{m-1} = 1$, and sample the remaining points $(t_1,\ldots,t_{m-2})$ from the uniform distribution over $(0, t_{m-1})$. The snapshots are corrupted with additive noise whose variance is determined by the magnitude of the state at each point in space, i.e.,
\begin{align}
    \vb{s}_j
    = \vb{s}(t_j) + \vb*{\varepsilon}_j
    \in\mathbb{R}^{n_x}\,,
    \qquad
    \vb*{\varepsilon}_j \sim \mathcal{N}(\vb{0},\operatorname{diag}(\xi\vb{s}(t_j))^2),
\end{align}
where $\vb{s}(t_j)$ is the numerical solution for $q$ at time $t_j\in\vb*{t}$ and $\xi = 5\%$ is the noise level.
Noise is not added at $t = t_0$ or to the spatial boundaries since the initial conditions and Dirichlet boundary conditions are given.
\Cref{fig:heat3-samples} visualizes the $(a,b)$ parameter space and plots an example noised snapshot in space.

\subsubsection{Prediction at new input parameters}

\begin{figure}[!ht]
    \centering
    \vspace{.25in}
    \includegraphics[width=\textwidth]{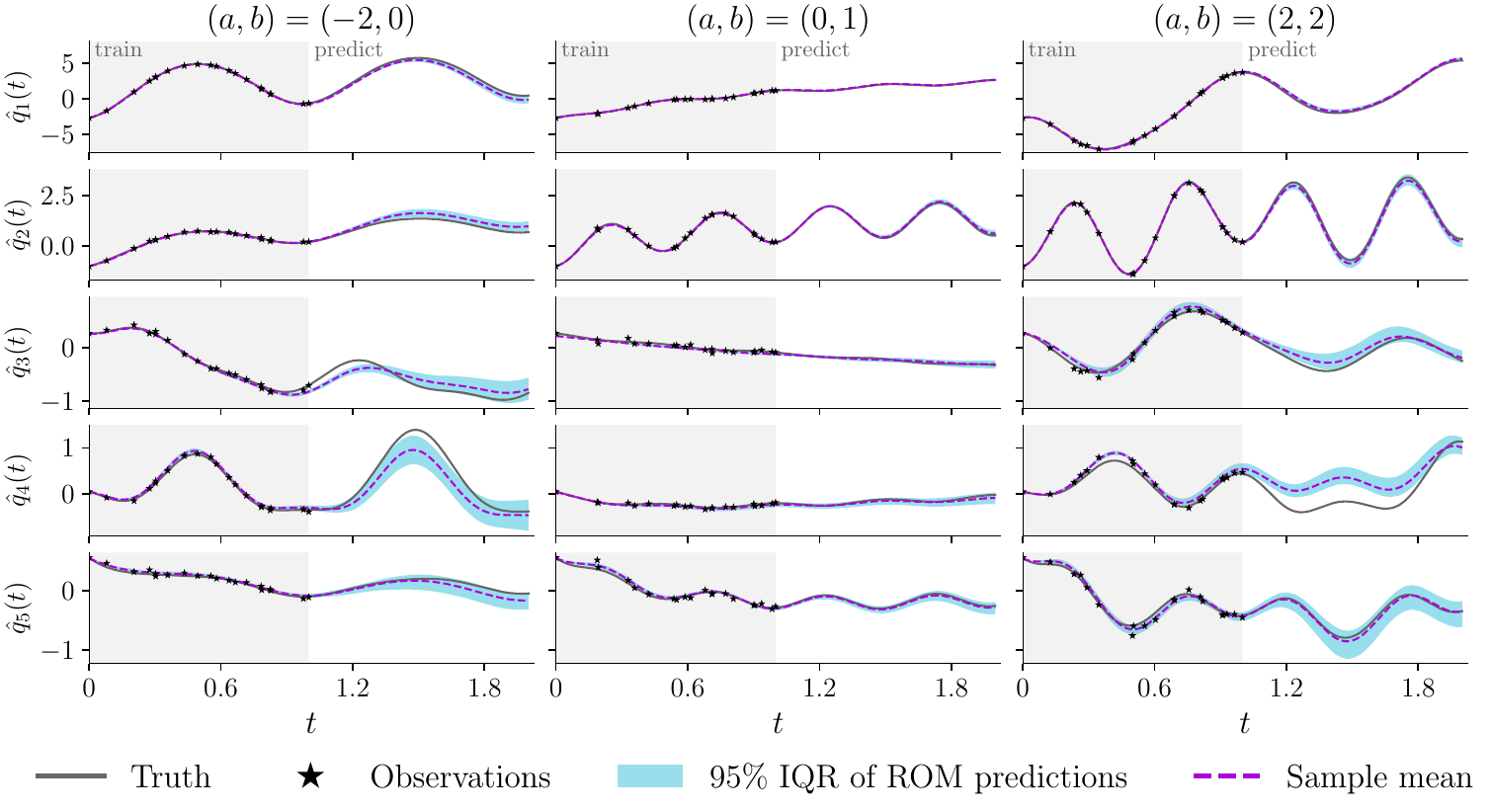}
    \vspace{-0.75cm}
    \caption{Bayesian ROM predictions for the nonlinear diffusion-reaction problem \cref{eq:heat-cubic}--\cref{eq:heat-initial} for three values of the input parameters $(a,b)$. For each set of input parameters, data are observed at $m=20$ points over the time domain $[0,1]$ with $\xi = 5\%$ relative noise. Five hundred posterior samples are used to generate ROM solutions, and the resulting distribution is characterized by the pointwise sample mean and $95\%$ interquantile range.}
    \label{fig:heat3-reduced}
\end{figure}

\begin{figure}[!ht]
    \centering
    \includegraphics[width=\textwidth]{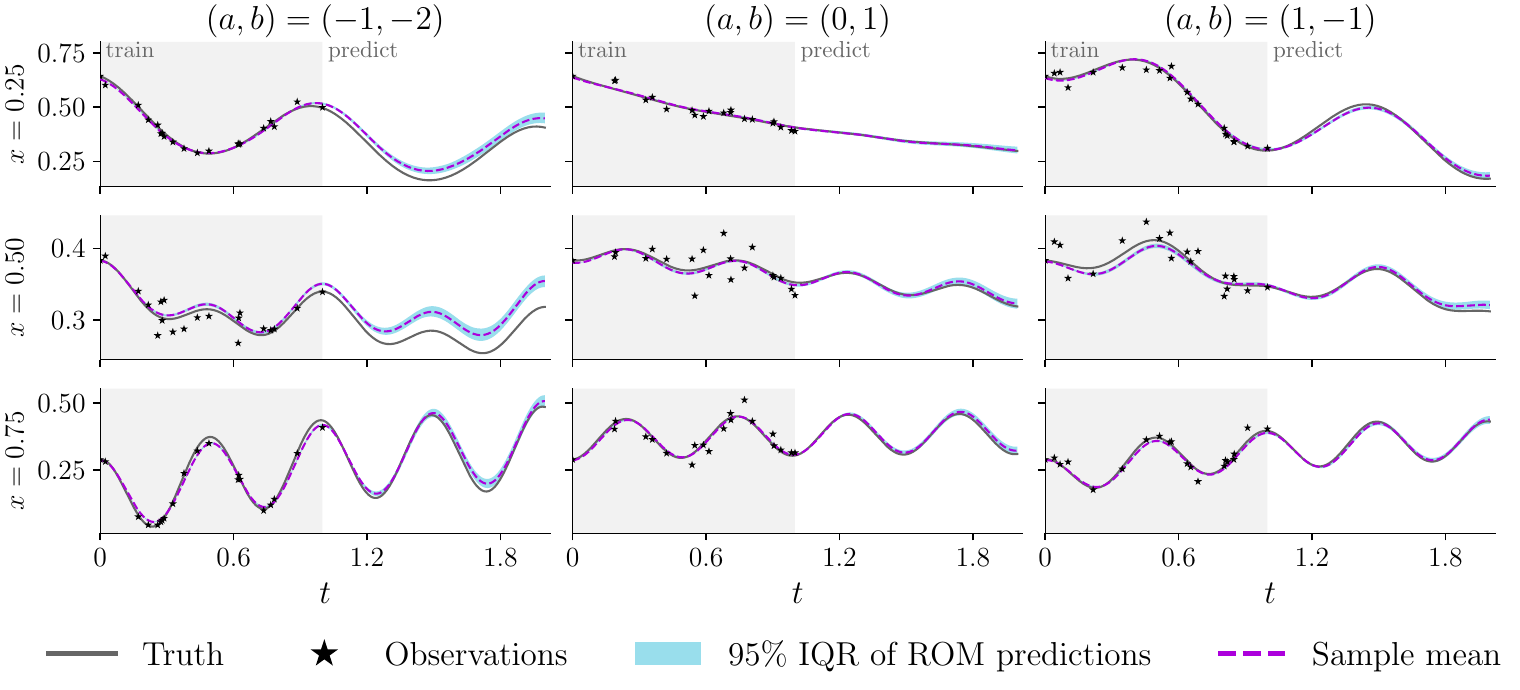}
    \vspace{-0.75cm}
    \caption{Bayesian ROM predictions of \Cref{fig:heat3-reduced} expressed in the original state variable $q$ for three values of the input parameters $(a,b)$. Results are shown at three spatial locations.}
    \label{fig:heat3-full}
\end{figure}

We apply \Cref{alg:GP-BayesOpInfMulti} (the analog of \Cref{alg:GP-BayesOpInf} for multiple trajectories) with $m' = 4m = 80$ GP~estimates to obtain a posterior distribution for $\hat{\vb{O}}$ with $r = 5$ POD modes.
We then sample the posterior $500$ times and solve the corresponding ROM~\cref{eq:heat-rom} with the same backwards-differences time stepper as the full-order model.
For three choices of the input parameters $(a,b)$, \Cref{fig:heat3-reduced} shows results in the reduced space, while \Cref{fig:heat3-full} shows results in the original state space at select spatial locations.
The accuracy of ROM solutions varies across the inputs because some trajectories are more easily captured by the low-dimensional state approximation than others.
Note that the width of the prediction distribution increases whenever the sample mean is less accurate.

\newpage

\newpage

\begin{figure}[t]
    \centering
    \includegraphics[width=\textwidth]{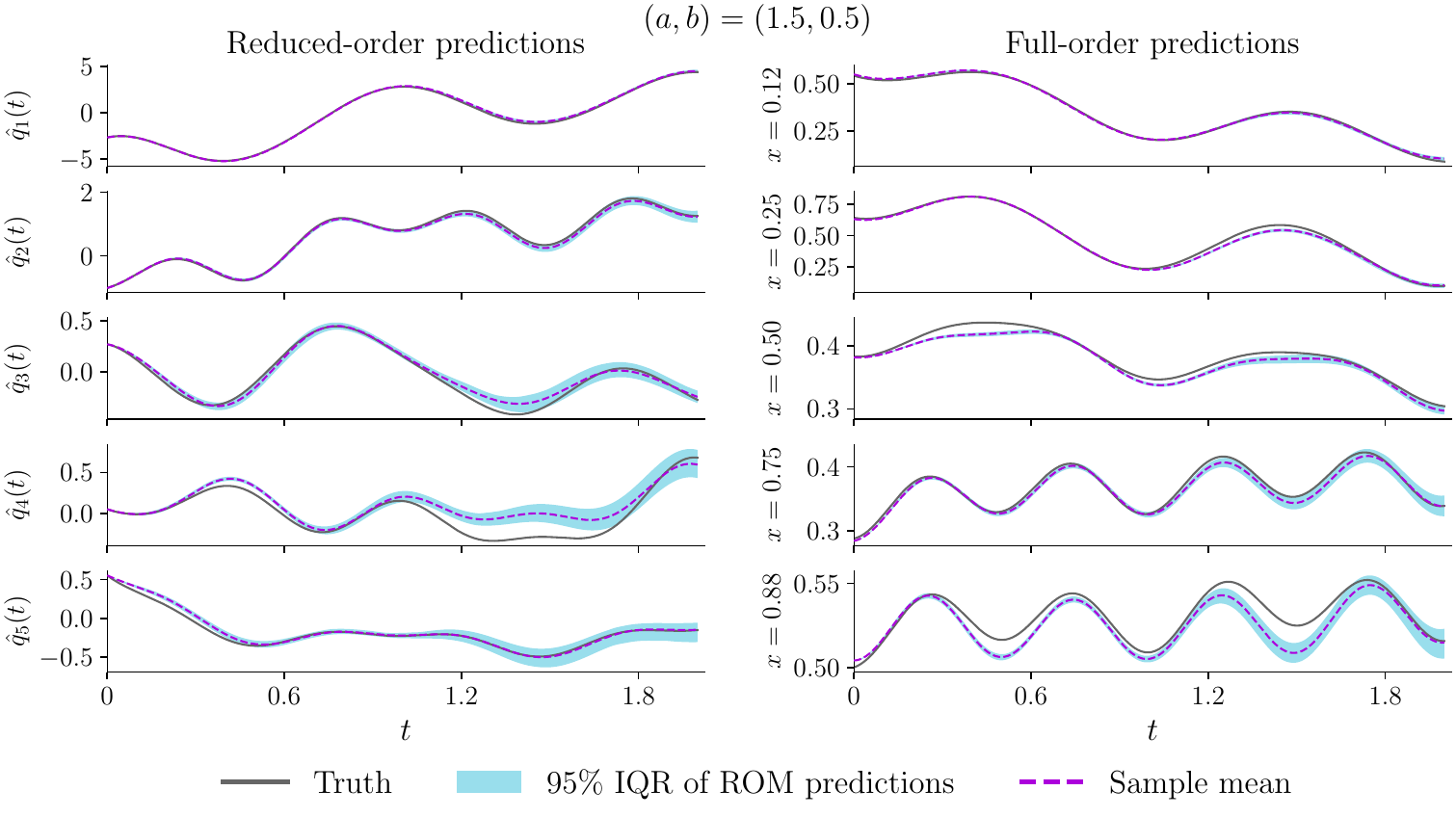}
    \vspace{-0.75cm}
    \caption{Bayesian ROM predictions for the nonlinear diffusion-reaction problem \cref{eq:heat-cubic}--\cref{eq:heat-initial} for new input parameter values $(a,b) = (1.5,0.5)$ not used in the training set, both in the reduced space (left) and pointwise in the original state space (right).}
    \label{fig:heat3-newparams}
\end{figure}

It is worth noting that, in this example, aggregating a few samples each from multiple trajectories tends to result in a more diverse and informative POD basis than when using many samples from a single trajectory, which is why the results are high quality even with the relatively small $m = 20$. This diversity is also helpful for generalization: \Cref{fig:heat3-newparams} shows results for posterior ROM samples using the new input parameter values $(a,b) = (1.5, 0.5)$, i.e., using an input function $\vb{u}(t)$ that was not part of the training set. This represents a prediction in the $(a,b)$ input parameter space as well as extrapolation in time. As in previous experiments, the width of the prediction distribution increases in time, but here the transition from training regime to prediction regime is less pronounced.

\section{GP-based parameter estimation for ODEs}
\label{sec:odeparamestimate}

The GP-BayesOpInf framework can be applied, with a minor adjustment, to a wide class of systems of ODEs that are not the result of spatially discretizing a PDE. Our approach allows for parameters to appear in multiple equations, treats all equations jointly, and leverages the prior variance selection framework developed in \Cref{sec:procedure}.

\subsection{Derivation} \label{sec:ode-methodology}
We consider dynamical systems modeled by ODEs that may be nonlinear in state but which are linear in a set of scalar parameters, i.e.,
\begin{align}
    \label{eq:ode-setup}
    \frac{\text{d}}{\text{d}t}\hat{\vb{q}}(t)
    &=
\vb{S}
    (\hat{\vb{q}}(t),\vb{u}(t))\hat{\vb{o}},
\end{align}
where $\hat{\vb{q}}(t)\in\mathbb{R}^{r}$ and $\vb{u}(t)\in\mathbb{R}^{p}$ are the state and the input (as in the PDE case),
$\vb{S}:\mathbb{R}^{r}\times\mathbb{R}^{p}\to\mathbb{R}^{r\times d(r,p)}$ describes the (potentially nonlinear) system dynamics,
and $\hat{\vb{o}}\in\mathbb{R}^{d(r,p)}$ are the system parameters.
We use $\hat{\vb{q}}(t)$ to be consistent with previous notation, but in this setting the state $\hat{\vb{q}}(t)$ is \emph{not} the result of compressing a high-dimensional state variable.
Our objective is to construct a probability distribution for the parameter vector $\hat{\vb{o}}$ given (noisy) observations of the state and the corresponding inputs.

If data are available for $\hat{\vb{q}}(t)$, a Bayesian inference for $\hat{\vb{o}}$ can be formulated similar to~\cref{eq:posterior-inference}.
Let $\vb*{t}_{i}\in\mathbb{R}^{m_i}$ denote the $m_i\in\mathbb{N}$ times at which the $i$-th component of $\hat{\vb{q}}(t)$ is observed as the (possibly noisy) vector $\vb{y}_{i}\in\mathbb{R}^{m_i}$, $i = 1,\ldots,r$, and define the concatenated state data vector
\begin{align}
    \vb{y}
    = \left[\begin{array}{c}
        \vb{y}_1 \\ \vdots \\ \vb{y}_{r}
    \end{array}\right]
    \in \mathbb{R}^{m},
    \qquad
    m = \sum_{i=1}^{r}m_i.
\end{align}
As in \cref{eq:gp-single}, let $\hat{q}_{i}(t)$ denote a GP for the $i$-th component of $\hat{\vb{q}}(t)$.
Define $\vb*{Y}$ to be the random $m$-vector whose entries are Gaussian random variables produced by evaluating $\hat{q}_{i}$ over $\vb*{t}_i$ for each $i = 1,\ldots, r$ and concatenating the results.
Hence, $\vb{y}$ is an observation of $\vb*{Y}$.
Similarly, for estimation times $\vb*{t}_i'\in\mathbb{R}^{m'}$, let $\vb*{Z}$ be the random vector with $rm'$ entries given by concatenating the time derivative GPs $\dot{\hat{q}}_{i}$ evaluated at $\vb*{t}'$, $i=1,\ldots, r$.
Finally, let $\tilde{\vb{D}}\in\mathbb{R}^{rm'\times d(r,p)}$ contain state and input data at the estimation times~$\vb*{t}'$
following the structure of $\vb{S}$
and define $\vb*{R} = \vb*{Z} - \tilde{\vb{D}}\hat{\vb{o}}$.
We then have the parameter vector inference
\begin{align}
    \begin{aligned}
    \pi_\texttt{post}(\hat{\vb{o}})
    = p(\hat{\vb{o}}\mid \vb*{Y}=\vb{y},\vb*{R}=\vb{0},\tilde{\vb{D}},\vb*{\theta}_{1},\ldots,\vb*{\theta}_{r},\vb*{\gamma})
    \propto \pi_\texttt{prior}(\hat{\vb{o}})\pi_\texttt{like}(\hat{\vb{o}}),
    \quad\text{with}
    \\
    \pi_\texttt{prior}(\hat{\vb{o}})
    = p(\hat{\vb{o}}\mid\vb*{\gamma})
    \propto \exp\left(-\frac{1}{2}\|\operatorname{diag}(\vb*{\gamma})\hat{\vb{o}}\|_{2}^{2}\right)\,,
    \\
    \pi_\texttt{like}(\hat{\vb{o}})
    = p(\vb*{Y}=\vb{y},\vb*{R}=\vb{0}\mid\hat{\vb{o}},\tilde{\vb{D}},\vb*{\theta}_{1},\ldots,\vb*{\theta}_{r})\,,
    \end{aligned}
\end{align}
where, as before, $\vb*{\theta}_{i}$ denotes the hyperparameters of the GP $\hat{q}_{i}$.
This likelihood aligns with the gradient matching likelihood used in computational statistics \cite{calderhead2008odegps,tronarp2019probabilisticodes,ye2024gaussian}.

The mean and covariance of the posterior $\hat{\vb{o}}\sim\mathcal{N}(\vb*{\mu},\vb*{\Sigma})$ are given by
\begin{subequations}
\begin{align}
    \vb*{\mu}
    &= \vb*{\Sigma}\tilde{\vb{D}}\trp\vb{W}\tilde{\vb{z}}\,,
    &
    \vb*{\Sigma}
    &= \left(\tilde{\vb{D}}\trp\vb{W}\tilde{\vb{D}} + \vb*{\Gamma}\trp\vb*{\Gamma}\right)^{-1}\,,
\end{align}
in which
\begin{align}
    \tilde{\vb{z}}
    = \left[\begin{array}{c}
        \tilde{\vb{z}}_{1} \\ \vdots \\ \tilde{\vb{z}}_{r}
    \end{array}\right]\in\mathbb{R}^{rm'},
    \qquad
    \vb{W}
    = \operatorname{diag}\left(\vb{W}_{1}^{zz},\ldots,\vb{W}_{r}^{zz}\right)
    \in\mathbb{R}^{rm'\times rm'},
    \qquad
    \vb*{\Gamma} = \operatorname{diag}(\vb*{\gamma})\in\mathbb{R}^{d(r,p)\times d(r,p)},
\end{align}
\end{subequations}
where $\vb{W}_{i}^{zz}\in\mathbb{R}^{m'\times m'}$ and $\tilde{\vb{z}}_{i}\in\mathbb{R}^{m'}$ are defined by the GP $\hat{q}_{i}$ as in \cref{eq:KsandWs} and \cref{eq:gpreconstruct}, respectively.
The generalized least-squares problem \cref{eq:generalizedlstsq} for the posterior mean $\vb*{\mu}$ is
\begin{align}
    \min_{\vb*{\eta}}\left\|
        \tilde{\vb{D}}\vb*{\eta}
        - \tilde{\vb{z}}
    \right\|_{\vb{W}}^{2}
    + \|\vb*{\Gamma\xi}\|_{2}^{2}.
\end{align}

\subsection{SEIRD epidemiological model}

In this section, we consider a system of five ODEs with five scalar parameters that defines the SEIRD compartmental model for infectious disease spreading through a population~\cite{korolev2021identification}:
\begin{align}
    \label{eq:seird-equations}
    \begin{aligned}
    \frac{\textup{d}}{\textup{d}t}\hat{q}_S(t)
    &= -\beta \hat{q}_S(t) \hat{q}_I(t),
    &
    \frac{\textup{d}}{\textup{d}t}\hat{q}_R(t)
    &= (1 - \alpha) \gamma \hat{q}_I(t),
    \\
    \frac{\textup{d}}{\textup{d}t}\hat{q}_E(t)
    &= \beta \hat{q}_S(t) \hat{q}_I(t) - \delta \hat{q}_E(t),
    &
    \frac{\textup{d}}{\textup{d}t}\hat{q}_D(t)
    &= \alpha \rho \hat{q}_I(t),
    \\
    \frac{\textup{d}}{\textup{d}t}\hat{q}_I(t)
    &= \delta \hat{q}_E(t) - (1 - \alpha) \gamma \hat{q}_I(t) - \alpha \rho \hat{q}_I(t).
    \end{aligned}
\end{align}
The state variables represent the proportions of the population that are susceptible to the disease ($\hat{q}_{S}(t)$), exposed but not yet infected ($\hat{q}_E(t)$), currently infected ($\hat{q}_I(t)$), recovered (and immune) post-infection ($\hat{q}_{R}(t)$), and deceased post-infection ($\hat{q}_D(t)$).
The system parameters describe the rates at which individuals transition from one state to another:
$\beta > 0$ is the expected number of people an infected person exposes to the disease per day,
$\delta > 0$ is the incubation period for transition from exposed to infected,
$\gamma > 0$ is the recovery rate of an infected individual,
$\alpha > 0$ is the infection fatality rate, and
$\rho > 0$ is the inverse of the average number of days for an infected person to die if they do not recover.
Our goal is to give a probabilistic estimate of the system parameters given noisy and/or sparse state observations.

To use the methodology of \Cref{sec:ode-methodology}, we recognize that \cref{eq:seird-equations} can be written as
\begin{align}
    \begin{aligned}
    \frac{\textup{d}}{\textup{d}t}\left[\begin{array}{c}
        \hat{q}_{S}(t) \\ \hat{q}_{E}(t) \\ \hat{q}_{I}(t) \\ \hat{q}_{R}(t) \\ \hat{q}_{D}(t)
    \end{array}\right]
    = \left[\begin{array}{cccc}
        -\hat{q}_{S}(t)\hat{q}_{I}(t) & 0 & 0 & 0 \\
        \hat{q}_{S}(t)\hat{q}_{I}(t) & -\hat{q}_{E}(t) & 0 & 0 \\
        0 & \hat{q}_{E}(t) & -\hat{q}_{I}(t) & -\hat{q}_{I}(t) \\
        0 & 0 & \hat{q}_{I}(t) & 0 \\
        0 & 0 & 0 & \hat{q}_{I}(t)
    \end{array}\right]
    \left[\begin{array}{c}
        \beta \\ \delta \\ (1 - \alpha)\gamma \\ \alpha \rho
    \end{array}\right],
    \end{aligned}
\end{align}
which is \cref{eq:ode-setup} with parameter vector $\hat{\vb{o}} = (\beta,\delta,(1-\alpha)\gamma,\alpha\rho)\trp\in\mathbb{R}^{4}$.
Note that an estimate of $\hat{\vb{o}}$ does not explicitly provide estimates for $\alpha$, $\gamma$, or $\rho$, but it does enable prediction in time and for new initial conditions.

\subsection{Data generation}
To generate data for this problem, we set initial conditions
$
\hat{\vb{q}}(0) = (0.994, 0.005, 0.001, 0, 0)\trp
$
and true parameter values $\beta = 0.25$, $\delta = 0.1$, $\gamma = 0.1$, $\alpha = 0.05$, and $\rho = 0.05$. Hence, the true parameter vector is $\hat{\vb{o}} = (0.25, 0.1, 0.095, 0.0025)\trp$. The system \cref{eq:seird-equations} is then solved with the explicit adaptive-step Runge--Kutta time stepping method of order 5 over the time domain $[0, 199]$ days.
We conduct experiments with different numbers of observations $m$, final observation times $t_{m-1}$, and noise levels $\xi > 0$.
The state variables are observed at integer sample times chosen uniformly from $\mathbb{Z}\cap [t_{0},t_{m-1}]$ (i.e., sample times are rounded to the nearest day). Thus, if $m = t_{m-1} + 1$, all state variables are observed daily; when $m < t_{m-1}+1$, each state variable can be observed over a potentially different set of $m$ days within time $t_{m-1}$.
Letting $y_{i}(t_j)\in\mathbb{R}$ denote the true solution of the $i$-th state variable at sample time $t_j$, for each $(i,j)$ pair we record noisy observations drawn from a truncated normal distribution centered at $y_{i}(t_j)$ with magnitude-dependent variance $\xi y_{i}(t_j)$ and truncation limits $[0,1]$.
Hence, the observations are contained in the interval $[0,1]$, with more noise applied to states of larger magnitude.
We thereby obtain snapshot vectors $\vb{y}_{i} = (y_{i,0},y_{i,1},\ldots,y_{i,m-1})\trp$ for each state variable, which are used as inputs to \Cref{alg:GP-BayesODEs} (the analog of \Cref{alg:GP-BayesOpInf} for this setting).

\subsection{Parameter estimation results}
\begin{figure}[t]
    \centering
    \includegraphics[width=\textwidth]{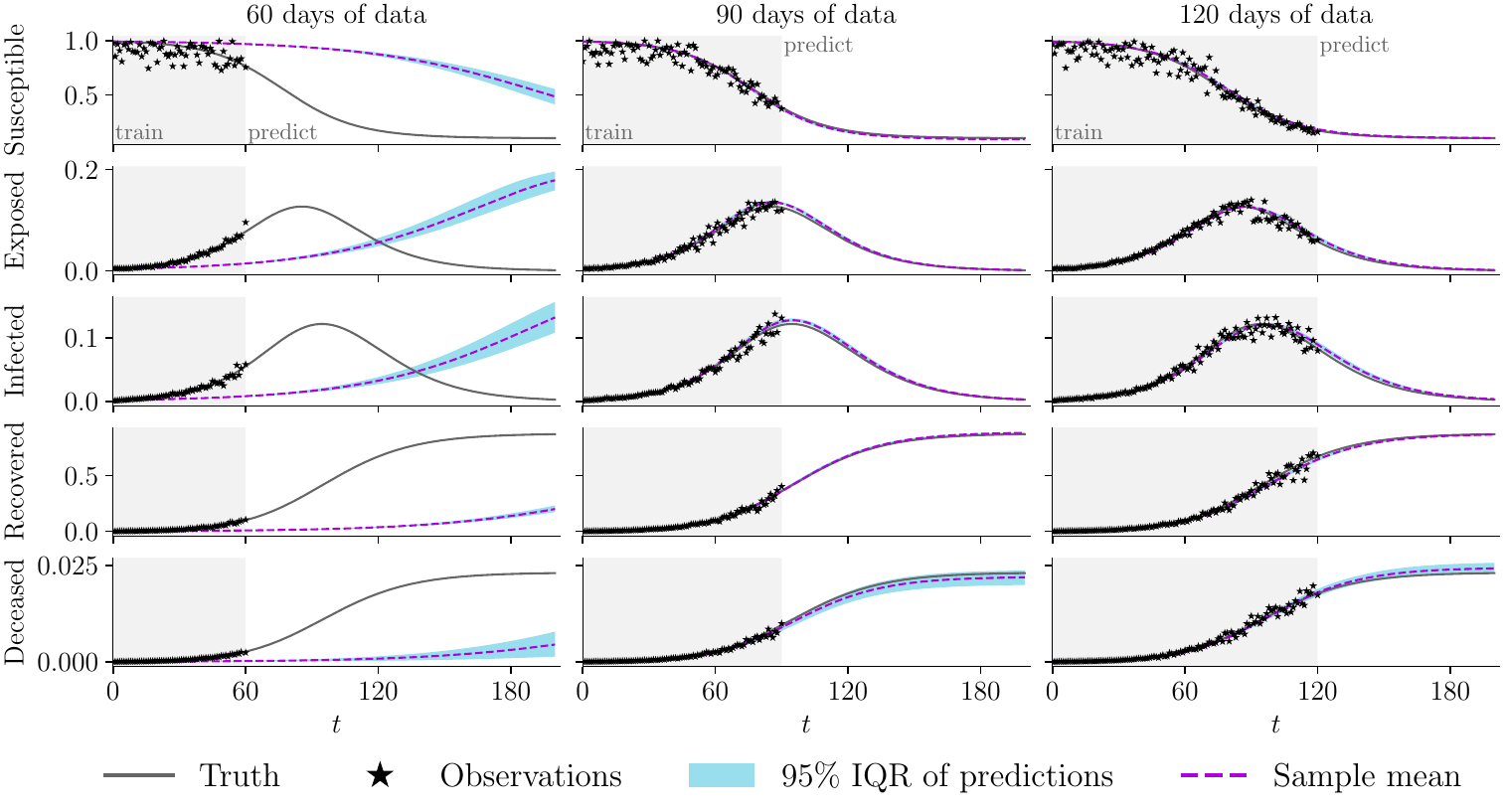}
    \vspace{-0.75cm}
    \caption{Bayesian ODE predictions for the SEIRD system \cref{eq:seird-equations} with varying amounts of data and $\xi = 10\%$ relative noise.}
    \label{fig:seird-noisy}
\end{figure}

We first set $t_{m-1} = m-1$ and relatively high noise level $\xi = 10\%$ (i.e., one noisy observation per day). We then apply \Cref{alg:GP-BayesODEs} with $m = 60$, $90$, and $120$ days of data, using $m' = 4m$ GP~estimates in each case. Sampling the posterior distributions for $\hat{\vb{o}}$ $500$~times and solving the corresponding SEIRD system produces the results shown in \Cref{fig:seird-noisy}. Using data for only 60 days leads to poor results and the proposed method is not able to learn the dynamics of the underlying system due to the high level of noise. However, as the number of days for which data are available increases to $90$ and $120$, the predictions in the test regime become increasingly accurate. \Cref{table:seird-results} reports the mean of the posterior distribution for each case and shows a similar trend. This shows training data through a certain amount of time $t_{m-1}$ is required to capture the dynamics and accurately extrapolate in time.

\begin{table}[t]
    \centering
    \begin{tabular}{lll|c|c|c|c}
        &&& $\beta$ & $\delta$ & $(1-\alpha)\gamma$ & $\alpha\rho$ \\
        \hline
        \Cref{fig:seird-noisy} && \textbf{Truth:} & \textbf{0.250} & \textbf{0.100} & \textbf{0.095} & \textbf{0.0025} \\
$m = 60$ & $t_{m-1}=m-1$ & $\xi=10\%$ &
        0.097 & 0.029 & 0.0269 & 0.00061 \\
        $m = 90$ & $t_{m-1}=m-1$ & $\xi=10\%$ &
        0.252 & 0.096 & 0.0922 & 0.00228 \\
        $m = 120$ & $t_{m-1}=m-1$ & $\xi=10\%$ &
        0.243 & 0.098 & 0.0926 & 0.00259 \\
        \hline
        \Cref{fig:seird-sparse} && \textbf{Truth:} & \textbf{0.250} & \textbf{0.100} & \textbf{0.095} & \textbf{0.0025} \\
$m = 10$ & $t_{m-1}+1=60$ & $\xi=5\%$ &
        0.123 & 0.080 & 0.0769 & 0.00162 \\
        $m = 10$ & $t_{m-1}+1=90$ & $\xi=5\%$ &
        0.222 & 0.098 & 0.0725 & 0.00253 \\
        $m = 10$ & $t_{m-1}+1=120$ & $\xi=5\%$ &
        0.245 & 0.099 & 0.0915 & 0.00244 \\
    \end{tabular}
    \caption{Parameter posterior means for the SEIRD system \cref{eq:seird-equations}.}
    \label{table:seird-results}
\end{table}

\begin{figure}[t]
    \centering
    \includegraphics[width=\textwidth]{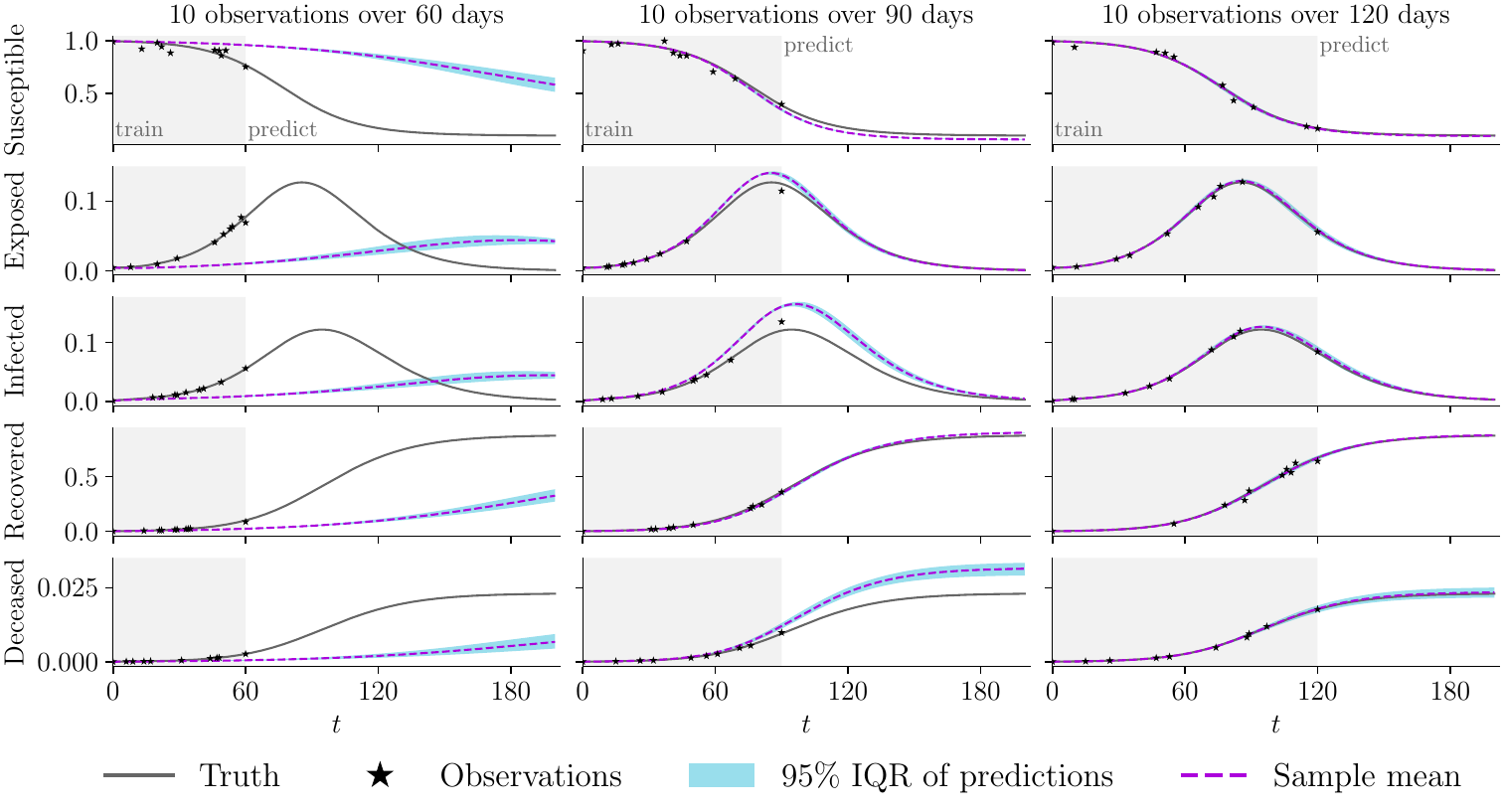}
    \vspace{-0.75cm}
    \caption{Bayesian ODE predictions for the SEIRD system \cref{eq:seird-equations} with $m = 10$ observations, $\xi = 5\%$ relative noise, and varying observation periods. Data observation times are different for each state.}
    \label{fig:seird-sparse}
\end{figure}

To further analyze the effect of $t_{m-1}$ on the proposed GP-based parameter estimation method, we limit the number of observations to $m = 10$ and reduce noise level to $\xi = 5\%$ while changing the time over which data are collected to $t_{m-1} + 1 = 60$, $90$, and $120$ days. Note that data observation times are not aligned across the states. In each case $m' = 4(t_{m-1}+1)$ GP~estimates are employed in \Cref{alg:GP-BayesODEs}. The means of the posterior parameter distribution are shown in \Cref{table:seird-results} and sampled predictions are provided in \Cref{fig:seird-sparse}. The results of the Bayesian fit are accurate with even $m = 10$ observations as long as they are spread over $120$ days, but are inaccurate when training data are only available for $60$ days. Hence, for this experiment, the performance of the method is more sensitive to the quality of the training data than to the amount of data.

\begin{remark}[Equation fitting vs.~trajectory fitting]
For Bayesian parameter estimation in ODEs, the likelihood setting in this work essentially fits given equations to data. Time derivatives on the left-hand side of the ODEs are estimated using GPs, and the parameters on the right-hand side are estimated through a linear regression problem. This is fundamentally different from a likelihood that compares the evolution of the dynamical model---time integration of the ODEs---with the observed data, i.e., fitting the trajectories. The latter may also improve the robustness to data noise and sparsity, but the training process often demands more computational effort. A typical example is the (Bayesian) neural ODE approach \cite{chen2018neural,dandekar2020bayesian}, which requires solving adjoint state equations to evaluate the gradients of loss functions.
\end{remark}

\section{Concluding remarks}
\label{sec:conclusion}

This work presents a new Bayesian method for OpInf that uses a probabilistic formulation for the reduced-order operators, thereby combining physics-based and data-driven modeling. In particular, GP emulation incorporates both reduced state data and their differential equation structure into the definition of the likelihood function for the reduced operators. The power of GPs in smoothing and interpolating noisy and/or sparse data enables probabilistic ROMs in scenarios where standard OpInf is not feasible.
The closed-form posterior estimates of the reduced operators provide an interpretable low-dimensional representation of complex physics-based models.

The numerical results demonstrate that GP-BayesOpInf is effective at constructing probabilistic surrogates for high-dimensional dynamical systems with polynomial structure given only data and knowledge of the model structure. The method was shown to perform well with abundant but noisy data, as well as in low-noise scenarios with sparse data. GP-BayesOpInf also applies to problems with multiple trajectories stemming from different initial conditions or input functions, as well as to a large class of differential systems with parameters shared across the equations.

There are several avenues for future work related to GP-BayesOpInf.
Our experiments used squared-exponential GP kernels and linear dimensionality reduction by POD, but the methodology is flexible and can be applied with any sufficiently differentiable GP kernel and with any type of compression (including nonlinear manifold approximators such as in \cite{barnett2022quadmanifold,geelen2024learning,geelen2023quadmanifold,jain2017quadratic,schwerdtner2024greedy}). It would be interesting to explore if embedding desirable properties in the kernel can improve model performance for problems with specific characteristics (for example, periodicity).
In the presence of large amounts of noise, the number of reliable reduced modes that can be obtained through POD is limited. To reduce the projection error by including more reduced modes in GP-BayesOpInf, it may be advantageous to leverage filtering strategies as shown in~\cite{farcas2022filtering} prior to applying dimensionality reduction. We also plan to explore more advanced GP techniques, such as sparse GPs \cite{smola2000sparsegp,snelson2005sparsegp} to reduce the trajectory data needed for regression (i.e., selecting $m'$ automatically), multifidelity GPs \cite{alvarez2012vectorkernels,chaudhuri2021mfegra,le2013multi,ohagan2000multigps,perdikaris2017nonlinear,poloczek2017multi} to fuse data from different sources, as well as a single vector-valued GP for all reduced modes to incorporate physically meaningful spatial correlations. Another promising direction is to leverage the quantified uncertainties to guide the active selection of initial conditions and/or input functions for parsimoniously generating informative training trajectories.

\section*{Acknowledgments}

S.M. was supported by the John von Neumann postdoctoral fellowship, a position sponsored by the Applied Mathematics Program of the U.S.~Department of Energy Office of Advanced Scientific Computing Research.
This article has been authored by an employee of National Technology \& Engineering Solutions of Sandia, LLC under Contract No.~DE-NA0003525 with the U.S.~Department of Energy (DOE). The employee owns all right, title and interest in and to the article and is solely responsible for its contents.
This paper describes objective technical results and analysis. Any subjective views or opinions that might be expressed in the paper do not necessarily represent the views of the U.S.~Department of Energy or the United States Government.

A.C.~and K.E.W.~acknowledge support from DARPA Automating Scientific Knowledge Extraction and Modeling (ASKEM) program under Contract No. HR0011262087 (award number 653002). The views, opinions, and/or findings expressed are those of the author(s) and should not be interpreted as representing the official views or policies of the Department of Defense or the U.S.~Government.

This work was in part carried out when M.G.~held a position at the University of Twente (NL), for which he acknowledges financial support from \emph{Sectorplan Bèta} under the focus area \emph{Mathematics of Computational Science}.

\newpage

\begin{appendices}

\section{Gaussian process regression}
\label{appendix:gps}

A Gaussian process is a collection of random variables, any finite number of which obeys a joint normal distribution. In Gaussian process regression, the prior on the regression function is assumed to follow a Gaussian process, i.e., for $(\vb{x},\vb{x}')\in \Omega \times \Omega$,
\begin{equation}
f(\vb{x})\sim \mathcal{GP}(b(\vb{x}),\kappa(\vb{x},\vb{x}'))\,,
\end{equation}
where $\Omega$ is the input domain of interest, $b: \Omega\to \mathbb{R}$ is a given mean function, and $\kappa:\Omega\times\Omega\to \mathbb{R}$ is a covariance function which often depends on kernel hyperparameters $\vb*{\theta}$. We also assume that the observation of $f(\vb{x})$, denoted by $y(\vb{x})$, may be corrupted by an independent white noise term, i.e., $y(\vb{x})=f(\vb{x})+\epsilon$ with $\epsilon\sim\mathcal{N}(0,\chi)$.

Given a finite number $m$ of training input locations $\vb{X}=[~\vb{x}_1~|~\vb{x}_2~|~\cdots~|~\vb{x}_m~]\trp$ in the domain $\Omega$, a joint Gaussian is thus defined for the corresponding observed outputs $\vb{y}=\{y(\vb{x}_1),y(\vb{x}_2),\cdots,y(\vb{x}_m)\}\trp$:
\begin{equation}
\vb{y}|\vb{X}~\sim~\mathcal{N}(b(\vb{X}),\vb{K})\,,\quad \vb{K}=\kappa(\vb{X},\vb{X})+\chi \vb{I}_m\,.
\end{equation}
Using the standard rule for conditioning Gaussians, the (predictive) posterior distribution for any new test input $\vb{s}\in\Omega$ conditioned on the training data can be formulated as a new Gaussian process $f^*(\vb{s})$:
\begin{align}
\begin{aligned}
    & f^*(\vb{s})|\vb{X},\vb{y}~\sim~\mathcal{GP}(b^*(\vb{s}),\kappa^*(\vb{s},\vb{s}'))\,, \quad &&\text{with}\\
    & b^*(\vb{s}) = b(\vb{s})+\kappa(\vb{s},\vb{X})\vb{K}^{-1}(\vb{y}-b(\vb{X}))\,,\quad &&\text{and}\qquad \kappa^*(\vb{s},\vb{s}') = \kappa(\vb{s},\vb{s}')-\kappa(\vb{s},\vb{X})\vb{K}^{-1}\kappa(\vb{X},\vb{s}')\,.
\end{aligned}
\end{align}
The values of hyperparameters can be tuned to optimize the predictive performance. A widely used technique is the empirical Bayes approach of maximizing marginal likelihood. With a standard gradient-based optimizer, the optimal hyperparameters $\vb*{\theta}^*$ can be estimated through the following maximization problem:
\begin{equation}
\begin{split}
\vb*{\theta}^* & =\arg\max_{\vb*{\theta}}~\log p(\vb{y}|\vb{X},\vb*{\theta}) \\
& =\arg\max_{\vb*{\theta}}\left\{-\frac{1}{2}(\vb{y}-b(\vb{X}))^{\mathrm{T}}\vb{K}^{-1}(\vb*{\theta})
(\vb{y}-b (\vb{X}))-\frac{1}{2}\log\abs{\vb{K}(\vb*{\theta})}-\frac{m}{2}\log(2\pi)\right\}\,.
\end{split}
\end{equation}
Here, $p(\vb{y}|\vb{X},\vb*{\theta})$ is the marginal likelihood function of $\vb*{\theta}$, which is given by the density function value of the observations $\vb{y}$ over the input locations $\vb{X}$ given the hyperparameters $\vb*{\theta}$.

\section{Algorithmic extensions}
\label{appendix:multiple-trajectories}

Algorithms \ref{alg:OperatorPosterior}--\ref{alg:GP-BayesOpInf} are written for problems in the PDE setting where a single trajectory of data is available for training.
This appendix generalizes the algorithms to the special cases discussed in \Cref{sec:multitrajectory} and \Cref{sec:odeparamestimate}.

\subsection{Multiple trajectories}

If data are available from $\ell > 1$ trajectories due to varying initial conditions or inputs, as discussed in \Cref{sec:multitrajectory}, the algorithmic procedure requires minor modifications.
First, one GP is trained per reduced mode per trajectory, resulting in $r\cdot \ell$ total executions of \Cref{alg:GP-Fit}.
Second, for each reduced mode, GP estimates from all $\ell$ trajectories are concatenated as described in \cref{eq:multitrajectory-concatenation}.
Third, the prior variances are selected by minimizing the average reconstruction error across all $\ell$ trajectories.
\Cref{alg:OpInfErrorMulti} and \Cref{alg:GP-BayesOpInfMulti} adapt \Cref{alg:OpInfError} and \Cref{alg:GP-BayesOpInf}, respectively, to this setting.
\Cref{alg:GP-BayesOpInfMulti} uses $\vb*{t}^{(l)}\in\mathbb{R}^{m^{(l)}}$ to denote the $m^{(l)}\in\mathbb{N}$ observation times for the $l$-th trajectory, which may be different for each $l$.

\begin{algorithm}[t]
\begin{algorithmic}[1]
\Procedure{OpInfErrorMulti}{
    \newline\phantom{---}Prior variances $\vb*{\gamma}_1,\ldots,\vb*{\gamma}_r \in \mathbb{R}^{d}$;
    \newline\phantom{---}Regression arguments $\vb*{\alpha} = \{\tilde{\vb{D}},\tilde{\vb{z}}_{1},\ldots,\tilde{\vb{z}}_{r},\vb{W}_{1}^{1/2},\ldots,\vb{W}_{r}^{1/2}\}$,
\newline\phantom{---}GP state estimates $\tilde{\vb{Q}}^{(1)},\ldots,\tilde{\vb{Q}}^{(\ell)}\in\mathbb{R}^{r\times m'}$
    \newline\phantom{---}Input functions $\vb{u}^{(1)},\ldots,\vb{u}^{(\ell)}:\mathbb{R}\to\mathbb{R}^{p}$,
    \newline\phantom{---}GP estimation times $\vb*{t}' \in\mathbb{R}^{m'}$,
    \newline\phantom{---}Stability margin parameter $\varphi > 0$,
    \newline\phantom{---}Number of samples $n_s \in \mathbb{N}$,
    \newline\phantom{---}Final prediction time $t_f > t'_0$
    \newline}

\State $(\vb*{\mu}_1,\vb*{\Sigma}_1),\ldots,(\vb*{\mu}_r,\vb*{\Sigma}_r) \gets
        \textproc{OpPost(}
            \vb*{\alpha},
\vb*{\gamma}_1,\ldots,\vb*{\gamma}_r
        \textproc{)}$

    \vspace{.1cm}
    \LineComment{Check stability and compute the training error for each trajectory.}
    \For{$l = 1,\ldots,\ell$}
        \State $\hat{\vb{q}}_0^{(l)} \gets \tilde{\vb{Q}}_{:,0}^{(l)}$
\State $B^{(l)} \gets \varphi \cdot \max_{ij}|\tilde{\vb{Q}}^{(l)}_{ij}|$
\For{$k = 1, \ldots, n_s$}
            \State $\hat{\vb{O}}_{k}^{(l)} \gets~\text{sample from}~~[~\mathcal{N}(\vb*{\mu}_1,\vb*{\Sigma}_1)~~\cdots~~\mathcal{N}(\vb*{\mu}_r,\vb*{\Sigma}_r)~]\trp$
\State $\breve{\vb{Q}} \gets~\text{solve}~~\frac{\textrm{d}}{\textrm{d}t}\hat{\vb{q}}(t) = \hat{\vb{O}}_{k}^{(l)}\vb{d}(\hat{\vb{q}}(t),\vb{u}^{(l)}(t)),~\hat{\vb{q}}(t'_0) = \hat{\vb{q}}_0^{(l)}~~\text{over}~~t\in[t'_0,t_f]$
            \If{$\max_{ij}|\breve{\vb{Q}}_{ij}| > B$}
\State \textbf{return} $\infty$
            \EndIf
            \State $\hat{\vb{Q}}_{k}^{(l)} \gets~\text{solve}~~\frac{\textrm{d}}{\textrm{d}t}\hat{\vb{q}}(t) = \hat{\vb{O}}_{k}^{(l)}\vb{d}(\hat{\vb{q}}(t),\vb{u}^{(l)}(t)),~\hat{\vb{q}}(t'_0) = \hat{\vb{q}}_0^{(l)}~~\text{over}~~t\in\vb*{t}'$
        \EndFor
        \State $\hat{\vb{Q}}^{(l)} = \frac{1}{n_s}\sum_{k=1}^{n_s}\hat{\vb{Q}}_{k}^{(l)}$
\EndFor
    \State \textbf{return} $\frac{1}{n_s}\sum_{l=1}^{\ell}\|\tilde{\vb{Q}}^{(l)} - \hat{\vb{Q}}^{(l)}\|$
        \Comment{Average error compared to GP state estimates.}
\EndProcedure
\end{algorithmic}
\caption{Estimate the mean error of the Bayesian OpInf model over multiple trajectories.}
\label{alg:OpInfErrorMulti}
\end{algorithm}

\begin{algorithm}[t]
\begin{algorithmic}[1]
\Procedure{GPBayesOpInfMulti}{
\newline\phantom{---}Observation times $\vb*{t}^{(l)}\in\mathbb{R}^{m^{(l)}}$, $l=1,\ldots,\ell$,
    \newline\phantom{---}Snapshots $\vb{Q}^{(l)}\in\mathbb{R}^{N\times m^{(l)}}$, $l=1,\ldots,\ell$,
    \newline\phantom{---}Reduced dimension $r \ll N$,
\newline\phantom{---}Estimation times $\vb*{t}'\in\mathbb{R}^{m'}$,
    \newline\phantom{---}Numerical regularization $0 < \tau \ll 1$,
    \newline\phantom{---}Input functions $\vb{u}^{(1)},\ldots,\vb{u}^{(\ell)}:\mathbb{R}\to\mathbb{R}^{p}$
    \newline
    }
    \LineComment{GP regression on compressed training data.}
    \State $\vb{V} \gets \texttt{basis}(\vb{Q}^{(1)},\ldots,\vb{Q}^{(\ell)}; r)$
        \Comment{Global rank-$r$ basis (e.g., POD).}
    \For{$l=1,\ldots,\ell$}
        \State $[~\vb{y}_1^{(l)}~~\cdots~~\vb{y}_r^{(l)}~]\trp \gets \vb{V}\trp\vb{Q}^{(l)}$
        \For{$i = 1, \ldots, r$}
            \State $\tilde{\vb{y}}_i^{(l)},\tilde{\vb{z}}_i^{(l)},(\vb{W}^{zz\,(l)}_i)^{1/2} \gets \textproc{GPFit(}\vb*{t}^{(l)},\vb{y}_i^{(l)},\vb*{t}',\tau\textproc{)}$
        \EndFor
        \LineComment{Prepare regression data for the $l$-th trajectory.}
        \State $\tilde{\vb{Q}}^{(l)} \gets [~\tilde{\vb{y}}_1^{(l)}~~\cdots~~\tilde{\vb{y}}_r^{(l)}~]\trp$
        \State $\tilde{\vb{D}}^{(l)} \gets \texttt{data\_matrix}(\tilde{\vb{Q}}^{(l)},\vb{u}^{(l)}(\vb*{t}'))$
    \EndFor

    \LineComment{Concatenate regression data across all trajectories.}
    \State $\tilde{\vb{D}} \gets [~(\tilde{\vb{D}}^{(1)})\trp~~\cdots~~(\tilde{\vb{D}}^{(\ell)})\trp~]\trp$
        \Comment{Concatenate data matrices.}
    \For{$i=1,\ldots,r$}
        \State $\tilde{\vb{z}}_{i} \gets [~(\tilde{\vb{z}}_{i}^{(1)})\trp~~\cdots~~(\tilde{\vb{z}}_{i}^{(\ell)})\trp~]\trp$
            \Comment{Concatenate time derivatives.}
        \State $\vb{W}_{i} \gets \operatorname{diag}(\vb{W}_{i}^{zz\,(1)},\cdots,\vb{W}_{i}^{zz\,(\ell)})$
            \Comment{Block diagonal weights.}
    \EndFor
    \State $\vb*{\alpha} \gets \{\tilde{\vb{D}},\tilde{\vb{z}}_{1},\ldots,\tilde{\vb{z}}_{r},\vb{W}_{1}^{1/2},\ldots,\vb{W}_{r}^{1/2}\}$
        \Comment{Regression arguments.}

    \LineComment{Operator inference regression to construct a probabilistic reduced model.}
    \State $\vb*{\gamma}_{1}^{*},\ldots,\vb*{\gamma}_{r}^{*} \,\gets\,
\operatorname{argmin}_{\vb*{\gamma}_{1},\ldots,\vb*{\gamma}_{r}}
    \textproc{OpInfErrorMulti}(\vb*{\gamma}_{1},\ldots,\vb*{\gamma}_{r}; \vb*{\alpha},\tilde{\vb{Q}}^{(1)},\ldots,\tilde{\vb{Q}}^{(\ell)},\vb{u}^{(1)},\ldots,\vb{u}^{(\ell)},\vb*{t}')$
    \State $(\vb*{\mu}_1,\vb*{\Sigma}_1),\ldots,(\vb*{\mu}_r,\vb*{\Sigma}_r) \,\gets\, \textproc{OpPost(}\vb*{\alpha},\vb*{\gamma}_1^{*},\ldots,\vb*{\gamma}_r^{*}\textproc{)}$
    \State \textbf{return} $[~\mathcal{N}(\vb*{\mu}_1,\vb*{\Sigma}_1)~~\cdots~~\mathcal{N}(\vb*{\mu}_r,\vb*{\Sigma}_r)~]\trp$
\EndProcedure
\end{algorithmic}
\caption{Bayesian Operator Inference with Gaussian processes for multiple data trajectories.}
\label{alg:GP-BayesOpInfMulti}
\end{algorithm}

\subsection{Parameter estimation for ordinary differential equations}

The parameter estimation setting of \Cref{sec:odeparamestimate} requires small algorithmic adjustments and simplifications.
\Cref{alg:OpInfErrorODEs} and \Cref{alg:GP-BayesODEs} adapt \Cref{alg:OpInfError} and \Cref{alg:GP-BayesOpInf}, respectively, to this scenario.

\begin{algorithm}[t]
\begin{algorithmic}[1]
\Procedure{OpInfErrorODEs}{
    \newline\phantom{---}Prior variance $\vb*{\gamma} \in \mathbb{R}^{d}$;
    \newline\phantom{---}Regression arguments $\vb*{\alpha} = \{\tilde{\vb{D}},\tilde{\vb{z}},\vb{W}^{1/2}\}$,
    \newline\phantom{---}GP state estimates $\tilde{\vb{Q}}\in\mathbb{R}^{r\times m'}$
    \newline\phantom{---}Input functions $\vb{u}:\mathbb{R}\to\mathbb{R}^{p}$,
    \newline\phantom{---}GP estimation times $\vb*{t}' \in\mathbb{R}^{m'}$,
    \newline\phantom{---}Stability margin parameter $\varphi > 0$,
    \newline\phantom{---}Number of samples $n_s \in \mathbb{N}$,
    \newline\phantom{---}Final prediction time $t_f > t'_0$
    \newline}

    \State $\vb*{\mu},\vb*{\Sigma} \gets
        \textproc{OpPost(}
            \vb*{\alpha}, \vb*{\gamma}
        \textproc{)}$
        \Comment{\Cref{alg:OperatorPosterior} with $r=1$.}

    \State $\hat{\vb{q}}_0 \gets \tilde{\vb{Q}}_{:,0}$
    \State $B \gets \varphi \cdot \max_{ij}|\tilde{\vb{Q}}_{ij}|$
    \For{$k = 1, \ldots, n_s$}
        \State $\hat{\vb{o}}_{k} \gets~\text{sample from}~\mathcal{N}(\vb*{\mu},\vb*{\Sigma})$
        \State $\breve{\vb{Q}} \gets~\text{solve}~~\frac{\textrm{d}}{\textrm{d}t}\hat{\vb{q}}(t) = \breve{\vb{D}}(\hat{\vb{q}}(t),\vb{u}(t))\hat{\vb{o}}_{k},~\hat{\vb{q}}(t'_0) = \hat{\vb{q}}_0~~\text{over}~~t\in[t'_0,t_f]$
        \If{$\max_{ij}|\breve{\vb{Q}}_{ij}| > B$}
            \State \textbf{return} $\infty$
        \EndIf
        \State $\hat{\vb{Q}}_{k} \gets~\text{solve}~~\frac{\textrm{d}}{\textrm{d}t}\hat{\vb{q}}(t) = \breve{\vb{D}}(\hat{\vb{q}}(t),\vb{u}(t))\hat{\vb{o}}_{k},~\hat{\vb{q}}(t'_0) = \hat{\vb{q}}_0~~\text{over}~~t\in\vb*{t}'$
    \EndFor
    \State $\hat{\vb{Q}} = \frac{1}{n_s}\sum_{k=1}^{n_s}\hat{\vb{Q}}_{k}$
    \State \textbf{return} $\|\tilde{\vb{Q}} - \hat{\vb{Q}}\|$
\EndProcedure
\end{algorithmic}
\caption{Estimate the mean error of the Bayesian OpInf model for ODE parameter estimation.}
\label{alg:OpInfErrorODEs}
\end{algorithm}

\begin{algorithm}[!ht]
\begin{algorithmic}[1]
\Procedure{GPBayesODEs}{
\newline\phantom{---}Observation times $\vb*{t}_{i}\in\mathbb{R}^{m_{i}}$, $i=1,\ldots,r$,
    \newline\phantom{---}Snapshots $\vb{y}_{i}\in\mathbb{R}^{m_i}$, $i=1,\ldots,r$,
\newline\phantom{---}Estimation times $\vb*{t}'\in\mathbb{R}^{m'}$,
    \newline\phantom{---}Numerical regularization $0 < \tau \ll 1$,
    \newline\phantom{---}Input functions $\vb{u}:\mathbb{R}\to\mathbb{R}^{p}$
    \newline
    }

\LineComment{GP regression on state observations.}
    \For{$i = 1, \ldots, r$}
        \State $\tilde{\vb{y}}_i,\tilde{\vb{z}}_i,(\vb{W}^{zz}_i)^{1/2} \gets \textproc{GPFit(}\vb*{t},\vb{y}_i,\vb*{t}',\tau\textproc{)}$
\EndFor
    \LineComment{Prepare regression data.}
    \State $\tilde{\vb{Q}} \gets [~\tilde{\vb{y}}_1~~\cdots~~\tilde{\vb{y}}_r~]\trp$
    \State $\tilde{\vb{D}} \gets \texttt{data\_matrix}(\tilde{\vb{Q}},\vb{u}(\vb*{t}'))$
\State $\tilde{\vb{z}} \gets [~\tilde{\vb{z}}_{1}\trp~~\cdots~~\tilde{\vb{z}}_{r}\trp~]\trp$
        \Comment{Concatenate time derivatives.}
    \State $\vb{W} \gets \operatorname{diag}\left(\vb{W}_{1}^{zz},\cdots,\vb{W}_{r}^{zz}\right)$
        \Comment{Block diagonal weights.}
    \State $\vb*{\alpha} \gets \{\tilde{\vb{D}},\tilde{\vb{z}},\vb{W}^{1/2}\}$
        \Comment{Regression arguments.}

    \LineComment{Operator inference regression to construct a probabilistic reduced model.}
    \State $\vb*{\gamma}^{*}\,\gets\,
\operatorname{argmin}_{\vb*{\gamma}}
    \textproc{OpInfErrorODEs}(\vb*{\gamma}; \vb*{\alpha},\tilde{\vb{Q}},\vb{u},\vb*{t}')$
        \Comment{\Cref{alg:OpInfErrorODEs}.}
    \State $\vb*{\mu},\vb*{\Sigma} \,\gets\, \textproc{OpPost(}\vb*{\alpha},\vb*{\gamma}^{*}\textproc{)}$
        \Comment{\Cref{alg:OperatorPosterior} with $r=1$.}
    \State \textbf{return} $\mathcal{N}(\vb*{\mu},\vb*{\Sigma})$

\EndProcedure
\end{algorithmic}
\caption{Bayesian parameter estimation for ODEs with Gaussian processes.}
\label{alg:GP-BayesODEs}
\end{algorithm}

\clearpage
\section{Main Nomenclature}
\label{appendix:nomenclature}

\textcolor{white}{Notation}
\begin{tabular}{rl}
    \textbf{Full order state} \\
    $N\in\mathbb{N}$ & full-order state dimension \\
    $\vb{q}(t) \in \mathbb{R}^{N}$ & full-order state vector \\
    $\vb{u}(t) \in \mathbb{R}^{p}$ & input function (time-dependent forcing / boundary conditions, etc.) \\[.0625in]

    \textbf{Reduced order state} \\
    $r\in\mathbb{N}$ & reduced-order state dimension, number of basis vectors \\
    $\hat{\vb{q}}(t) \in \mathbb{R}^{r}$ & reduced-order state vector with entries $(\hat{q}_{1}(t),\ldots,\hat{q}_{r}(t))$ \\
    $\vb{V}\in\mathbb{R}^{N\times r}$ & rank-$r$ basis, i.e., $\vb{q}(t) \approx \vb{V}\hat{\vb{q}}(t) + \bar{\vb{q}}$ \\
    $\vb{d}:\mathbb{R}^{r}\times\mathbb{R}^{p}\to\mathbb{R}^{d}$ & data vector defining reduced-order model structure \\
    $\hat{\vb{O}}\in\mathbb{R}^{r \times d}$ & operator matrix defining reduced-order model structure \\
    $\hat{\vb{o}}_{i}\in\mathbb{R}^{d}$ & operators for the $i$-th reduced mode, the $i$-th row of $\hat{\vb{O}}$ \\[.0625in]

    \textbf{Training data} \\
    $m\in\mathbb{N}$ & number of available snapshots \\
    $\vb*{t}\in\mathbb{R}^{m}$ & times $t_0,\ldots,t_{m-1}$ corresponding to the observed snapshots \\
    $\vb{q}_{j}\in\mathbb{R}^{N}$ & possibly noisy full-order snapshot at time $t = t_{j}$, i.e., $\vb{q}_{j} \approx \vb{q}(t_j)$ \\
    $\bar{\vb{q}} \in \mathbb{R}^{N}$ & average training snapshot\\
    $\hat{\vb{q}}_{j}\in\mathbb{R}^{r}$ & reduced-order snapshot at time $t = t_{j}$, i.e., $\hat{\vb{q}}_{j} = \vb{V}\trp(\vb{q}_j - \bar{\vb{q}})$ \\
    $\vb{y}_{i}\in\mathbb{R}^{m}$ & reduced-order snapshot data for the $i$-th reduced mode,
        \\ & i.e., the $i$-th entry of each $\hat{\vb{q}}_{j}$, $j=0,\ldots,m-1$
\\[.0625in]

    \textbf{Gaussian process} \\
    $\kappa_i:\mathbb{R}\times\mathbb{R}\to\mathbb{R}$ & kernel function for the $i$-th reduced mode with hyperparameters $\vb*{\theta}_i$ \\
    $m'\in\mathbb{N}$ & number of GP reduced state / time derivative estimates \\
    $\vb*{t}'\in\mathbb{R}^{m'}$ & times corresponding to the GP estimates \\
    $\tilde{\vb{y}}_{i} \in \mathbb{R}^{m'}$ & GP reduced state estimate for the $i$-th reduced mode \\
$\tilde{\vb{z}}_{i} \in \mathbb{R}^{m'}$ & GP reduced time derivative estimate for the $i$-th reduced mode \\
$\vb{K}_{i}^{yy} \in \mathbb{R}^{m \times m}$ & kernel matrix for $\vb{y}$-$\vb{y}$ correlation (observation times $\vb*{t}$) \\
    $\vb{K}_{i}^{zz} \in \mathbb{R}^{m' \times m'}$ & kernel matrix for $\vb{z}$-$\vb{z}$ correlation (estimation times $\vb*{t}'$)\\
    $\vb{K}_{i}^{zy} \in \mathbb{R}^{m' \times m}$ & kernel matrix for $\vb{z}$-$\vb{y}$ correlation ($\vb*{t}'$ vs $\vb*{t}$) \\
    $\vb{W}_{i}^{zz} \in \mathbb{R}^{m' \times m'}$ & block of the inverse of the kernel matrix for observation times \\[.0625in]

    \textbf{Operator inference} \\
    $\vb{D} \in \mathbb{R}^{m\times d}$ & data matrix for standard OpInf \\
    $\tilde{\vb{D}} \in \mathbb{R}^{m' \times d}$ & GP-enhanced data matrix \\
    $\vb*{\mu}_{i} \in \mathbb{R}^{d}$ & posterior mean of $\hat{\vb{o}}_{i}$ \\
    $\vb*{\Sigma}_{i}\in\mathbb{R}^{d\times d}$ & posterior covariance matrix of $\hat{\vb{o}}_{i}$ \\
    $\vb*{\Gamma}_{i}\in\mathbb{R}^{d\times d}$ & prior covariance for $\hat{\vb{o}}_i$\\
    $\vb*{\gamma}_i\in\mathbb{R}^{d}$ & diagonals of the prior covariance $\vb*{\Gamma}_i$\\
\end{tabular}

\end{appendices}

\newpage
\bibliographystyle{abbrv}
\bibliography{references}

\end{document}